\documentclass[10pt]{article}

\usepackage{epsfig,subfigure}
\usepackage{amsmath,amssymb}

\usepackage{exscale}
\usepackage[mathscr]{eucal}
\usepackage{fancyhdr}
\usepackage{graphicx}
\usepackage{multirow,bigstrut}
\usepackage{labelfig}

\usepackage[a4paper]{geometry}
\newtheorem{Lemme}{Lemma}[section] 
\newtheorem{theorem}[Lemme]{Theorem}
\newtheorem{proposition}[Lemme]{Proposition}
\newtheorem{lemma}[Lemme]{Lemma}   
\newtheorem{corollary}[Lemme]{Corollary}
\newtheorem{remark}[Lemme]{Remark}	
\newtheorem{definition}[Lemme]{Definition}

\def\Box{\leavevmode\vrule height 5pt width 4pt depth 0pt\relax}
\date{\today}
\date{}

\begin{document}

\title{On generalized Newtonian
 fluids in curved pipes\thanks{This work was partially
supported by FCT project PEst-OE/MAT/UI4032/2011}}

\author{Nadir Arada\thanks{Departamento de Matem\'atica, Faculdade de Ci\^encias e Tecnologia da Universidade Nova de Lisboa, Quinta da Torre, 2829-516 Caparica, Portugal ({\tt naar@fct.unl.pt}).}}

\maketitle

\noindent \begin{abstract} \noindent This paper is concerned with steady, fully developed motion of a Navier-Stokes fluid with shear-dependent viscosity in a curved pipe under a given axial pressure gradient.  We establish existence and uniqueness results, derive appropriate estimates and prove a characterization of the secondary flows. The approximation, with respect to the curvature ratio, of the full governing systems by some Dean-like equation is studied. \vspace{3mm}\vspace{4mm}\vspace{4mm}
\end{abstract} 

\section{\large{Introduction}}
\label{introduction}\setcounter{equation}{0}
There is a great interest in the study of curved pipe flows due to its wide range of applications in engineering (e.g. hydraulic pipe systems related to corrosion failure) and in biofluid dynamics, such as blood flow in the vascular system. It is known since the pioneer experimental works of  Eustice (\cite{eustice1}, \cite{eustice2}) that these flows are challenging and much more complex than flows in straight pipes. Among their distinguishing features is the existence of 
secondary flows induced by the centrifugal force and which appear even for the slightest curvature.\vspace{2mm}\\
Fully developed viscous flow in a curved pipe with circular cross-section was first studied theoretically by Dean (\cite{dean1}, \cite{dean2}) in the case of Newtonian fluids by applying regular perturbation methods, the perturbation parameter being the curvature ratio. He simplified the governing equations, by neglecting all the effects due to pipe curvature except the centrifugal forces, 
and showed that for small curvature ratio the flow depends only on a single parameter, the so-called {\it Dean number}.
Following this fundamental work, the results based on perturbation solutions have been extended for a larger range of curvature ratio and Reynolds number, showing in particular the existence of additional pairs of vortices and multiple solutions (see e.g.  
\cite{yang}, \cite{daska}). Different geometries including circular, elliptical and annular cross-sections have also been considered  by several authors (see e.g.  \cite{berger}, \cite{ito}, \cite{Rob}, \cite{soh}, \cite{topa}). \vspace{2mm}\\
Despite the great interest in curved pipes, a rigorous analysis of the solvability of the Dean's equations and the full Navier-Stokes equations was not available prior to the work of Galdi and Robertson \cite{galdirobertson}, where existence and uniqueness results for steady, fully developed flows  are established. Existence of secondary motions is also studied, as well as the approximation of the Navier-Stokes equations
 by the Dean's equation for small curvature ratios. \vspace{1mm}\\
Flows of non-Newtonian fluids in curved pipes have also been studied by several authors  (see e.g.  \cite{Rob-coscia}, \cite{fan}, \cite{Rob-jichote}, \cite{Rob-mul}). Perturbation methods were used by Robertson and Muller \cite{Rob-mul} to study steady, fully developed flow of Oldroyd-B fluids, and to compare 
the results for creeping and non-creeping flows. For a second order model, Jitchote and Robertson \cite{Rob-jichote} obtained analytical solutions to the perturbation equations and analyze the effects of non-zero second normal stress coefficient on the behaviour of the solution. Theoretical results regarding 
this problem were obtained by Coscia and Robertson in \cite{Rob-coscia}, where existence and uniqueness of a strong solution for small non-dimensional pressure drop is established.\vspace{2mm}\\
The aim of the present paper is to extend the analysis carried-out in \cite{galdirobertson} to the class of quasi-Newtonian fluids. 
This class is described by partial differential equations of the quasilinear type (Navier-Stokes equations with a non constant viscosity that decreases with increasing shear rate in the case of shear-thinning flows and increases with increasing shear rate in the case of shear-thickening flows). It was first proposed and studied in bounded domains by Ladyzhenskaya in \cite{lady1}, \cite{lady2} and \cite{lady} as a modification of the Navier-Stokes system, and was similarly suggested by Lions in \cite{lions}.  Existence of weak solutions was proved by both authors using compactness arguments and the theory of monotone operators. Much work has been done since these pioneering results and, without ambition for completeness, we cite  Ne\v cas {\it et al.}  who established existence of weak solutions under less restrictive assumptions (see for example \cite{necas}). \vspace{2mm}\\
Since we are dealing with fully developed flows in curved pipes, the typical issues related to the nonlinear extra stress tensor and the convective term arise and can be handled as in the case of bounded domains. However, additional difficulties related with extra terms (depending on the curvature ratio) occur and need to be managed. The splitting method, consisting in two {\it coupled} formulations respectively associated to the secondary flows and to the axial flow and used in \cite{galdirobertson} for the study of the Newtonian case, cannot be applied. Because of the nonlinearity of the shear stress tensor, the coercivity property of the corresponding bilinear forms in the shear-thinning case, and the monotonicity property of these forms in both shear-thinning and shear thickening cases fail to be satisfied. To overcome these difficulties, we consider a {\it global} formulation in an appropriate functional setting and adapt some standard tools, such as the Korn inequality. An existence result is established for arbitrary values of the Reynolds number, of the pressure drop and for any curvature ratio of the pipe, and a uniqueness result for small Reynolds numbers. The {\it global} formulation allows also to derive uniform estimates independent of the Dean number. Using a posteriori the splitting method, we establish other estimates for the secondary flows that highlight the connection with the Dean number. Following \cite{galdirobertson}, we also prove that there are no secondary flows if, and only if, the Dean number is equal to zero. Finally, we consider an approximation problem (that can be seen as a generalization of the Dean's equation), study its solvability, establish corresponding estimates and evaluate the relationship between its solutions and those of the full governing equations.\vspace{2mm}\\
The plan of the paper is as follows. The governing equations are given in Section 2. Notation, and some preliminary results are given in Section 3. Section 4 is devoted to the statement and discussion of the main results. 
In Section 5, we consider the case of the shear-thickening flows. We establish a version of the Korn inequality more appropriate for our framework and derive some estimates for the convective term and the extra stress tensor. The existence and uniqueness results for the full governing equations are then established and the approximation analysis is carried out. The shear-thinning case is treated in a similar way in Section 6.

\section{\large{Governing Equations}}
\label{section2_1}\setcounter{equation}{0}
\noindent We are concerned with steady flows of incompressible generalized Newtonian
 fluids. For these fluids, the Cauchy stress tensor $\widetilde {\boldsymbol T}$  is related to the kinematic variables through 
$$\widetilde 
	{\boldsymbol T}=-\widetilde \pi I+2\mu\left(1+|\widetilde D\widetilde { \boldsymbol u}|^2
\right)^{\frac{p-2}{2}}\widetilde D\widetilde { \boldsymbol u},$$
where $\widetilde  {\boldsymbol u }$ is the velocity field,
$\widetilde D\widetilde { \boldsymbol u}=\frac{1}{2}\left(\widetilde\nabla \widetilde { \boldsymbol u}+\widetilde\nabla \widetilde { \boldsymbol u}^\top\right)$ denotes the symmetric part of the velocity gradient, $\mu >0$ is the kinematic viscosity, and $\widetilde\pi$ represents the pressure. The notation $\sim$ denotes a dimensional quantity.
The equations of conservation of momentum and mass, relative to a rectangular coordinate system, are 
\begin{equation}\label{equation_dim}\left\{ 
\begin{array}{ll}\rho \left(\tfrac{\partial \widetilde { \boldsymbol u}}{\partial \widetilde t}+\widetilde { \boldsymbol u}\cdot \widetilde\nabla \widetilde { \boldsymbol u}\right) 
                             +\widetilde\nabla \widetilde \pi = \widetilde\nabla\cdot\left(2\mu\left(1+|\widetilde D\widetilde { \boldsymbol u}|^2
\right)^{\frac{p-2}{2}}\widetilde D\widetilde { \boldsymbol u}\right),\vspace{3mm}\\
\widetilde\nabla\cdot \widetilde { \boldsymbol u} = 0,\end{array}
\right.\end{equation}
where $\rho>0$ is the constant density of the fluid. In this work, we consider steady flow of generalized Newtonian fluids through a curved pipe of arbitrary shaped cross-section $\Sigma$ with  constant centerline radius $R$. Due to the geometric characteristics of the curved pipe, it is convenient to write system $(\ref{equation_dim})$ in the rectangular toroidal 
coordinates $(\widetilde x_i)$ defined with respect to the rectangular Cartesian coordinates $(\widetilde y_i)$ through the relations  
	\begin{equation}\label{transf_1}\widetilde x_1=\widetilde y_3,\qquad \widetilde x_2=\sqrt{\widetilde y_1^2+\widetilde y_2^2}
	  -R,\qquad \widetilde x_3=
	 R\arctan 
	  \tfrac{\widetilde y_2}{\widetilde y_1}
	\end{equation}
and inverse relations
	\begin{equation}\label{transf_2}\widetilde y_1=\left(R+\widetilde x_2\right)\cos\left(\tfrac{\widetilde x_3}{R}\right),\qquad 
	\widetilde y_2=\left(R+\widetilde x_2\right)\sin\left(\tfrac{\widetilde x_3}{R}\right),\qquad 
	\widetilde y_3=\widetilde x_1.\end{equation}
 More details on the toroidal coordinate system, and on the corresponding formulation of the operators involved in (\ref{equation_dim}), can be found in Appendix \ref{appendix_a}.\vspace{5mm}

\begin{figure}[htbp]
\SetLabels
\L (0.85*0.55)  \scriptsize $\widetilde{x}_{1}$ \\
\L (0.915*0.29)  \scriptsize $\widetilde{x}_{2}$ \\
\L (0.755*0.40)  \scriptsize $\widetilde{x}_{3}$ \\
\endSetLabels
\leavevmode
\epsfxsize=5cm
\epsfysize=4cm
$$
\strut\AffixLabels{\mbox{\psfig{figure=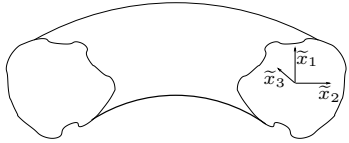,width=4.5cm}}}
$$
\caption{\small Rectangular toroidal coordinates in a curved pipe}
\end{figure}
\noindent We restrict our study to curved pipes with cross section independent of $\tilde x_3$, and consider flows which are steady and fully developed, i.e. the components of the velocity vector with respect to the new basis are independent of both time and axial variable $x_3$. For such flows
 	\begin{equation}\label{fdev_velocity}\tfrac{\partial  \widetilde u_i}{\partial \widetilde x_3}
	  =0 \qquad i=1,2,3\end{equation}
 and the axial component of the pressure gradient
	\begin{equation}\label{fdev_pressure}\tfrac{\partial \widetilde \pi}{\partial \widetilde x_3} =-\widetilde G \end{equation}
 is a constant. We consider the non-dimensional form of this system by
 introducing the following quantities 
	\begin{equation}\label{adim_var}r_0=\sup_{\widetilde
	x\in\overline\Sigma}|\widetilde x|, \qquad 
	x_i=\tfrac{\widetilde x_i}{r_0},\qquad \boldsymbol u =\tfrac{\widetilde {\boldsymbol  u }}{U_0}, \qquad 
	\pi=\tfrac{\widetilde \pi r_0}{\mu U_0}, \qquad G=\tfrac{\widetilde G r_0^2}{\mu U_0},\end{equation}
where $\delta\in [0,1[$ is the pipe curvature ratio, and
 $U_0$ represents a characteristic velocity of the flow.
 We also introduce the Reynolds number ${\cal R}e=\tfrac{\rho U_0 r_0}{\mu}$ and the pipe curvature ratio $\delta=\tfrac{r_0}{R}\in [0,1[$.\vspace{2mm}\\
  The  governing equations can then be written with respect to the non-dimensional quantities as
\begin{equation}\label{equation}
\left\{ \begin{array}{ll} -\nabla^\star \cdot\left({\boldsymbol\tau}(D^\star{\boldsymbol u})
	\right)+{\cal R}e \,   
                {\boldsymbol u}\cdot \nabla^\star {\boldsymbol u}+ \nabla^\star {\pi} =0 &\quad \mbox{in} \ 
	\Sigma,\vspace{3mm} \\
\nabla^\star \cdot {\boldsymbol u}=0&\quad \mbox{in} \ 
	\Sigma,\vspace{3mm} \\
	\boldsymbol u=0 &\quad \mbox{on} \ \partial\Sigma\end{array}\right.
\end{equation}
with
	$${\boldsymbol\tau}(D^\star{ \boldsymbol u})=2\left(1+ \dot{\gamma}^2|D^\star{ \boldsymbol u}|^2
\right)^{\frac{p-2}{2}}D^\star{ \boldsymbol u},$$
where $\dot{\gamma}=\tfrac{U_0}{r_0}$ is a caracteristic shear-rate. In order to simplify the redaction, we will assume without lost of generality that $\dot{\gamma}=1$.  The operators involved in the definition of (\ref{equation}), are defined by
	$$\begin{array}{ll}
	\nabla^\star\cdot {\boldsymbol\tau}=\nabla   \cdot {\boldsymbol\tau}+
        \tfrac{\delta}{B} \left({\tau}_{12}\,{\boldsymbol a }_1+
	\left({\tau}_{22}-{\tau}_{33}\right)\, {\boldsymbol a }_2+2{\tau}_{23}\,  {\boldsymbol a }_3
	\right),\vspace{3mm}\\ 
	\nabla^\star{ \boldsymbol u}=\nabla  { \boldsymbol u}+\tfrac{\delta}{B}\left(u_2\,
	    {\boldsymbol a }_3 \otimes{\boldsymbol a }_3-u_3\, {\boldsymbol a }_3\otimes{\boldsymbol a }_2\right),\vspace{3mm}\\
	D^\star{ \boldsymbol u}=\tfrac{1}{2}
	\left(\nabla^\star{ \boldsymbol u}+\nabla^\star{ \boldsymbol u}^\top\right),\vspace{3mm}\\
	{ \boldsymbol u}\cdot \nabla^\star
	   { \boldsymbol u}={ \boldsymbol u}\cdot \nabla  
	   { \boldsymbol u}+\tfrac{\delta}{B}\left(u_3u_2\, {\boldsymbol a }_3-
	   u_3^2\,{\boldsymbol a }_2\right),\vspace{3mm}\\
	\nabla^\star {\pi}=\nabla \pi-\tfrac{G}{B}\,{\boldsymbol a }_3\vspace{3mm}\\
	\nabla^\star \cdot { \boldsymbol u}=\tfrac{\partial u_1}{\partial x_1}+\tfrac{\partial u_2}{\partial x_2}+\tfrac{\delta}{B}\, u_2=\tfrac{1}{B}\, \nabla\cdot\left(B\boldsymbol u\right)\end{array}$$
with $$\nabla=\big(\tfrac{\partial}{\partial x_1},\tfrac{\partial}{\partial x_2},0\big), \qquad \nabla \cdot=\tfrac{\partial}{\partial x_1}+\tfrac{\partial}{\partial x_2},\qquad B=1+\delta x_2,$$and where  $({\boldsymbol a }_1,{\boldsymbol a }_2,{\boldsymbol a }_3)$
denotes the orthonormal basis in the toroidal coordinates. (For a detailed derivation of the dimensionless equation (\ref{equation}), see  Appendix \ref{appendix_b}.)

\section{Notation, assumptions and preliminary results}
\subsection{Algebraic results}
Throughout the paper, if ${ \boldsymbol u}=(u_1,u_2,u_3)$ in the rectangular toroidal coordinates, we denote by $u$ the vector with  toroidal components $(u_1,u_2,0)$. Similarly, if $\boldsymbol S=(S_{ij})_{i,j=1,2,3}$ is a tensor in $\mathbb{R}^{3\times 3}$, we denote by $S$ the tensor in $\mathbb{R}^{2\times 2}$ with toroidal components $S_{ij}$, $i=1,2$. \vspace{1mm}\\
For $\eta, \zeta\in \mathbb{R}^{d\times d}$, we define the scalar product and the corresponding norm by
	$$\eta:\zeta=\sum_{i,j=1}^d \eta_{ij}\zeta_{ij} \quad 
	\mbox{and} \quad \left|\eta\right|=
	\left(\eta:\eta\right)^{\frac{1}{2}}.$$
In the next two results, we state well known continuity, coercivity and monotonicity properties for $\boldsymbol\tau$.  The corresponding proof is standard (see e.g. \cite{necas}). We first consider the case $p\geq 2$ corresponding to the shear-thickening flows.
\begin{lemma} \label{tensor_proper2} Assume that $p\geq 2$ and let  $\eta\in \mathbb{R}_{\mathrm{sym}}^{3\times3}$. Then the  tensor $\boldsymbol\tau$  satisfies the following properties
\begin{description}
 \item {Continuity.}
	$$\left|{\boldsymbol\tau}(\eta)-\boldsymbol\tau(\zeta)\right|\leq 
	(p-1)\left(1+ |\eta|^2+|\zeta|^2\right)^{\frac{p-2}{2}}
	|\eta-\zeta|,$$
\item {Coercivity.}
	$${\boldsymbol\tau}(\eta):\eta\geq 2|\eta|^2,\qquad \quad
	 {\boldsymbol\tau}(\eta):\eta\geq 2|\eta|^p, $$
\item {Monotonicity.}
	$$\left\{\begin{array}{ll}
	\left({\boldsymbol\tau}(\eta)-{\boldsymbol\tau}(\zeta)\right):\left(\eta-\zeta\right)\geq\left|\eta-\zeta\right|^2,\vspace{3mm}\\
	\left({\boldsymbol\tau}(\eta)-{\boldsymbol\tau}(\zeta)\right):\left(\eta-\zeta\right)\geq \tfrac{1}{2^{p-1}(p-1)} 
	\left|\eta-\zeta\right|^p.\end{array}\right.$$
\end{description}
 \end{lemma} 
Next we consider the case $1<p<2$ corresponding to the shear-thinning flows.
\begin{lemma} \label{tensor_proper1} Assume that $1<p<2$ and let  $\eta\in \mathbb{R}_{\mathrm{sym}}^{3\times3}$. Then the  tensor $\boldsymbol\tau $  satisfies the following properties
\begin{description}
 \item {Continuity.}
$$\left|{\boldsymbol\tau }(\eta)-\boldsymbol\tau (\zeta)\right|\leq 
	C_p|\eta-\zeta|^{p-1} \qquad \mbox{with} \ 
	C_p=1+2^{\frac{2-p}{2}},$$
\item {Coercivity.}
$${\boldsymbol\tau }(\eta):\eta\geq 
	2\left(1+|\eta|^2\right)^{\frac{p-2}{2}}|\eta|^2,$$
\item {Monotonicity.}
	$$\left({\boldsymbol\tau }(\eta)-{\boldsymbol\tau }(\zeta)\right):
	\left(\eta-\zeta\right)
	\geq 2(p-1)\left(1+ |\eta|^2+ |\zeta|^2
	\right)^{\frac{p-2}{2}}\left|\eta-\zeta\right|^2.$$
\end{description}
	 \end{lemma} 
\subsection{Functional setting}
Throughout the paper $\Sigma$ is a bounded domain in $\mathbb{R}^2$, with a boundary $\partial \Sigma$. Even though several of our results are valid for an arbitrary bounded domain, we will assume without loss of generality that $\Sigma$ is of class $C^2$.
By $\boldsymbol W^{k,p}(\Sigma)$ ($k\in \mathbb{N}$ and $1<p<\infty$), we denote the standard Sobolev spaces 
and we denote the associated norms by $\|\cdot\|_{k,p}$. We set 
$\boldsymbol W^{0,p}(\Sigma)\equiv \boldsymbol L^{p}(\Sigma)$, $\|\cdot\|_{L^{p}}\equiv \|\cdot\|_{p}$,  $\boldsymbol L^{p}_0(\Sigma)=\left\{\boldsymbol u\in  L^{p}(\Sigma)\mid \int_\Sigma \boldsymbol u(x)\,dx=0\right\}$, and 
	$$\left\|\cdot\right\|_{p,B}=\left\|B^{\frac{1}{p}}\cdot\right\|_{p}, \qquad m=\left\|\tfrac{1}{B}\right\|_\infty  \qquad \mbox{and} \qquad
	n=\|B\|_\infty.$$
The dual space of  $\boldsymbol W^{1,p}_0(\Sigma)$ is denoted by
 $\boldsymbol W^{-1,p'}(\Sigma)$, where  $p'=\tfrac{p}{p-1}$  is the dual exponent to $p$, 
 its norm is denoted by $\|\cdot\|_{-1,p'}$ and the duality pairing between these spaces by $\langle\cdot,\cdot\rangle$.
We will also use the following notation
	$$\left(\boldsymbol u,\boldsymbol v\right)=\int_\Sigma { \boldsymbol u}(x)\cdot { \boldsymbol v}(x)\,dx, \qquad { \boldsymbol u}\in L^{p}(\Sigma), \ { \boldsymbol v}\in 
	\boldsymbol L^{p'}(\Sigma),$$
$$\left(\eta,\zeta\right)=\int_\Sigma \eta(x):\zeta(x)\,dx, \qquad \eta\in 
	\boldsymbol L^p(\Sigma,\mathbb{R}^{3\times 3}), \ \zeta\in 
	\boldsymbol L^{p'}(\Sigma,\mathbb{R}^{3\times 3}).$$
The space of infinitely differentiable functions with compact support in $\Sigma$ will be denoted by ${\cal D}(\Sigma)$. In order to eliminate the pressure in the weak formulation of our equation, we will work in divergence-free spaces. Consider
        $$ {\cal V}=\big\{\boldsymbol\varphi\in {\cal D}(\Sigma)
        \mid \nabla\cdot \boldsymbol\varphi=0\big\},$$
	$${\cal V}_B=\big\{\boldsymbol\varphi\in {\cal D}(\Sigma)
        \mid \nabla \cdot \left(B\boldsymbol\varphi\right)=0\big\},$$
and denote by $\boldsymbol V^{p}$ and $\boldsymbol V_B^p$ the closure of ${\cal V}$ and ${\cal V}_B$ in the $L^p$-norm of gradients, i.e.
	$${\boldsymbol V}^{p}=
        \left\{\boldsymbol\varphi\in \boldsymbol W_0^{1,p}(\Sigma)
        \mid \nabla\cdot \boldsymbol\varphi=0\right\}, $$
	$${\boldsymbol V}^{p}_B=
        \left\{\boldsymbol\varphi\in \boldsymbol W_0^{1,p}(\Sigma)
        \mid \nabla\cdot \left(B\boldsymbol\varphi\right)=0\right\}.$$
Next, we recall the Sobolev inequality.
\begin{lemma}  Let  $r$ and $q$ be such that $r\geq q$ if $q\geq 2$ and  
$1<r\leq q^\ast=\tfrac{2q}{2-q}$ if $1<q<2$. Then for all
all $ \boldsymbol u\in \boldsymbol W^{1,q}_0(\Sigma)$, we have
	\begin{equation}\label{sobolev0}
	\left\|{ \boldsymbol u}\right\|_r\leq  S_{q,r}
	 \left\|\nabla \boldsymbol u\right\|_q,
	\end{equation}
where 
	$$S_{q,r}=\left\{\begin{array}{ll}
\tfrac{\max\left(q,\frac{r}{2}\right)}{2\sqrt{2}}\, |\Sigma|^{\frac{1}{2}+\frac{1}{r}-\frac{1}{q}}& \quad \mbox{if} \ r\geq q \ \mbox{and} \ 
	q\geq 2, \vspace{2mm}\\
\tfrac{q^\ast}{4\sqrt{2}}\,
 |\Sigma|^{\frac{1}{2}+\frac{1}{r}-\frac{1}{q}}& \quad \mbox{if} \ 1<r\leq q \ \mbox{and} \ 1<q<2,\vspace{2mm}\\
\tfrac{\max\left(q,\frac{r}{2}\right)}{2\sqrt{2}}\, |\Sigma|^{\frac{1}{2}+\frac{1}{r}-\frac{1}{q}}&  \quad  \mbox{if} \ q<r\leq q^\ast  \ \mbox{and} \ 1<q<2.\end{array}\right.
	$$
\end{lemma}
{\bf Proof.} The Sobolev inequality (\ref{sobolev0}) is classical and follows by combining the following interpolation (see Lemma II.3.2 in \cite{galdi})
           $$\left\|\boldsymbol u\right\|_r\leq \left(\tfrac{\max\left(q,\frac{r}{2}\right)}{2\sqrt{2}}\right)^\lambda \left\|\boldsymbol u\right\|_q^{1-\lambda}
	\left\|\nabla\boldsymbol u\right\|_q^\lambda \qquad \mbox{with} \ \lambda=\tfrac{2(r-q)}{rq} $$
with the inequality $\left\|\boldsymbol u\right\|_q\leq |\Sigma|^{\frac{1}{q}-\frac{1}{r}} \left\|\boldsymbol u\right\|_r$.$
\hfill\Box$\vspace{2mm}\\
Next, we consider a particular version of the Poincar\'e  inequality.
\begin{lemma}  For all ${ \boldsymbol u}\in { \boldsymbol W}^{1,q}_0(\Sigma)$, $1\leq q\leq +\infty$, the following estimate holds
	\begin{equation} \label{poincare}\left\|{ \boldsymbol u}\right\|_q\leq
	 \left\|\tfrac{\partial { \boldsymbol u}}{\partial x_1}\right\|_q.\end{equation}
\end{lemma}
{\bf Proof.} The result follows from Theorem II.5.1 in \cite{galdi}  by observing that 
$\Sigma\subset ]-1,1[\times \mathbb{R}$.$\hfill\Box$\vspace{2mm}\\
As well known (see \cite{galdi}), the standard Poincar\'e inequality is given by
	\begin{equation}
	\label{poincare_opt}\left\|{ \boldsymbol u}\right\|_q\leq
	 S_{q,q} \left\|\nabla   { \boldsymbol u}\right\|_q \end{equation}
and consequently, the Poincar\'e constant in $(\ref{poincare})$ is not necessarily optimal.  However, since $B$ is independent of the variable $x_1$, a direct consequence of $(\ref{poincare})$ is that
	\begin{equation}\label{poincare2}\left\|{ \boldsymbol u}\right\|_{q,B}\leq \left\|\tfrac{\partial { \boldsymbol u}}{\partial x_1}\right\|_{q,B}\leq 
	\left\|\nabla   { \boldsymbol u}\right\|_{q,B}.\end{equation}
This property is particularly useful to establish the Korn inequality in the shear-thinning case. \vspace{2mm}\\
Finally, we establish some estimates useful in the approximation analysis of the solutions with respect to the parameter $\delta$.
\begin{lemma} \label{est_psi} Let $u$, $v$ be in $W^{1,q}_0(\Sigma)$ with $q\geq 2$, and let $\sigma$ be a continous function such that
	\begin{equation}\label{psi_assumption}|\sigma(\lambda)|\leq c_0|\lambda|^\alpha 
	\qquad \mbox{with} \ c_0>0.\end{equation}
Then for every $\alpha \geq 0$, we have
	\begin{equation} 
	\label{est_psi_1}\left|\left(\sigma(u),v\right)\right|\leq c_0 D_{q,\alpha}	
	\left\|\nabla u\right\|_{q}^\alpha
	 \left\|\nabla v\right\|_{q},\end{equation}
	\begin{equation} \label{est_psi_2}\left\|\sigma(u)\right\|_{q'}\leq c_0 E_{q,\alpha}
	\left\|\nabla u\right\|_2^\alpha,\end{equation}
where $D_{q,\alpha}=\max\left(\tfrac{q}{2\sqrt{2}},\tfrac{\alpha q}{4\sqrt{2}}\right)^{\alpha+1} |\Sigma|^{\frac{1}{q'}+\frac{\alpha+1}{2}-
\frac{\alpha}{q}}$ and $E_{q,\alpha}=\max\left(\tfrac{1}{\sqrt{2}},\tfrac{\alpha}{2\sqrt{2}}\right)^{\alpha} |\Sigma|^{\frac{1}{q'}}$.
\end{lemma}
{\bf Proof.} Let us first assume that $\alpha\geq 1$. Then, by using the H\"older inequality and the Sobolev inequality (\ref{sobolev0}) we obtain
	$$\begin{array}{ll}
	\left|\left(\sigma(u),v\right)\right|&\leq
	\left\|\sigma(u)\right\|_{q}\|v\|_{q'}\leq c_0 
	\left\|u\right\|_{\alpha q}^\alpha\|v\|_{q'}\leq c_0 
	|\Sigma|^{\frac{1}{q'}-\frac{1}{\alpha q}}\left\|u\right\|_{\alpha q}^\alpha\|v\|_{\alpha q}\vspace{2mm}\\
	&\leq c_0 |\Sigma|^{\frac{1}{q'}-\frac{1}{\alpha q}}
	\left( S_{q,\alpha q}\right)^{\alpha+1}
	 \left\|\nabla u\right\|_{q}^\alpha \|\nabla v\|_{q}\end{array}$$
which gives (\ref{est_psi_1}). Similarly, 
	$$\begin{array}{ll}\left\|\sigma(u)\right\|_{q'}
	&\leq |\Sigma|^{\frac{1}{q'}-\frac{1}{2}}
	\left\|\sigma(u)\right\|_{2}
	\leq c_0 |\Sigma|^{\frac{1}{q'}-\frac{1}{2}}
	\left\|u\right\|_{2\alpha}^\alpha\vspace{2mm}\\
	&\leq c_0|\Sigma|^{\frac{1}{q'}-\frac{1}{2}}
	\left( S_{2,2\alpha}\right)^{\alpha}
	 \left\|\nabla u\right\|_{2}^\alpha\end{array}$$
and we obtain (\ref{est_psi_2}).
If $\alpha\leq 1$ then
	$$\begin{array}{ll}
	\left|\left(\sigma(u),v\right)\right|&\leq
	\left\|\sigma(u)\right\|_{\frac{q}{\alpha}}\|v\|_{\left(\frac{q}{\alpha}\right)'}\leq c_0 
	\left\|u\right\|^\alpha_{q}\|v\|_{\frac{q}{q-\alpha}}\leq c_0 |\Sigma|^{\frac{1}{q'}-\frac{\alpha}{q}}
	\left\|u\right\|^\alpha_{q}\|v\|_{q}\vspace{2mm}\\
	&\leq c_0 |\Sigma|^{\frac{1}{q'}
	-\frac{\alpha}{q}}\left( S_{q,q}
	\right)^{\alpha+1} \left\|\nabla u\right\|_{q}^{\alpha}
	 \|\nabla v\|_{q}\end{array}$$
and 
	$$\begin{array}{ll}\left\|\sigma(u)\right\|_{q'}
	&\leq |\Sigma|^{\frac{1}{q'}-\frac{\alpha}{2}}
	\left\|\sigma(u)\right\|_{\frac{2}{\alpha}}
	\leq c_0 |\Sigma|^{\frac{1}{q'}-\frac{\alpha}{2}}
	\left\|u\right\|^\alpha_{2}\vspace{2mm}\\
	&\leq	c_0 |\Sigma|^{\frac{1}{q'}-\frac{\alpha}{2}}\left( S_{2,2}
	\right)^{\alpha} \left\|\nabla u\right\|_{2}^{\alpha}\end{array}$$
and the claimed result is proved.$\hfill\Box$
\begin{lemma} \label{est_psi_thin} Let $u,v\in W^{1,q}_0(\Sigma)$ with $\tfrac{3}{2}\leq q<2$, and let $\sigma$ be a continous function satisfying 
$(\ref{psi_assumption})$. Then 
\begin{equation}\label{est_psi_1_thin}\left|\left(\sigma(u),v\right)\right|\leq c_0 D_{q,\alpha}	
	\left\|\nabla u\right\|_{q}^\alpha
	 \left\|\nabla v\right\|_{q} \qquad \mbox{for every} \ \
	0\leq \alpha\leq q^\ast-1,\end{equation}
	\begin{equation}\label{est_psi_2_thin}
	\left\|\sigma(u)\right\|_{q'}\leq c_0 E_{q,\alpha}
	\left\|\nabla u\right\|_q^\alpha\qquad \mbox{for every} \ \ \tfrac{1}{q'}
	< \alpha\leq \tfrac{q^\ast}{q'},\end{equation}
where $D_{q,\alpha}=\left( S_{q,q^\ast}\right)^{\alpha+1}$, $E_{q,\alpha}=\left( S_{q,\alpha q'}\right)^{\alpha}$ and $q^\ast=\tfrac{2q}{2-q}$.
\end{lemma}
{\bf Proof.} Notice first that the Sobolev inequality (\ref{sobolev0}) is valid for $1<r\leq q^\ast$ if $1\leq q<2$. If $0<\alpha\leq q^\ast-1$, then 
$\sigma(u) v$ belongs to $L^1(\Sigma)$ and by using (\ref{sobolev0}), we obtain
	$$\begin{array}{ll}
	\left|\left(\sigma(u),v\right)\right|&\leq c_0
	\left\|u^\alpha\right\|_{\frac{q^\ast}{\alpha}}
	\|v\|_{q^\ast}\leq c_0
	\left\|u\right\|^\alpha_{q^\ast}
	\|v\|_{q^\ast}\vspace{2mm}\\
	&\leq c_0 \left( S_{q,q^\ast}\right)^{\alpha+1}
	\left\|\nabla u\right\|_{q}^\alpha
	\left\|\nabla v\right\|_{q}\end{array}$$
which gives (\ref{est_psi_1_thin}).
To prove the last estimate, we observe that if $\tfrac{1}{q'}< \alpha\leq \tfrac{q^\ast}{q'}$ then $1<\alpha q'\leq q^\ast$. By using (\ref{sobolev0}) we obtain
	$$\left\|\sigma(u)\right\|_{q'}\leq c_0
	\left\|u\right\|_{\alpha q'}^\alpha
	\leq c_0 \left( S_{q,\alpha q'}\right)^{\alpha}
	 \left\|\nabla u\right\|_{q}^\alpha$$
and the claimed result is proved.$\hfill\Box$
\section{Weak formulation and statement of the main results}
\setcounter{equation}{0}
\label{section_shear_thic_introd}
To give a sense to the weak solution of (\ref{equation}), let us recall that $B=1+\delta x_2$ does not depend on $x_1$ and notice that if $\boldsymbol S$ is a symmetric tensor, then we have
	$$\begin{array}{ll}
	  \left(\nabla^\star \cdot \boldsymbol S\right)_{1}& =
	    \tfrac{\partial S_{11}}{\partial x_1}+
	    \tfrac{\partial S_{21}}{\partial x_2}
	    +\tfrac{\delta}{B}S_{21}=\tfrac{1}{B}
	\left(\tfrac{\partial \left(B S_{11}\right)}{\partial x_1}+
	    \tfrac{\partial \left(B S_{21}\right)}{\partial x_2}\right)\vspace{1mm}\\
	&=
	    \tfrac{1}{B}\left(\nabla\cdot 
	    \left(B\boldsymbol S\right)\right)_1,\vspace{2mm}\\
          \left(\nabla^\star \cdot \boldsymbol S \right)_2 &=
	    \tfrac{\partial S_{12}}{\partial x_1}+
	    \tfrac{\partial S_{22}}{\partial x_2}+
	    \tfrac{\delta}{B}\left(S_{22}-S_{33}\right)
	=\tfrac{1}{B}\left(\tfrac{\partial \left(B S_{12}\right)}{\partial x_1}+
	    \tfrac{\partial \left(BS_{22}\right)}{\partial x_2}\right)-\tfrac{\delta}{B}\,S_{33}\vspace{1mm}\\
	&=
\tfrac{1}{B}\left(\nabla\cdot 
	    \left(B\boldsymbol S\right)\right)_2-
	    \tfrac{\delta}{B}\, S_{33},\vspace{2mm}\\
          \left(\nabla^\star \cdot \boldsymbol S \right)_3 &=
	    \tfrac{\partial S_{13}}{\partial x_1}+
	    \tfrac{\partial S_{23}}{\partial x_2}+
	    \tfrac{2\delta}{B}\,S_{23}=\tfrac{1}{B}\left(\tfrac{\partial \left(B S_{13}\right)}{\partial x_1}+
	    \tfrac{\partial \left(B S_{23}\right)}{\partial x_2}\right)+\tfrac{\delta}{B}\,S_{23}\vspace{1mm}\\
		&=
	    \tfrac{1}{B}\left(\nabla\cdot 
	    \left(B\boldsymbol S\right)\right)_3+
	    \tfrac{\delta}{B}\,S_{32},
	    \end{array}$$
that is
	$$\nabla^\star \cdot \boldsymbol S=\tfrac{1}{B}
	  \nabla\cdot \left(B\boldsymbol S\right)
	  +\tfrac{\delta}{B}\left(S_{23}\mathbf a_3-S_{33}
	  \mathbf a_2\right).$$
Therefore, if $\boldsymbol S$ belongs to in 
$\boldsymbol L^{p'}(\Sigma)$ and 
${\boldsymbol\varphi}=(\varphi_1,\varphi_2,\varphi_3)\in \boldsymbol V^p_B$, then an integration by parts shows that
	$$\begin{array}{ll}-\left(\nabla^\star\cdot \boldsymbol S,B{\boldsymbol\varphi}\right)
	&=-\left(\tfrac{1}{B}\nabla^\star\cdot 
	\left(B \boldsymbol S\right),B{\boldsymbol\varphi}\right)
	=-\left(\nabla^\star\cdot \left(B \boldsymbol S\right),
	{\boldsymbol\varphi}\right)\vspace{2mm}\\
	&=-\left(\nabla  \cdot \left(B\boldsymbol S\right),
	{\boldsymbol\varphi}\right)-
	\delta\left(S_{23},\varphi_3\right)+
	\delta\left(S_{33},\varphi_2\right)\vspace{1mm}\\
	&=-\displaystyle\sum_{i=1}^3
	\left(\tfrac{\partial (B S_{1i})}{\partial x_1}
	+\tfrac{\partial (B S_{2i})}{\partial x_2},
	\varphi_i\right)
	-\delta\left(S_{23},\varphi_3\right)+
	\delta\left(S_{33},\varphi_2\right)\vspace{1mm}\\
	&=\displaystyle\sum_{i=1}^3\sum_{j=1}^2\left(B S_{ij},
        \tfrac{\partial \varphi_i}{\partial x_j}\right)
	-\delta\left(S_{23},\varphi_3\right)+
	\delta\left(S_{33},\varphi_2\right)
	=\left(\boldsymbol S,B D^\star{\boldsymbol\varphi}\right).
	\end{array}$$
By taking the inner product of $(\ref{equation})_1$ and 
$B\boldsymbol\varphi$ and integrating over $\Sigma$, we obtain the following weak formulation. 
\begin{definition} Assume that $p\geq \tfrac{3}{2}$. A function ${ \boldsymbol u}\in \boldsymbol V^p_B$ is a weak solution of $(\ref{equation})$ if
	\begin{equation}\label{weak_formulation}\left({\boldsymbol\tau}(D^\star{ \boldsymbol u}),B D^\star{\boldsymbol\varphi}\right)+{\cal R}e 
	\left(B{ \boldsymbol u}\cdot \nabla^\star{ \boldsymbol u},
	{\boldsymbol\varphi}\right)=\left(G,\varphi_3\right)\qquad \mbox{for all} \ {\boldsymbol\varphi}\in {\boldsymbol V}_B^p.
	\end{equation}
\end{definition}
As in the case of bounded domains, this definition in meaningful for $p\geq \tfrac{3}{2}$ and will be used when considering both shear-thinning and shear-thickening flows. The restriction on the exponent $p$ ensures that the convective term belongs to $\boldsymbol L^1$ when considering test functions in $\boldsymbol V_B^p$. \vspace{1mm}\\
This formulation can be splitted into two {\it coupled} formulations, respectively associated to $u=(u_1,u_2,0)$ and to $(0,0,u_3)$. Indeed, by setting $\boldsymbol{\varphi}=(\varphi_1,\varphi_2,0)$ in (\ref{weak_formulation}), we can easily see that $u$ satisfies
	\begin{equation}\label{weak_formulation_u}
	\left(\tau(D^\star{ \boldsymbol u}),B D\varphi\right)+
	\delta\left({\tau}_{33}(D^\star{ \boldsymbol u}),\varphi_2\right)+
	{\cal R}e\, \left(Bu\cdot \nabla u,\varphi\right)
	={\cal R}e\,\delta\left(u_3^2,\phi_2\right),\end{equation}
where $\tau=({\tau}_{i,j})_{i,j=1,2}$. Similarly, by setting ${\boldsymbol\varphi}=(0,0,\varphi_3)$ in (\ref{weak_formulation}), we see that $u_3$ satisfies
	\begin{equation}\label{weak_formulation_u3}
	\begin{array}{ll}\left({\tau}_{13}(D^\star{ \boldsymbol u}),
	B\tfrac{\partial \varphi_3}{\partial x_1}\right)
	+\left({\tau}_{23}(D^\star{ \boldsymbol u}),
	B\tfrac{\partial \varphi_3}{\partial x_2}
	-\delta\varphi_3\right)+{\cal R}e\left((Bu\cdot \nabla u_3,
	\varphi_3)+\delta\left(u_3u_2,\varphi_3\right)\right)\vspace{2mm}\\
	=\left(G,\varphi_3\right).
	\end{array}
	\end{equation}
In the case of the Navier-Stokes equations, both formulations 
(\ref{weak_formulation}) and (\ref{weak_formulation_u})-(\ref{weak_formulation_u3}) can be used to establish the existence and uniqueness of a weak solution. In the present paper, the corresponding existence and uniqueness results are obtained as a particular case (by setting $p=2$) using the first formulation, while similar results were established in \cite{galdirobertson} by using the second formulation. One notable difference however is related with the corresponding estimates: involving  $\boldsymbol u$ and $w_3$ and independent of the Reynolds number in the first case, involving $u$ and $w_3$ (and thus $\boldsymbol u$) and depending on $\delta{\cal R}e$ (at least for $u$ and $\boldsymbol u$) in the second one. This fact has no influence on the study of the solvability of our problem, but a different estimate on $u$ is missed if we do not consider (\ref{weak_formulation_u}). When dealing with the generalized Navier-Stokes equations (\ref{equation}), the 
{\it global} formulation proves much more suitable than the splitting one. Because of the nonlinearity of the shear-stress, the coercivity property in the case $p<2$  and the monotonicity property in both cases $p<2$ and $p>2$ fail to be satisfied separately for $u$ and $u_3$. We will use the formulation (\ref{weak_formulation_u}) a posteriori to derive an additional estimate for $u$ that highlights the dependance on the curvature ratio and the Reynolds number.\vspace{1mm}\\
Existence and uniqueness of a weak solution in the shear-thickening case, together with associated estimates, is considered in the next result.
\begin{theorem}\label{main1} Assume that $p\geq 2$. Then problem 
$(\ref{equation})$ admits at least a weak solution $\boldsymbol u\in \boldsymbol V^p_B$, and this solution satisfies
	\begin{align}\label{main_estimates_1_thic}&
	\left\|D^\star{ \boldsymbol u}\right\|_{2,B}\leq \kappa_1,\\
	\label{main_estimates_2_thic}
	&\left\|D^\star{ \boldsymbol u}\right\|_{p,B}
	\leq\kappa_1^\frac{2}{p},\\
	\label{est_u3_thic}&\left\|\nabla u_3\right\|_{2,B}\leq 2^{\frac{1}{2}}\kappa_1,\\
	\label{est_u_2_thic}
	&\left\|Du\right\|_{2,B}\leq \kappa_2 \delta{\cal R}e,\\
	\label{est_u_p_thic}&\left\|Du\right\|_{p,B}\leq \left(\kappa_2 \delta{\cal R}e\right)^\frac{2}{p}\end{align}
with
	$$\kappa_1=\left(\tfrac{m}{2}\right)^{\frac{1}{2}}|G||\Sigma|\qquad \mbox{and} \qquad  
 \kappa_2=\tfrac{m^{\frac{3}{2}}|\Sigma|}{2}\kappa_1^2.
	$$
 Moreover, if the Reynolds number ${\cal R}e$ is such that
	\begin{equation}\label{Re_cond}
	{\cal R}e<\tfrac{2}{\kappa_1\kappa_3} 
	\qquad \mbox{with} \ \kappa_3=
	nm^\frac{3}{2} |\Sigma|^{\frac{3}{4}}\end{equation}
then problem $(\ref{equation})$ admits a unique weak solution.
\end{theorem}
 Unlike the estimates given in $(\ref{est_u_2_thic})$-$(\ref{est_u_p_thic})$, the estimates in $(\ref{main_estimates_1_thic})$-$(\ref{main_estimates_2_thic})$, and consequently in $(\ref{est_u3_thic})$, are uniformly bounded and neither depend on $\delta$ nor on ${\cal R}e$. On the other hand, the estimates in $(\ref{est_u_2_thic})$-$(\ref{est_u_p_thic})$ show that the secondary shear-thickening flows, if they exist, depend simultaneously on the pipe curvature ratio and on the Reynolds number or equivalently, and after introducing the standard adimensionalization $(\tilde u,\tilde u_3)=
(\sqrt{\delta} u,u_3)$, that they depend on the Dean number 
${\cal D}e=\sqrt{\delta}\, {\cal R}e$.  These results are specially useful when dealing with small Dean numbers and imply, in this situation, that the shear-thickening secondary flows are necessarily small. \vspace{2mm}\\
Similar results are obtained in the shear-thinning case.
\begin{theorem}\label{main2} Assume that $\tfrac{3}{2}\leq p<2$. Then problem 
$(\ref{equation})$ admits at least a weak solution $\boldsymbol u\in \boldsymbol V^p_B$, and this solution satisfies
	\begin{align}\label{main_estimates_1_thin}
	&\left\|D^\star{ \boldsymbol u}\right\|_{p,B}\leq \kappa_2,\\
	&\label{main_estimates_2_thin}
	\left\|D^\star{ \boldsymbol u}\right\|_{p,B}^p
	\leq \kappa_3,\\
	&\label{est_u3_thin}\left\|\nabla u_3\right\|_{p,B}\leq \kappa_4,\\
	&\label{est_u_thin}\left\|D u\right\|_{p,B}\leq \kappa_5 \,\delta {\cal R}e
	\end{align}
with
	$$\begin{array}{ll}\kappa_1=\tfrac{m^{\frac{1}{p}}}{2}|G|
	|\Sigma|^{\frac{1}{p'}}, \qquad 
	 \kappa_2=\kappa_1 \left(\tfrac{\|B\|_1}{p-1}+\kappa_1
	^{p'}\right)^{\frac{2-p}{p}}, 
	\qquad \kappa_3=p'\|B\|_1+\left(2^{\frac{2-p}{2}}\kappa_1\right)^{p'}, \vspace{2mm}\\
	\kappa_4=2^{2-p}\left(1+\delta m \right)\kappa_2,\end{array}$$
and where $\kappa_5$ is a positive constant only depending on $\Sigma$, $p$, $G$ and $m$.
Moreover,  there exists a positive constant $\kappa_6$  only depending on $\Sigma$, $p$, $G$, $m$ and $n$ such that if
	\begin{equation}\label{restriction_control}{\cal R}e<\tfrac{1}{2\kappa_1\kappa_6}\left(\tfrac{\|B\|_1}{p-1}+\kappa_1^{p'}\right)^{-\frac{2(2-p)}{p}},
	\end{equation}
then problem $(\ref{equation})$ admits a unique weak solution.
\end{theorem}
As in the case of shear-thickening flows, the estimates in $(\ref{main_estimates_1_thin})$-$(\ref{est_u3_thin})$ are uniformly bounded and independent of $\delta$ and ${\cal R}e$ while estimate $(\ref{est_u_thin})$ depends simultaneously on these two parameters and shows that if the Dean number is small, then the shear-thinning secondary flows are necessarily small.
\begin{remark} In the statement of Theorem $\ref{main2}$, the dependence of the constants $\kappa_5$ and $\kappa_6$ on the  parameters $\delta$ and $m$ is explicitely known. More precisely, we have
	$$\begin{array}{ll} \kappa_5=\tfrac{m^{\frac{3}{p}}}{C_{K,1}}\left( S_{p,2p'}\right)^3 \kappa_4^2\left(\tfrac{\|B\|_1}{p-1}+\kappa_1^{p'}\right)^{\frac{2-p}{p}},\vspace{2mm}\\
	\kappa_6=\tfrac{8n^{\frac{p+1}{p}} m^{\frac{6}{p}}(1+\delta m)^3}{C_{K,1}^3}\left( S_{p,2p'}\right)^2,\end{array}$$
where $C_{K,1}$ is the classical Korn inequality in $\boldsymbol W^{1,p}_0(\Sigma)$. Notice also that the constants $\kappa_2$, $\kappa_3$, $\kappa_4$, $\kappa_5$ and $\kappa_6$ depend on $\kappa_1$. 
\end{remark}
Having a weak solution, the corresponding term $\nabla\pi$ can be constructed by the same way as in the linear case. The pressure is determined up to a constant and becomes unique under the additional condition $\int_\Sigma \pi\,dx=0$.
\begin{corollary}\label{pressure_1} Assume that $p\geq 2$ and that $\boldsymbol u\in \boldsymbol V_B^p$ is a weak solution of 
$(\ref{equation})$. Then there exists a unique $\pi \in L^{p'}_0(\Sigma)$ such that $(\ref{equation})_1$ holds in $\boldsymbol W^{-1,p'}(\Sigma)$. Moreover, we have the following estimate
	$$\left\|\pi\right\|_{p'}\leq \kappa \left(\left\|Du\right\|_{p,B}+\delta^2\left\|u_2\right\|_{p,B}+{\cal R}e \left(\left\|D u\right\|_{p,B}^2+	\delta \left\|\nabla u_3\right\|_{p,B}^2\right)\right)$$
where  $\kappa$ is a positive constant only depending on $\Sigma$, $p$, $m$ and $n$.
\end{corollary}
Similarly, existence of the pressure in the shear-thinning case is considered below.
\begin{corollary}\label{pressure_1_thin} Assume that $\tfrac{3}{2}\leq p<2$ and that $\boldsymbol u\in \boldsymbol V_B^p$ is a weak solution of 
$(\ref{equation})$. Then there exists a unique $\pi \in L^{p'}_0(\Sigma)$ such that $(\ref{equation})_1$ holds in $\boldsymbol W^{-1,p'}(\Sigma)$. Moreover, we have the following estimate
	$$\left\|\pi\right\|_{p'}\leq \kappa \left(\left\|Du\right\|_{p,B}^{p-1}+\delta^p   \left\|u_2\right\|_{p,B}^{p-1}+{\cal R}e \left(\left\|D u\right\|_{p,B}^2+
	\delta \left\|\nabla u_3\right\|_{p,B}^2\right)\right),$$
where  $\kappa$ is a positive constant only depending on $\Sigma$, $p$, $m$ and $n$.
\end{corollary}
The next results deal with special properties of the solutions of $(\ref{equation})$. We recall that a solution $(\boldsymbol u=(u,u_3),\pi)$ is unidirectional flow if 
$u_1=u_2=0$.
\begin{proposition} Assume that $p\geq \tfrac{3}{2}$. If $\delta{\cal R}e=0$, then all the solutions of $(\ref{equation})$ are unidirectional flows.
\end{proposition}
{\bf Proof.} Let us first assume that $p\geq 2$. Taking into account (\ref{est_u_p_thic}), we deduce that
	$\left\|\nabla u\right\|_p=0$ and 
since $u_{\mid \partial \Sigma}=0$, by using the Poincar\'e inequality  it follows that 
$u=0$ and thus $\boldsymbol u=(0,0,u_3)$ 
with $u_3$ satisfying
	\begin{equation}\label{u3_uni}\left\{\begin{array}{ll}-\nabla^\star\left(\left(1+\tfrac{1 }{2}|\nabla^\star u_3|^2\right)^{\frac{p-2}{2}}\nabla^\star u_3\right)=\tfrac{G}{B} & \mbox{in} \ \Sigma,\vspace{2mm}\\
	u_3=0 & \mbox{in} \ \partial\Sigma.\end{array}\right.\end{equation}
The result corresponding to the shear-thinning case can be proved similarly using the estimate
 (\ref{est_u_thin}) instead of the estimate (\ref{est_u_p_thic}). $\hfill\Box$
\begin{proposition} Assume that $p\geq \tfrac{3}{2}$. If $\delta{\cal R}e>0$, then the solutions of $(\ref{equation})$  are not  unidirectional flows.
\end{proposition}
{\bf Proof.} Let us assume that $(u_1,u_2)=(0,0)$. System (\ref{equation}) reduces to
	$$\nabla\pi={\cal R}e\,\tfrac{\delta}{B}u_3{\boldsymbol a}_2	\quad\mbox{and} \quad 
	u_3 \ \mbox{satisfies} \  (\ref{u3_uni}).$$
If $\delta{\cal R}e>0$, the first equation implies that $\pi$ (and thus $u_3$) does not depend on the variable $x_1$. Following the arguments developed in \cite{galdirobertson}, we can prove that 
	\begin{equation}\label{u3zero}u_3=0 \qquad \mbox{in} \ \Sigma.
	\end{equation}
Indeed, let $(x_1,x_2)$ be an arbitrary point in $\Sigma$ and denote by 
$(x_0,x_2)$ a point on $\partial \Sigma$ that is the intersection of the straight line with origin $(x_1,x_2)$, parallel to $\boldsymbol a_1$ and such that the segment $\mbox{seg}= \cup_{\alpha\in [0,1]}$ lies in $\overline \Sigma$. Since $u_3$ is independent of $x_1$, it follows that 
$u(x_1^\alpha,x_2)$ is independent of $\alpha$ for $x_1^\alpha\in \mbox{seg}$. Therefore, by taking into account the boundary condition $u_{\mid \Sigma}=0$, we obtain
	$$u_3(x_1,x_2)=u_3(x_1^1,x_2)=u_3(x_1^0,x_2)=u_3(x_0,x_2)=0.$$
The point $(x_1,x_2)$ being arbitrary in $\Sigma$, we deduce that (\ref{u3zero}) holds, and this contradicts  $(\ref{u3_uni})$ and completes the proof.$\hfill\Box$\vspace{2mm}\\
Let us now analyze the behavior of the weak solutions of $(\ref{equation})$ with respect to the parameter $\delta$. The objective would be to use these results when $\delta$ is small to approximate a solution $\boldsymbol u$ of (\ref{equation}) by a solution of a similar but simpler system. More precisely, we consider the following problem
$$(E_\sigma) \qquad 
\left\{ \begin{array}{ll} -\nabla \cdot\left({\boldsymbol\tau}(D\boldsymbol w)
	\right)+{\cal R}e \,   
                {\boldsymbol w}\cdot \nabla {\boldsymbol w}+ \nabla {\pi} =\tfrac{G}{B}\boldsymbol a_3+\delta \sigma(w_3) \boldsymbol a_2&\quad \mbox{in} \ 
	\Sigma,\vspace{2mm} \\
\nabla \cdot {\boldsymbol w}=0&\quad \mbox{in} \ 
	\Sigma,\vspace{2mm} \\
	\boldsymbol w=0 &\quad \mbox{on} \ \partial\Sigma,\end{array}\right.
$$
where $\sigma$ is a non constant function and we aim to estimate the difference $\boldsymbol u-\boldsymbol w$. Obviously, the considerations concerning the monotonicity and coercivity properties of the {\it global} and the {\it coupled} formulations described above and the difficulties encountered in the treatment of the full governing equations arise in a similar way for $(E_\sigma)$. Moreover, in the derivation of the {\it a priori} estimate, the term involving $\sigma$ induces an additional difficulty that can be overcome by carrying out a careful analysis. \vspace{1mm}\\
Begining with the shear-thickening case, we summarize the properties of the solutions of $(E_\sigma)$. 
\begin{proposition} \label{w_exist_thick} Assume that $p\geq 2$ and let 
$\sigma$ be a non constant continous function satisfying 
$(\ref{psi_assumption})$  for some $\alpha\geq 0$. Then problem $(E_\sigma)$ admits at least a weak solution $\boldsymbol w\in \boldsymbol V^p$, and this solution satisfies
	$$\begin{array}{ll}\left\|\nabla w_3\right\|_2\leq c_1,
	&\qquad \left\|\nabla w_3\right\|_p^p\leq 2^{\frac{p-2}{2}}
	c_1^2,\vspace{2mm}\\
	\left\|D w\right\|_2\leq c_0c_2\delta,&\qquad 
	\left\|D w\right\|_p^p\leq \left(c_0c_2\delta\right)^2,
	\end{array}$$
where 
	$$c_1=\tfrac{m|\Sigma||G|}{\sqrt{2}} \quad \mbox{and} \quad 
	c_2= \tfrac{D_{2,\alpha}c_1^\alpha}{\sqrt{2}} \quad
	 \mbox{with} \ D_{2,\alpha} \ \mbox{given in Lemma} \ \ref{est_psi}.$$
Moreover, if $\delta>0$ then the solutions of $(E_\sigma)$ are not unidirectional flows.
\end{proposition}
Next we state the corresponding approximation result.
\begin{proposition} \label{delta_vs_zero} Assume that the assumptions of Proposition $\ref{w_exist_thick}$ are fulfilled and let $\boldsymbol u$, $\boldsymbol w$ be the solutions of $(\ref{equation})$ and $(E_\sigma)$. There exists ${\cal R}e_0>0$ such that if ${\cal R}e\leq {\cal R}e_0$ then 
	$$\left\|D\left(\boldsymbol u-\boldsymbol w\right)\right\|_{p}^p+
	\left\|D\left(\boldsymbol u-\boldsymbol w\right)
	\right\|_{2}^2\leq 
	\kappa\,\delta^{p'}$$
where $\kappa$ depends only on $p$, $\Sigma$, $m$, $n$, $c_0$ and $\alpha$.
\end{proposition}
Similarly, we consider the solvability of $(E_\sigma)$ in the shear-thinning case and the corresponding approximation result.
\begin{proposition} \label{w_exist_thin} Assume that $\tfrac{3}{2}\leq p<2$ and let $\sigma$ is a non constant continuous function  satisfying $(\ref{psi_assumption})$ for some $\alpha$ such that $\tfrac{1}{p'}<\alpha< \tfrac{p^\ast}{2p'}$. 
 Then problem $(E_\sigma)$ admits at least a weak solution $\boldsymbol w\in \boldsymbol V^p$, and this solution satisfies
	$$\begin{array}{ll}\left\|\nabla w_3\right\|_p\leq c_1
\left(|\Sigma|+c_4^p\right)^{\frac{2-p}{p}},\vspace{2mm}\\
	\left\|D w\right\|_p\leq 
	c_0c_1^\alpha c_2 
	\left(|\Sigma|+c_4^p\right)^{\frac{(2-p)(\alpha+1)}{p}}\, \delta,
	\end{array}$$
where 
	$$\begin{array}{ll}c_1=m|\Sigma|^{\frac{1}{p'}}|G|, \qquad 
	c_2=\tfrac{D_{p,\alpha}}{2C_{K,1}} \quad
	 \mbox{with} \ D_{p,\alpha} \ \mbox{given in Lemma} \ \ref{est_psi_thin},
	 \vspace{2mm}\\
	c_3=\left(c_1+
	c_0c_1^{\alpha}c_2 \right)
	\left(1+|\Sigma|\right)^{\frac{(2-p)(\alpha+1)}{p}},\vspace{2mm}\\
	 c_4=
	 c_3\left(\tfrac{1}{1-(2-p)(\alpha+1)}+
	c_3^{\frac{1}{1-(2-p)(\alpha+1)}}
	\right)^{\frac{(2-p)(\alpha+1)}{p}}.\end{array}$$ 
Moreover, if $\delta>0$ then the solutions of $(E_\sigma)$ are not unidirectional flows.
\end{proposition}
\begin{proposition} \label{delta_vs_zero_thin} Assume that the assumptions of Proposition $\ref{w_exist_thin}$ are fulfilled and let $\boldsymbol u$, $\boldsymbol w$ be the solutions of $(\ref{equation})$ and $(E_\sigma)$. There exists ${\cal R}e_0>0$ such that if ${\cal R}e\leq {\cal R}e_0$ then 
	$$\left\|D\left(\boldsymbol u-\boldsymbol w\right)\right\|_p\leq 
	\kappa\,\delta^{p-1}$$
where $\kappa$ depends only on $p$, $\Sigma$, $m$, $n$, $c_0$ and $\alpha$.
\end{proposition}
Propositions \ref{w_exist_thick}, \ref{delta_vs_zero}, \ref{w_exist_thin} and \ref{delta_vs_zero_thin} 
 show that for $\delta>0$ sufficiently small, a solution of (\ref{equation}) can be approximated by a solution of $(E_\sigma)$, whose secondary flows exist even if they are proportionately weak. It is worth observing that this result is valid for a relatively large class of functions $\sigma$ and raises an interesting question related with the possible choices for $c_0$ and $\alpha$ that would guarantee an optimal approximation, in a sense to be correctly and adequately defined.\vspace{2mm}\\
We finish this section by considering the case of Navier-Stokes equations obtained by setting $p=2$. Notice that the constants $\kappa_1$, $\kappa_2$  and $\kappa_3$ in the statement of Theorem $\ref{main1}$ are independent of the exponent $p$ and that the condition that guarantees the uniqueness of weak solutions only depends on
$\Sigma$, $G$, $m$ and $n$. As a consequence, the estimates and the sufficient condition on ${\cal R}e$ are identical  in the particular case of Newtonian fluids.
\begin{theorem}\label{main_NS} The Navier-Stokes problem 
	$$\left\{ \begin{array}{ll} 
	-\nabla^\star \cdot\left(2D^\star{\boldsymbol u}
	\right)+{\cal R}e \,   
                {\boldsymbol u}\cdot \nabla^\star {\boldsymbol u}+ \nabla {\pi} =\tfrac{G}{B} \, \boldsymbol a_3 &\quad \mbox{in} \ 
	\Sigma,\vspace{2mm} \\
\nabla^\star \cdot {\boldsymbol u}=0&\quad \mbox{in} \ 
	\Sigma,\vspace{2mm} \\
	\boldsymbol u=0 &\quad \mbox{on} \ \partial\Sigma,\end{array}\right.$$
 admits at least a weak solution $\boldsymbol u\in \boldsymbol V^2_B$.
This solution satisfies the estimates $(\ref{main_estimates_1_thic})$, $(\ref{est_u3_thic})$, $(\ref{est_u_2_thic})$ and if ${\cal R}e$ satisfies $(\ref{Re_cond})$, then the solution is unique. If $\delta{\cal R}e=0$, then all the solutions are unidirectional flows, otherwise they are not unidirectional flows. Finally, let $\sigma$ be a non constant continuous function satisfying $(\ref{psi_assumption})$ for some $\alpha\geq 0$. Then the following problem
	$$\left\{ \begin{array}{ll} 
	-\Delta \boldsymbol w+{\cal R}e \,   
                {\boldsymbol w}\cdot \nabla {\boldsymbol w}+ \nabla {\pi} =\tfrac{G}{B}\boldsymbol a_3+\delta \sigma(w_3) \boldsymbol a_2 &\quad \mbox{in} \ 
	\Sigma,\vspace{2mm} \\
\nabla \cdot {\boldsymbol w}=0&\quad \mbox{in} \ 
	\Sigma,\vspace{2mm} \\
	\boldsymbol w=0 &\quad \mbox{on} \ \partial\Sigma,\end{array}\right.$$
admits a weak solution in $\boldsymbol V^2$. Moreover, there exists ${\cal R}e_0>0$ such that if ${\cal R}e\leq {\cal R}e_0$ then
	$$\left\|D(\boldsymbol u-\boldsymbol w)\right\|_2\leq
	 \kappa \delta,$$
where $\kappa$ depends on $\Sigma$, $m$, $n$, $c_0$ and $\alpha$.
	\end{theorem}
The case of Navier-Stokes equations has been fully studied in \cite{galdirobertson}. In the previous result, we recover similar results with some differences concerning the analysis with respect to $\delta$. On the one hand, our estimate is valid for a class of problems larger than the classical Dean problem obtained by setting 
	\begin{equation}\label{sigma_NS}
	\sigma(\lambda)={\cal R}e \lambda^2.\end{equation}
On the other hand, the estimate corresponding to the secondary flows is less accurate. Indeed, after introducing the adimensionalization 
$(\tilde u,\tilde u_3)=(\sqrt{\delta} u,u_3)$, we obtain
	$$\left\|D(\tilde u-\tilde w)\right\|_2\leq
	 \kappa \sqrt{\delta}, \qquad 
	 \left\|\nabla(\tilde u_3-\tilde w_3)\right\|_2\leq
	 \kappa \delta$$
while the estimate obtained in \cite{galdirobertson} reads as
	$$\left\|D(\tilde u-\tilde w)\right\|_2\leq
	 \kappa \delta, \qquad 
	 \left\|\nabla(\tilde u_3-\tilde w_3)\right\|_2\leq
	 \kappa \delta.$$
This is due to some technical difficulties mainly related with the combined effect of $Du$ and $\nabla u_3$ in the shear-rate and its consequences on the monotonicity properties of the tensor $\boldsymbol\tau$. Indeed, in the case of Navier-Stokes equations with $\sigma$ given by 
(\ref{sigma_NS}), the corresponding {\it coupled formulations} allow to derive, in a first step, estimates for $\left\|D(\tilde u-\tilde w)\right\|_2$ and $\left\|\nabla(\tilde u_3-\tilde w_3)\right\|_2$ dependent on one another. The combination of these estimates in a second step gives the result. In the case of a shear-dependent viscosity, and as already observed concerning the solvability of problem (\ref{equation}), the lack of monotonicity of the tensores $\tau=(\tau_{ij})_{i,j=1,2}$, 
$\tau_{13}$ and 
$\tau_{23}$ prevents from using the same arguments, and the global estimates we obtain come with a cost.
\section{Shear-thickening flows}
\setcounter{equation}{0}
The aim of this section is to study the case of shear-thickening flows (corresponding to $p\geq 2$). To achieve this goal, we first establish a  Korn inequality, and then estimate the convective term as well as the extra stress tensor in an adequate setting. We finally prove the corresponding main results given above.
\subsection{On the Korn inequality}\label{korn_thick}
The next result deals with an inequality of Korn's  type in ${ \boldsymbol H}^1_0(\Sigma)$, very similar to the classical one but involving the operators $\nabla^\star$ and $D^\star$.
\begin{lemma} \label{lm-korn}Let $ { \boldsymbol u}=(u_1,u_2,u_3) \in  { \boldsymbol H}^1_0(\Sigma)$. Then
\begin{equation}\label{korn-l2}\left\|\nabla^\star{ \boldsymbol u}\right\|_{2,B}^2=2\left\|D^\star{ \boldsymbol u}\right\|_{2,B}^2-
       \left\|\tfrac{1}{B}\nabla  \cdot\left(B u \right)\right\|_{2,B}^2.\end{equation}
\end{lemma}
{\bf Proof.} The definition of $D^\star$ together with standard calculations show that
	\begin{align}\label{korn1}
	2\left\|D^\star{ \boldsymbol u}\right\|_{2,B}^2&=2\left\|D  u \right\|_{2,B}^2
	+2\left\|\tfrac{\delta}{\sqrt{B}}\,u_2\right\|_2^2
	+\left\|\tfrac{\partial u_3}{\partial x_1}\right\|_{2,B}^2
	+\left\|\sqrt{B}\,\tfrac{\partial u_3}{\partial x_2}
	-\tfrac{\delta}{\sqrt{B}} u_3\right\|_2^2\nonumber\\
	&=2\left\|D u \right\|_{2,B}^2
	+2\left\|\tfrac{\delta}{\sqrt{B}} u_2\right\|_2^2+
	\left\|\nabla   u_3\right\|_{2,B}^2
	+\left\|\tfrac{\delta}{\sqrt{B}}u_3\right\|_2^2\nonumber\\
	&=2\left\|Du \right\|_{2,B}^2
	+2\left\|\tfrac{\delta}{B} u_2\right\|_{2,B}^2+
	\left\|\nabla   u_3\right\|_{2,B}^2
	+\left\|\tfrac{\delta}{B}u_3\right\|_{2,B}^2,\end{align}
where $u=(u_1,u_2,0)$. Since
	$$\nabla  \cdot\left(2D u \right)
         -\Delta  u =\nabla  \cdot \left((\nabla   u )^T
	 \right)=\nabla  \left(\nabla  \cdot  u \right),$$
we deduce that for every $ \boldsymbol\varphi\in  { \boldsymbol H}^1_0(\Sigma)$ we have
	$$ \left(\nabla   u , \nabla  \boldsymbol\varphi\right)
	 = 2\left(D  u ,
	 \nabla  \boldsymbol\varphi\right)-\left(\nabla  \cdot u ,\nabla  
	 \cdot\boldsymbol\varphi\right).$$
On the other hand, easy calculations show that $Du :\nabla  \boldsymbol\varphi =Du :D\boldsymbol\varphi$. Combining these identities, we obtain
	$$
       \left(\nabla   u , \nabla  \boldsymbol\varphi\right)
        =2\left( D   u ,D  \boldsymbol\varphi\right)-
	\left(\nabla  \cdot u ,\nabla  
	 \cdot\boldsymbol\varphi\right),$$
and thus 
	$$\begin{array}{ll}\left(\nabla u,
	 \nabla  \left(B\boldsymbol\varphi\right)\right)&= \left(B\nabla   u ,\nabla  \boldsymbol\varphi\right)
	+\delta\left(\tfrac{\partial  u }{\partial x_2},
	\boldsymbol\varphi\right)\vspace{2mm}\\
	&=2\left(D   u ,
	 D  \left(B\boldsymbol\varphi\right)\right)-
	 \left(\nabla  \cdot u ,\nabla  
	 \cdot\left(B\boldsymbol\varphi\right)\right)\vspace{2mm}\\
	&=2\left(B D   u ,
	 D  \boldsymbol\varphi\right)+\delta\left(\tfrac{\partial  u }
	 {\partial x_2}+\nabla  u_2,\boldsymbol\varphi\right)-
	 \left(\nabla  \cdot u ,\nabla  
	 \cdot\left(B\boldsymbol\varphi\right)\right)\vspace{2mm}\\
	&=2\left(B D   u ,
	 D  \boldsymbol\varphi\right)+\delta\left(\tfrac{\partial  u }
	 {\partial x_2},\boldsymbol\varphi\right)
	-\delta\left(u_2,\nabla  \cdot
	 \boldsymbol\varphi\right)-
	 \left(\nabla  \cdot u ,\nabla  
	 \cdot\left(B\boldsymbol\varphi\right)\right)\vspace{2mm}\\
	&=2\left(B D   u ,
	 D  \boldsymbol\varphi\right)+\delta\left(\tfrac{\partial  u }
	 {\partial x_2},\boldsymbol\varphi\right)+
	 \left(\tfrac{\delta^2}{B}u_2,\varphi_2\right)-\left(\tfrac{1}{B}\,\nabla  \cdot\left(B u \right),\nabla  
	 \cdot\left(B\boldsymbol\varphi\right)\right).\end{array}$$
The last equality implies that
	$$\left(B\nabla   u ,\nabla  \boldsymbol\varphi\right)=
       2\left(B D   u , D  \boldsymbol\varphi\right)-
       \left(\tfrac{1}{B}\nabla  \cdot\left(B u \right),\nabla  
	 \cdot\left(B\boldsymbol\varphi\right)\right)+
	 \left(\tfrac{\delta^2}{B}u_2,v_2\right)$$
for all $ \boldsymbol\varphi\in  { \boldsymbol H}^1_0(\Sigma)$. Setting $\boldsymbol\varphi= u $, we obtain
	\begin{align}\label{korn2}\left\|\nabla   u \right\|_{2,B}^2&= 
       2\left\| D   u \right\|_{2,B}^2+
	 \left\|\tfrac{\delta}{\sqrt{B}}\,u_2\right\|_2^2-
       \left\|\tfrac{1}{\sqrt{B}}\nabla  \cdot\left(B u \right)\right\|_2^2\nonumber\\
	&= 2\left\| D   u \right\|_{2,B}^2+
	 \left\|\tfrac{\delta}{B}\,u_2\right\|_{2,B}^2-
       \left\|\tfrac{1}{B}\nabla  \cdot\left(B u \right)\right\|_{2,B}^2.\end{align}
On the other hand, by taking into acount the definition of $\nabla^\star$ we have
	\begin{equation}\label{grad_prime}\left\|\nabla^\star{ \boldsymbol u}\right\|_{2,B}^2=
	\left\|\nabla  { \boldsymbol u}\right\|_{2,B}^2+\left\|\tfrac{\delta}{B}u_2\right\|_{2,B}^2
	+\left\|\tfrac{\delta}{B}u_3\right\|_{2,B}^2.\end{equation}
The conclusion follows from (\ref{korn1}), (\ref{korn2}) and (\ref{grad_prime}).$\hfill\Box$
\begin{remark} \label{remark-korn}A direct consequence of Lemma $\ref{lm-korn}$ is that
	$$\left\|\nabla^\star { \boldsymbol u} \right\|_{2,B}^2\leq 
       2\left\|D^\star { \boldsymbol u} \right\|_{2,B}^2\qquad \mbox{for all} \  { \boldsymbol u} \in  { \boldsymbol H}^1_0(\Sigma) $$
and the equality holds if $\nabla\cdot\left(B  \boldsymbol u \right)=0$.
\end{remark}
\subsection{Estimates on the convective term}
We point out some notable facts related with the trilinear forms $a$ and $a_\star$ defined by
	$$\begin{array}{ll}
	a({ \boldsymbol u},{ \boldsymbol v},{ \boldsymbol w})
	=\left({ \boldsymbol u}\cdot \nabla { \boldsymbol v},
	{ \boldsymbol w}\right),\vspace{2mm}\\
	a_\star({ \boldsymbol u},{ \boldsymbol v},{ \boldsymbol w})
	=\left({ \boldsymbol u}\cdot \nabla^\star { \boldsymbol v},{ \boldsymbol w}\right)=
	a({ \boldsymbol u},{ \boldsymbol v},{ \boldsymbol w})+
	\left(\tfrac{\delta}{B}v_3u_2,w_3\right)
	-\left(\tfrac{\delta}{B}v_3u_3,w_2\right).\end{array}$$
\begin{lemma} \label{prop_a1}Assume that $p\geq \tfrac{3}{2}$. For every ${ \boldsymbol u}\in {\boldsymbol V}_B^p$ and every ${ \boldsymbol v}$, ${ \boldsymbol w}\in  { \boldsymbol W}^{1,p}_0(\Sigma)$, we have
	$$a_\star(B { \boldsymbol u},{ \boldsymbol v},{ \boldsymbol v})=0\quad \mbox{and} \quad 
	a_\star(B { \boldsymbol u},{ \boldsymbol v},{ \boldsymbol w})
	=-a_\star(B { \boldsymbol u},{ \boldsymbol w},{ \boldsymbol v}).$$
\end{lemma}
{\bf Proof.} Taking into account the definition of $a_\star$ and the fact that $\nabla\cdot \left(B { \boldsymbol u}\right)=0$, we deduce that
	$$\begin{array}{ll}a_\star(B { \boldsymbol u},{ \boldsymbol v},{ \boldsymbol v})&=
	a(B { \boldsymbol u},{ \boldsymbol v},{ \boldsymbol v})+
	\delta\left(u_3v_2,v_3\right)-\delta\left(u_3,v_3v_2\right)=a(B { \boldsymbol u},{ \boldsymbol v},{ \boldsymbol v})=0.\end{array}$$
Similarly,
	$$\begin{array}{ll}a_\star(B { \boldsymbol u},{ \boldsymbol v},{ \boldsymbol w})&=a(B { \boldsymbol u},{ \boldsymbol v},{ \boldsymbol w})+\delta\left(u_3v_2,w_3\right)-\delta
	\left(u_3v_3,w_2\right)\vspace{2mm}\\
	&=-a(B { \boldsymbol u},{ \boldsymbol w},{ \boldsymbol v})+\delta\left(u_3v_2,w_3\right)-\delta
	\left(u_3v_3,w_2\right)=-a_\star(B { \boldsymbol u},{ \boldsymbol w},{ \boldsymbol v})\end{array}$$
and the proof is complete.$\hfill\Box$
\begin{lemma} \label{convective_2} Let ${ \boldsymbol u}$, ${ \boldsymbol v}$ and ${ \boldsymbol w}$ be in $\boldsymbol H^{1}_0(\Sigma)$. Then the following estimate holds
	$$\left|a_\star({ \boldsymbol u},{ \boldsymbol v},B{ \boldsymbol w})\right|\leq 
	\kappa_3
	\left\|D^\star{ \boldsymbol u}\right\|_{2,B}\left\|D^\star{ \boldsymbol v}\right\|_{2,B}
	\left\|D^\star{ \boldsymbol w}\right\|_{2,B},$$
where $\kappa_3=nm^{\frac{3}{2}} |\Sigma|^{\frac{3}{4}}$.
\end{lemma}
{\bf Proof.}  Setting $r=4$ and $q=2$ in the Sobolev inequality (\ref{sobolev0}), we obtain
	$$\begin{array}{ll}\left|a_\star({ \boldsymbol u},{ \boldsymbol v},B{ \boldsymbol w})\right|&\leq n
	\left\|{ \boldsymbol u}\right\|_{4}\left\|\nabla^\star { \boldsymbol v}\right\|_{2}
	\left\|{ \boldsymbol w}\right\|_{4}\leq n \left( S_{2,4}\right)^2
	 \left\|\nabla  { \boldsymbol u}\right\|_{2}  \left\|\nabla^\star { \boldsymbol v}\right\|_{2} 
	\left\|\nabla  { \boldsymbol w}\right\|_{2}\vspace{2mm}\\
	&\leq nm^{\frac{3}{2}}\left( S_{2,4}\right)^2
	 \left\|\nabla^\star  { \boldsymbol u}\right\|_{2,B}  \left\|\nabla^\star { \boldsymbol v}\right\|_{2,B} 
	\left\|\nabla^\star  { \boldsymbol w}\right\|_{2,B}. \end{array}$$
The estimate follows then by taking into account Remark \ref{remark-korn}. $\hfill\Box$
\subsection{Estimates on the extra stress tensor}
Our aim now is to establish some continuity, coercivity and 
monotonicity results for the stress tensor $\boldsymbol\tau$. 
\begin{proposition} \label{tensor_proper3} Assume that $p\geq 2$ and let  $\boldsymbol f, \boldsymbol g\in { \boldsymbol L}^p(\Sigma,\mathbb{R}^{3\times 3})$. Then the following estimates hold
	\begin{description}
 \item {Continuity.}
\begin{equation}
	\label{estim_S3_3}\Big\|\left(1+|\boldsymbol f|^2\right)^{\frac{p-2}{2}}\boldsymbol g\Big\|_{p',B}\leq  {\cal F}_B\big(\left\|\boldsymbol f\right\|_{p,B}
\big)\left\|\boldsymbol g\right\|_{p,B},\end{equation}
\begin{equation}
	\label{estim_S3_2}\left\|\boldsymbol\tau\left(\boldsymbol f\right)-\boldsymbol\tau\left(\boldsymbol g\right)\right\|_{p',B}\leq  (p-1){\cal F}_B\big(\left\|\boldsymbol f\right\|_{p,B}
+\left\|\boldsymbol g\right\|_{p,B}\big)
\left\|\boldsymbol f-\boldsymbol g\right\|_{p,B},\end{equation}

 \item {Coercivity.}
	\begin{equation}\label{estim_S2}\left({\boldsymbol\tau}(\boldsymbol f),B \boldsymbol f\right)\geq 2\left\|\boldsymbol f\right\|_{2,B}^2, \qquad \quad
	 \left({\boldsymbol\tau}(\boldsymbol f),B\boldsymbol f\right)\geq 2\left\|\boldsymbol f\right\|_{p,B}^p, \end{equation}
 \item {Monotonicity.}
\begin{equation}\label{estim_S3}\begin{array}{ll}\left({\boldsymbol\tau}(\boldsymbol f)-{\boldsymbol\tau}(\boldsymbol g),B\left(\boldsymbol f-\boldsymbol g\right)\right)
	\geq 2\left\|\boldsymbol f-\boldsymbol g\right\|_{2,B}^2,\vspace{1mm}\\
	\left({\boldsymbol\tau}(\boldsymbol f)-{\boldsymbol\tau}(\boldsymbol g),B\left(\boldsymbol f-\boldsymbol g\right)\right)\geq \tfrac{1}{2^{p-1}(p-1)}
	\left\|\boldsymbol f-\boldsymbol g\right\|_{p,B}^p.\end{array}\end{equation}
\end{description}
where ${\cal F}_B\big(\lambda\big)=\left(\big\|B^{\frac{1}{p}}\big\|_p+\lambda \right)^{p-2}$.
 \end{proposition} 
{\bf Proof.} Assume that $p>2$. Standard calculation show that
	$$\begin{array}{ll}\Big\|\left(1+|\boldsymbol f|^2\right)^{\frac{p-2}{2}}\boldsymbol g\Big\|_{p',B}&\leq \Big\|\left(1+|\boldsymbol f|^2\right)^{\frac{p-2}{2}}
\Big\|_{\frac{p}{p-2},B}\left\|\boldsymbol g\right\|_{p,B}\leq  \left(\big\|B^{\frac{1}{p}}\big\|_p+\left\|\boldsymbol f\right\|_{p,B}\right)^{p-2}\left\|\boldsymbol g\right\|_{p,B}\end{array}$$
which gives (\ref{estim_S3_3}). Estimates (\ref{estim_S2}) and $(\ref{estim_S3})$ are a direct consequence of the coercivity properties and the  monotonicity properties in Lemma   $\ref{tensor_proper2}$.  Finally, observing that
$$\begin{array}{ll}\left\|{\boldsymbol\tau}\left(\boldsymbol f\right)-{\boldsymbol\tau}\left(\boldsymbol g\right)\right\|_{p',B}
&\leq \displaystyle (p-1)\left\|\left(1+|\boldsymbol f|^2
+ |\boldsymbol g|^2\right)^{\frac{p-2}{2}}|\boldsymbol f-\boldsymbol g|
	\right\|_{p',B}\vspace{1mm}\\
&\leq \displaystyle (p-1)\left\|\left(1+  |\boldsymbol f|^2
+ |\boldsymbol g|^2\right)^{\frac{p-2}{2}}\right\|_{\frac{p}{p-2},B}
\|\boldsymbol f-\boldsymbol g\|_{p,B}\vspace{2mm}\\
&\leq  (p-1)\left(\big\|B^{\frac{1}{p}}\big\|_p+ \|\boldsymbol f\|_{p,B}
+\|\boldsymbol g\|_{p,B}\right)^{p-2}
\left\|\boldsymbol f-\boldsymbol g\right\|_{p,B}\end{array}$$
we obtain (\ref{estim_S3_2}). The case $p=2$ is direct. $\hfill\Box$ \vspace{2mm}\\
We finish the section by a result that will be useful in the sequel.
\begin{proposition} \label{prop_div_diff}Assume that $p\geq 2$ and let  $\boldsymbol f\in { \boldsymbol L}^p(\Sigma,\mathbb{R}^{3\times 3})$. Then the following estimate holds
	\begin{equation}\label{diff_div}\left\|\nabla^\star\cdot \boldsymbol\tau(\boldsymbol f)-\nabla\cdot \boldsymbol\tau(\boldsymbol f)\right\|_{p'}\leq 4\delta m\,{\cal F}_1\big(\|\boldsymbol f\|_p\big)
	\left(\left\|f\right\|_{p}
	+\left\|f_{33}\right\|_{p}+\left\|f_{23}\right\|_{p}\right)\end{equation}
where ${\cal F}_1$ is defined in Proposition $\ref{tensor_proper3}$ and 
$f=(f_{ij})_{i,j=1,2}$.
\end{proposition}
{\bf Proof.} Taking into account the definition of $\nabla^\star$, we obtain
	$$\begin{array}{ll}\left\|\nabla^\star\cdot \boldsymbol\tau(\boldsymbol f)-\nabla \cdot\boldsymbol\tau(\boldsymbol f)\right\|_{p'}&=\left\|\tfrac{\delta}{B} 
	\left({\tau}_{12}(\boldsymbol f)\,{\boldsymbol a }_1+
	\left({\tau}_{22}(\boldsymbol f)-{\tau}_{33}(\boldsymbol f)\right)\, {\boldsymbol a }_2+2{\tau}_{23}(\boldsymbol f)\,  {\boldsymbol a }_3
	\right)\right\|_{p'}\vspace{2mm}\\
	&\leq \delta m \left(\left\|{\tau}_{12}(\boldsymbol f)
	\right\|_{p'}+
	\left\|{\tau}_{22}(\boldsymbol f)\right\|_{p'}+
	\left\|{\tau}_{33}(\boldsymbol f)\right\|_{p'}+2
	\left\|{\tau}_{23}(\boldsymbol f)\right\|_{p'}
	\right)\vspace{2mm}\\
	&\leq 2\delta m \left(\left\|\tau(\boldsymbol f)\right\|_{p'}
	+\left\|{\tau}_{33}(\boldsymbol f)\right\|_{p'}
	+\left\|{\tau}_{23}(\boldsymbol f)\right\|_{p'}
	\right)
	\end{array}$$
and the conclusion follows from estimate  (\ref{estim_S3_3}).$\hfill \Box$
\subsection{Existence and uniqueness of shear-thickening flows}
The aim of this section is to prove the existence and uniqueness result Theorem \ref{main1}. As usual, we first derive some estimates that hold not only for the exact solution $\boldsymbol u$ of (\ref{weak_formulation}), but also for  corresponding standard Galerkin approximations $\boldsymbol u^k$. 
\begin{proposition} \label{estimate_2_alpha} Assume that $p\geq 2$ and let $\boldsymbol u$ be a weak solution of $(\ref{weak_formulation})$. Then, estimates
 $(\ref{main_estimates_1_thic})$-$(\ref{est_u_p_thic})$ hold.
\end{proposition}
{\bf Proof.} The proof is split into two steps.\vspace{1mm}\\
\underline{\it Step 1}. {\it Global estimates.}
Setting ${\boldsymbol\varphi}={ \boldsymbol u}$ in (\ref{weak_formulation}), and using Lemma \ref{prop_a1} and estimate (\ref{estim_S2}), we deduce that
	\begin{equation}\label{estimate1}2
	\left\|D^\star{ \boldsymbol u}\right\|_{2,B}^2 
	\leq \left({\boldsymbol\tau}(D^\star{ \boldsymbol u}),
	B D^\star{ \boldsymbol u}\right)=(G,u_3),\end{equation}
	\begin{equation}\label{estimate1_2}2
	\left\|D^\star{ \boldsymbol u}\right\|_{p,B}^p 
	\leq \left({\boldsymbol\tau}(D^\star{ \boldsymbol u}),
	B D^\star{ \boldsymbol u}\right)=(G,u_3).\end{equation}
Classical arguments together with  (\ref{sobolev0}) and  (\ref{korn1}) yield
	\begin{align}\left(G,u_3\right)
	&=\left|\left(\tfrac{G}{\sqrt{B}},
	\sqrt{B}\,u_3\right)\right|\nonumber\\
	&\leq |G|
	\left\|\tfrac{1}{\sqrt{B}}\right\|_2
	\left\|u_3\right\|_{2,B}= |G|
	\left\|B^{-1}\right\|_1^{\frac{1}{2}}
	\left\|u_3\right\|_{2,B}\nonumber\\
	&\leq \kappa_1\left\|\nabla
	  \left(\sqrt{B}\,u_3\right)\right\|_2\leq  \sqrt{2}\kappa_1
	\left(\left\|\nabla   u_3\right\|_{2,B}^2
	+\left\|\tfrac{\delta}{B}\,u_3\right\|_{2,B}^2
	\right)^{\frac{1}{2}}\nonumber\\
	&\label{estimate3_2}\leq 
	2\kappa_1\left\|D^\star{ \boldsymbol u}\right\|_{2,B},\end{align}
where $\kappa_1=\left(\tfrac{m}{2}\right)^{\frac{1}{2}}|G||\Sigma|$.
Due to (\ref{estimate1}) and (\ref{estimate3_2}), we obtain
(\ref{main_estimates_1_thic}) and \begin{equation}\label{estimate4}\left(G,u_3\right)\leq  2\kappa_1^2.\end{equation}
Estimate (\ref{main_estimates_2_thic}) follows then from (\ref{estimate1_2}) and (\ref{estimate4}). \vspace{1mm}\\
\underline{\it Step 2}. {\it Estimates for $u$ and $u_3$.} Let us now prove estimates (\ref{est_u3_thic})-(\ref{est_u_p_thic}). Notice first that (\ref{est_u3_thic}) is a direct consequence of $(\ref{main_estimates_1_thic})$ and (\ref{korn1}).
To derive the second estimate, we set $\varphi=u$ in the weak formulation $(\ref{weak_formulation_u})$ and get
	$$\left(\tau(D^\star{ \boldsymbol u}),B D u\right)+
	\delta\left({\tau}_{33}(D^\star{ \boldsymbol u}),u_2\right)=\delta{\cal R}e\,\left(u_3^2,u_2\right).$$
Therefore, by using the coercivity properties, we obtain
	\begin{equation}\label{est_u_1_2}\left\|Du\right\|_{2,B}^2+
	\left\|\tfrac{\delta}{B}u_2
	\right\|_{2,B}^2\leq \tfrac{\delta {\cal R}e }{2}\left(u_3^2,u_2\right),
	\end{equation}
	\begin{equation}\label{est_u_1_p} 
	\left\|Du\right\|_{p,B}^p+\left\|\tfrac{\delta}{B}u_2
	\right\|_{p,B}^p\leq \tfrac{\delta {\cal R}e}{2} \left(u_3^2,u_2\right).
	\end{equation}
On the other hand, taking into account  (\ref{est_psi_1}) with $\alpha=q=2$, (\ref{korn2}) and (\ref{est_u3_thic}), we have
	\begin{align}\label{u3_u2}
	\left|\left(u_3^2,u_2\right)\right|
	&\leq \left( S_{2,4}\right)^3 |\Sigma|^{\frac{1}{4}}
	\left\|\nabla u_3\right\|_{2}^2
	\left\|\nabla u_2\right\|_{2}\nonumber\\
	&\leq \tfrac{m^{\frac{3}{2}}|\Sigma|}{2\sqrt{2}}
	\left\|\nabla u_3\right\|_{2,B}^2
	\left\|\nabla u_2\right\|_{2,B}\nonumber\\
	&\leq \tfrac{m^{\frac{3}{2}}|\Sigma|}{2}
	  \left\|\nabla u_3\right\|_{2,B}^2
	\left\|D u\right\|_{2,B}\nonumber\\
	&\leq m^{\frac{3}{2}}|\Sigma|\,\kappa_1^2
	\left\|Du\right\|_{2,B}.\end{align}
Combining (\ref{est_u_1_2}) and (\ref{u3_u2}), we obtain
	 (\ref{est_u_2_thic}) and
	$$\left|\left(u_3^2,u_2\right)\right|\leq 2
	\left(\tfrac{m^{\frac{3}{2}}|\Sigma|\,\kappa_1^2}{2}
	 \right)^2\delta{\cal R}e.$$
Estimate (\ref{est_u_p_thic}) follows then from (\ref{est_u_1_p}).$\hfill\Box$\vspace{2mm}\\
\noindent{\bf Proof of Theorem \ref{main1}.} The proof, based on classical compactness and monotonicity arguments, is split into two steps.\vspace{2mm}\\
\underline{\it Step 1.} Let us prove the existence of a weak solution for (\ref{equation}). Let ${ \boldsymbol u}^k$ be a classical Galerkin approximation. Arguments similar to those used in the proof of Proposition \ref{estimate_2_alpha} show that
	$$\left\|D^\star{ \boldsymbol u}^k\right\|_{2,B}
	\leq \kappa_1, \qquad\quad
	\left\|D^\star{ \boldsymbol u}^k\right\|_{p,B}^p
	\leq\kappa_1^2$$
and imply that the sequences $(\sqrt{B}D^\star{ \boldsymbol u}^k)_k$ and
$(B^{\frac{1}{p}}D^\star{ \boldsymbol u}^k)_k$ are bounded in $\boldsymbol L^2(\Sigma)$ and $\boldsymbol L^p(\Sigma)$, respectively. By taking into account (\ref{grad_prime}), we deduce that $(\sqrt{B}\nabla{ \boldsymbol u}^k)_k$ is bounded in $\boldsymbol L^2(\Sigma)$ and thus  $(\nabla{ \boldsymbol u}^k)_k$ is also bounded in $\boldsymbol L^2(\Sigma)$. Moreover, 
due to estimate  (\ref{estim_S3_3}) we have
	$$\left\|{\boldsymbol\tau}\left(D^\star{ \boldsymbol u}^k\right)\right\|_{p',B}\leq 2{\cal F}_B\left(\left\|D^\star{ \boldsymbol u}^k\right\|_{p,B}\right)
	\left\|D^\star{ \boldsymbol u}^k\right\|_{p,B}$$
and thus $\left(B^{\frac{1}{p'}}{\boldsymbol\tau}\left(D^\star{ \boldsymbol u}^k\right)\right)_k$ is bounded in $\boldsymbol L^{p'}(\Sigma)$.
There then exist a subsequence, still indexed by $k$,  ${ \boldsymbol u}\in \boldsymbol V_B^p$ and $\widetilde {\boldsymbol\tau}\in \boldsymbol L^{p'}(\Sigma)$ such that 
	$$\nabla  { \boldsymbol u}^k\longrightarrow \nabla  { \boldsymbol u} \qquad \mbox{weakly} \ in \
	\boldsymbol L^2(\Sigma),$$
	$$B^{\frac{1}{p'}}{\boldsymbol\tau}\left(D^\star{ \boldsymbol u}^k\right)\longrightarrow \widetilde {\boldsymbol\tau} \qquad \mbox{weakly} \ in \
	\boldsymbol L^{p'}(\Sigma).$$
By using compactness results on Sobolev spaces,  we deduce that  $({ \boldsymbol u}^k)_k$  strongly converges
 to ${ \boldsymbol u}$ in $L^{4}(\Sigma)$ and thus
	$$\sqrt{B}D^\star{ \boldsymbol u}^k\longrightarrow \sqrt{B}D^\star{ \boldsymbol u} \qquad \mbox{weakly} \ in \
	\boldsymbol L^2(\Sigma).$$
Therefore, by taking into account Lemma \ref{prop_a1}, for all ${\boldsymbol\varphi}\in \boldsymbol V_B^p$ we have
        $$\begin{array}{ll}\left|a_\star\left({ \boldsymbol u}^k, 
	{ \boldsymbol u}^k,B{\boldsymbol\varphi}\right)
	-a_\star\left({ \boldsymbol u}, 
	{ \boldsymbol u},B{\boldsymbol\varphi}\right)\right|
	&\leq\left|a_\star\left({ \boldsymbol u}^k-{ \boldsymbol u},{ \boldsymbol u}^k,B{\boldsymbol\varphi}\right)\right|+
	\left|a_\star\left({ \boldsymbol u},{ \boldsymbol u}^k-{ \boldsymbol u},
	B{\boldsymbol\varphi}\right)\right|\vspace{1mm}\\
         &=\left|a_\star\left({ \boldsymbol u}^k-{ \boldsymbol u},{ \boldsymbol u}^k,B{\boldsymbol\varphi}\right)\right|
	+\left|a_\star\left(B{ \boldsymbol u},{\boldsymbol\varphi},{ \boldsymbol u}^k-{ \boldsymbol u}
	\right)\right|\vspace{1mm}\\
         & \leq \left(\left\|\nabla^\star { \boldsymbol u}^k\right\|_{2,B}
	\left\|{\boldsymbol\varphi}\right\|_{4,B}+
	\left\|{ \boldsymbol u}\right\|_{4,B} \left\|\nabla^\star 
	{\boldsymbol\varphi}\right\|_{2,B}\right) \left\|
	{ \boldsymbol u}^k-{ \boldsymbol u}\right\|_{4,B}\nonumber\\
          &\longrightarrow 0 \qquad \mbox{when}  \ k\rightarrow +\infty.\end{array}$$
By passing to the limit in
	$$\left({\boldsymbol\tau}(D^\star{ \boldsymbol u}^k),B D^\star{\boldsymbol\varphi}\right)+{\cal R}e\,a_\star\left({ \boldsymbol u}^k, { \boldsymbol u}^k,
	B{\boldsymbol\varphi}\right)=\left(G,\varphi_3\right)\qquad  \mbox{for all} \ {\boldsymbol\varphi}\in \boldsymbol V_B^{p},$$
we obtain
\begin{equation}\label{state_limit_weak3}\left(\widetilde {\boldsymbol\tau},B^{\frac{1}{p}} D^\star{\boldsymbol\varphi}\right)+{\cal R}e\, 
	a_\star\left({ \boldsymbol u},{ \boldsymbol u},
	B{\boldsymbol\varphi}\right)=\left(G,\varphi_3\right) \qquad \mbox{for all} \ {\boldsymbol\varphi}\in \boldsymbol V_B^{p}.\end{equation}
In particular, by settin $\boldsymbol\varphi=\boldsymbol u$ and using Lemma \ref{prop_a1} we deduce that
     \begin{equation}\label{8} \left(\widetilde{\boldsymbol\tau},B^{\frac{1}{p}} D^\star{ \boldsymbol u}\right)=
	\left(G,u_3\right).
  \end{equation}
On the other hand, $(\ref{estim_S3})_1$ implies 
                \begin{equation}\label{9}\left({\boldsymbol\tau}\left(
	D^\star{ \boldsymbol u}^k\right)-{\boldsymbol\tau}\left(D^\star{\boldsymbol\varphi}\right),B D^\star\left( { \boldsymbol u}^k-{\boldsymbol\varphi}\right)\right)\geq 0 
	\qquad \mbox{for all} \ 
                  {\boldsymbol\varphi}\in \boldsymbol V_B^{p}.\end{equation}
Since 
         $$\left({\boldsymbol\tau}\left(
	D^\star{ \boldsymbol u}^k\right),B D^\star{ \boldsymbol u}^k\right)=
	\left(G,u^k_3\right),$$ 
 by substituing in (\ref{9}), we obtain
	$$\left(G, u ^k_3\right)-\left({\boldsymbol\tau}\left(
	D^\star{ \boldsymbol u}^k\right),B D^\star{\boldsymbol\varphi}\right)-\left({\boldsymbol\tau}\left(D^\star{\boldsymbol\varphi}\right),B D^\star\left( { \boldsymbol u}^k-
	{\boldsymbol\varphi}\right)\right)\geq 0.$$
By passing to the limit, it follows that
          $$\left(G,u_3\right)-\left(B^{\frac{1}{p}}\,\widetilde{\boldsymbol\tau},D^\star{\boldsymbol\varphi}\right)-\left(B \,{\boldsymbol\tau}\left(D^\star{\boldsymbol\varphi}\right),D^\star \left({ \boldsymbol u}-{\boldsymbol\varphi}\right)\right)\geq 0\qquad \mbox{for all} \ 
                  {\boldsymbol\varphi}\in \boldsymbol V_B^{p}.$$
This inequality together with (\ref{8}) implies  that
         $$\left(B^{\frac{1}{p}}\,\widetilde{\boldsymbol\tau}
	-B \,{\boldsymbol\tau}\left(D^\star{\boldsymbol\varphi}\right)
                ,D^\star\left({ \boldsymbol u}-{\boldsymbol\varphi}\right)\right)\geq 0
            \qquad  \mbox{for all} \ \boldsymbol\varphi\in \boldsymbol V_B^{p}$$
and by setting ${\boldsymbol\varphi}={ \boldsymbol u}-t{\boldsymbol\psi}$ with $t>0$, we obtain
	$$\left(B^{\frac{1}{p}}\,\widetilde{\boldsymbol\tau}-B \,{\boldsymbol\tau}\left(D^\star{ \boldsymbol u}-tD{\boldsymbol\psi}\right)
                ,D^\star{\boldsymbol\psi}\right)
	\geq 0\qquad \mbox{for all} \ 
                  {\boldsymbol\psi}\in \boldsymbol V_B^{p}.$$ 
Letting $t$ tend to zero and using the continuity of ${\boldsymbol\tau}$, we get
	$$\left(B^{\frac{1}{p}}\,\widetilde{\boldsymbol\tau}-
	B\,{\boldsymbol\tau}\left(D^\star{ \boldsymbol u}\right),D^\star{\boldsymbol\psi}\right)\geq 0
            \qquad  \mbox{for all} \ {\boldsymbol\psi}
	\in \boldsymbol V_B^{p}$$ 
 and thus 
	\begin{equation}\label{tensor-beta}\left(B^{\frac{1}{p}}\,\widetilde{\boldsymbol\tau},D^\star{\boldsymbol\psi}\right)=\left(B \, {\boldsymbol\tau}\left(D^\star{ \boldsymbol u}\right),D^\star{\boldsymbol\psi}\right)
            \qquad  \mbox{for all} \ {\boldsymbol\psi}\in \boldsymbol V_B^{p}.
           \end{equation}
Combining (\ref{state_limit_weak3}) and (\ref{tensor-beta}), we deduce that
	$$\left(B\, {\boldsymbol\tau}\left(D^\star{ \boldsymbol u}\right),D^\star{\boldsymbol\varphi}\right)
	+{\cal R}e\, a_\star\left({ \boldsymbol u},{ \boldsymbol u},
	B{\boldsymbol\varphi}\right)=\left(G,\varphi_3\right)\qquad  \mbox{for all} \ {\boldsymbol\varphi}\in \boldsymbol V_B^{p}.$$
Hence ${ \boldsymbol u}$ is a solution of (\ref{weak_formulation}) \vspace{2mm}\\
\underline{\it Step 2.}  To prove the uniqueness result, let us assume that ${ \boldsymbol u}$ and ${ \boldsymbol v}$ are two weak solutions of $(\ref{equation})$. Substituting in the weak formulation of (\ref{weak_formulation}), setting $\boldsymbol\varphi={ \boldsymbol u}-{ \boldsymbol v}$ and taking into account Lemma \ref{prop_a1}, Lemma \ref{convective_2} and Proposition \ref{estimate_2_alpha}, we obtain
	\begin{align}\label{visco_convective3}\tfrac{1}{{\cal R}e}\left({\boldsymbol\tau}(D^\star{ \boldsymbol u})-
	{\boldsymbol\tau}(D^\star{ \boldsymbol v}),BD^\star\left({ \boldsymbol u}-{ \boldsymbol v}\right)\right)&=-a_\star( { \boldsymbol u}, { \boldsymbol u},B\left({ \boldsymbol u}-{ \boldsymbol v}\right))+
	a_\star\left( { \boldsymbol v}, { \boldsymbol v},B\left({ \boldsymbol u}-{ \boldsymbol v}\right)\right)\nonumber\\
	&=-a_\star\left( { \boldsymbol u},{ \boldsymbol u}-{ \boldsymbol v},B\left({ \boldsymbol u}-{ \boldsymbol v}\right)\right)
	-a_\star\left({ \boldsymbol u}-{ \boldsymbol v}, { \boldsymbol v},B\left({ \boldsymbol u}-{ \boldsymbol v}\right)\right)\nonumber\\
	&=-a_\star\left({ \boldsymbol u}-{ \boldsymbol v}, { \boldsymbol v},B\left({ \boldsymbol u}-{ \boldsymbol v}\right)\right)\nonumber\\
	&\leq  \kappa_3 
	\left\|D^\star\left({ \boldsymbol u}-{ \boldsymbol v}\right)\right\|_{2,B}^2
	\left\|D^\star  { \boldsymbol v}\right\|_{2,B}\nonumber\\
	&\leq \kappa_1\kappa_3
	\left\|D^\star\left({ \boldsymbol u}-{ \boldsymbol v}\right)\right\|_{2,B}^2. \end{align}
Combining (\ref{visco_convective3}) and $(\ref{estim_S2})_1$, we deduce that
	$$\left(2-\kappa_1\kappa_3\,{\cal R}e\right)\left\|D^\star\left({ \boldsymbol u}-{ \boldsymbol v}\right)\right\|_{2,B}^2\leq 0$$
and thus $ { \boldsymbol u}\equiv  { \boldsymbol v}$ if 
	${\cal R}e<\tfrac{2}{\kappa_1\kappa_3}$.$\hfill\Box$\vspace{2mm}\\
{\bf Proof of Corollary \ref{pressure_1}.} To simplify the redaction, let us set $\boldsymbol\tau(D^\star{ \boldsymbol u})=\boldsymbol\tau$. Notice first that if ${ \boldsymbol u}$ is a weak solution of (\ref{equation}) then
\begin{align}\label{weak_form_0}\left(\boldsymbol\tau,D \boldsymbol\varphi\right) -\left(\nabla^\star\cdot \boldsymbol\tau-\nabla\cdot \boldsymbol\tau,\boldsymbol\varphi\right)
	+{\cal R}e\,a_\star\left(\boldsymbol u, 
	  \boldsymbol u,\boldsymbol\varphi\right)=(\tfrac{G}{B}
	,\varphi_3)\qquad \mbox{for all} \ \boldsymbol\varphi\in \boldsymbol V^p.\end{align}
It follows that $u=(u_1,u_2,0)$ satisfies
$$\left(\tau,D \varphi\right)-
\left(\nabla^\star\cdot \tau-\nabla\cdot \tau,\varphi\right)+{\cal R}e\left(\left(u\cdot \nabla   u,\varphi\right)-\left(\tfrac{\delta}{B} u_3^2,\varphi_2\right)\right)=0$$
for all $\varphi=(\varphi_1,\varphi_2)\in \boldsymbol W^{1,p}_0(\Sigma,\mathbb R^2)$ such that $\nabla\cdot \varphi=0$, with $\tau=\left({\tau}_{ij}\right)_{i,j=1,2}$.
 Taking into account (\ref{estim_S3_3}), (\ref{diff_div}) and using standard arguments, we can prove that the mapping
	$${\cal G}: \varphi \mapsto \left(\tau,D \varphi\right)
-
\left(\nabla^\star\cdot \tau-\nabla\cdot \tau,\varphi\right)+{\cal R}e\left(\left(u\cdot \nabla   u,\varphi\right)
-\left(\tfrac{\delta}{B} u_3^2,\varphi_2\right)\right)$$
is a linear continuous functional on $\boldsymbol W^{1,p}_0(\Sigma,\mathbb{R}^2)$. By using a classical result (see \cite{amrouche}), we deduce that there exists $\pi\in \boldsymbol L^{p'}_0(\Sigma)$ such that
	$${\cal G}(\varphi)=-\left(\nabla  \pi,\varphi\right)=
	\left(\pi,\nabla  \cdot \varphi\right) \qquad \mbox{for all} \  \varphi\in \boldsymbol W^{1,p}_0(\Sigma,\mathbb{R}^2).$$
Moreover, there exists a positive constant $C$ depending only on $p$ and $\Sigma$ such that
	\begin{equation}\label{pression1}
	C\left\|\pi\right\|_{p'}\leq \left\|\nabla  \pi\right\|_{-1,p'}=\sup_{\varphi\in  \boldsymbol W^{1,p}_0(\Sigma,\mathbb{R}^2)}
	\tfrac{\left|\left(\pi,\nabla  \cdot \varphi\right)\right|}{\left\|\nabla   \varphi\right\|_p}.\end{equation}
On the other hand, using (\ref{sobolev0}), (\ref{estim_S3_3}), (\ref{diff_div}) and 
(\ref{main_estimates_2_thic}) we obtain
	\begin{align}\label{pression2}\left|\left(\tau,D \varphi\right)
-\left(\nabla^\star\cdot \tau-\nabla\cdot \tau,\varphi\right)\right|
		&\leq \left\|\tau
	\right\|_{p'}\left\|D \varphi\right\|_p +\left\|\nabla^\star\cdot \tau-\nabla\cdot \tau\right\|_{p'}
	\left\|\varphi\right\|_p\nonumber\\
		&\leq  \left(
	\left\|\tau\right\|_{p'}+S_{p,p}\left\|\nabla^\star\cdot \tau-\nabla\cdot \tau\right\|_{p'}
	\right)\left\|\nabla \varphi\right\|_p\nonumber\\
		&\leq  {\cal F}_1\left(\|D^\star{ \boldsymbol u}\|_p\right)\left(\left(1+4S_{p,p}\delta m\right)\left\|Du\right\|_{p}+
	4S_{p,p}\delta m \, \left\|\tfrac{\delta}{B}u_2\right\|_{p}\right)\left\|\nabla \varphi\right\|_p\nonumber\\
		&\leq  \tilde \kappa\left(
	\left\|Du\right\|_{p,B}+\delta^2 
	  \left\|u_2\right\|_{p,B}
	\right)\left\|\nabla \varphi\right\|_p,\end{align} 
where $\tilde \kappa$ only depends on $p$, $\Sigma$, $m$ and $n$. Similarly, we can easily see that
	 \begin{align}\label{pression3}\left|\left(u\cdot \nabla   u,\varphi\right)-\left(\tfrac{\delta}{B} u_3^2,\varphi_2\right)\right|&=
	\left|-\left(u\otimes u,\nabla\varphi\right)-\left(\tfrac{\delta}{B} u_3^2,\varphi_2\right)\right|\nonumber\\
	&\leq \left\|u\otimes u\right\|_{p'}\|\nabla\varphi\|_p+
	\delta m    \left\|u_3^2\right\|_{p'}
	\left\|\varphi_2\right\|_p\nonumber\\
	&\leq \hat\kappa \left(\left\|D u\right\|_{p,B}^2+
	\delta \left\|\nabla u_3\right\|_{p,B}^2\right)\left\|\nabla\varphi\right\|_p,	\end{align}
where $\hat \kappa$ only depends on $p$, $\Sigma$ and $m$.
Combining (\ref{pression1})-(\ref{pression3}), we deduce that
	$$\left\|\pi\right\|_{p'}\leq \kappa \left(\left\|Du\right\|_{p,B}+\delta^2   \left\|u_2\right\|_{p,B}+{\cal R}e \left(\left\|D u\right\|_{p,B}^2+
	\delta \left\|\nabla u_3\right\|_{p,B}^2\right)\right),$$
where $\kappa$ is a positive constant only depending on $\Sigma$, $p$, $m$ and $n$.\hfill\Box
\subsection{$\delta$-approximation}
{\bf Proof of Proposition \ref{w_exist_thick}.} Based on a standard Galerkin approximation of the corresponding {\it global} formulation, compactness and monotonicity arguments, the existence of a weak solution for $(E_\sigma)$ can be established once suitable a priori estimates are derived. However, because of the term involving $\sigma(w_3)$, the {\it global} formulation 
does not seem appropriate unless we restrict strongly the exponent $\alpha$ in (\ref{psi_assumption}). To overcome this difficulty, we consider the {\it coupled} formulations. Arguing as in the proof of Proposition \ref{estimate_2_alpha}, by setting $\boldsymbol\varphi=(0,0,w_3)$ in the corresponding weak formulation we obtain
	$$\int_\Sigma\left(1+|D\boldsymbol w|^2\right)^{\frac{p-2}{2}}
	|\nabla w_3|^2\,dx=\left(\tfrac{G}{B},w_3\right)$$
and thus
	$$\left\|\nabla w_3\right\|_2^2\leq 
	\left(\tfrac{G}{B},w_3\right)\leq c_1
	\|\nabla w_3\|_2 \qquad \mbox{with} \ c_1=
	\tfrac{m|G||\Sigma|}{\sqrt{2}} .$$
Hence
	  $$\left\|\nabla w_3\right\|_2\leq c_1, \qquad 
	\left(\tfrac{G}{B},w_3\right)\leq c_1^2$$
and 
	$$\left(\tfrac{1}{2}\right)^{\frac{p-2}{2}}
	\left\|\nabla w_3\right\|_p^p\leq 
	\left(\tfrac{G}{B},w_3\right)\leq c_1^2.$$  
Similarly, by setting $\boldsymbol\varphi=(w,0)=(w_1,w_2,0)$ in the corresponding weak formulation, we obtain
	$$2\int_\Sigma\left(1+|D\boldsymbol w|^2\right)^{\frac{p-2}{2}}
	|D w|^2\,dx=\delta \left(\sigma(w_3),w_2\right)$$
and thus
	$$\|Dw\|_{2}^2\leq \tfrac{\delta}{2}\left(\sigma(w_3),w_2\right),\qquad \|Dw\|_{p}^p\leq \tfrac{\delta}{2}\left(\sigma(w_3),w_2\right).$$
Estimate (\ref{est_psi_1}) together with the Korn inequality yield
	$$\|Dw\|_{2}^2\leq  \tfrac{c_0D_{2,\alpha}}{2}\, \delta
	\left\|\nabla w_3\right\|_2^\alpha
	\left\|\nabla w_2\right\|_2\leq\tfrac{c_0D_{2,\alpha}}{2}\,\delta
	\left\|\nabla w_3\right\|_2^\alpha
	\left\|\nabla w\right\|_2\leq c_0c_2 \delta
	\left\|Dw\right\|_2$$
and consequently
	$$\|Dw\|_{2}\leq c_0c_2 \delta, \qquad 
	\left(\sigma(w_3),w_2\right)\leq 
	2\delta\left(c_0c_2\right)^2.$$
Therefore
	$$\|Dw\|_{p}^p\leq \left(c_0c_2\delta\right)^2$$
and the a priori estimates are derived. The proof may be completed using arguments similar to those in the proof of Theorem \ref{main1}.   $\hfill\Box$
\begin{remark} Arguing as in the proof of Corollary $\ref{pressure_1}$ and using $(\ref{est_psi_1})$, we can prove the existence of $\pi\in L^{p'}(\Sigma)$ such that $(E_\sigma)_1$ holds in $\boldsymbol W^{-1,p'}(\Sigma)$. Moreover, the following estimate holds
	$$\left\|\pi\right\|_{p'}\leq \kappa\left(
	\|Dw\|_{p}+
	{\cal R}e\left\|Dw\right\|^2_{p}+
	 \delta\left\|\nabla w_3
	\right\|_p^\alpha\right),$$
where $\kappa$ is a positive constant only depending on $p$, $\Sigma$, $m$, $n$ and $\alpha$.
\end{remark}
{\bf Proof of Proposition \ref{delta_vs_zero}.} Let us first recall that $\boldsymbol u$ satisfies (\ref{weak_form_0}) and that $\boldsymbol w$ satisfies 
	$$\left({\boldsymbol\tau}\left(D  { \boldsymbol w}\right),D
	  {\boldsymbol\varphi}\right)
	+{\cal R}e\,a\left({ \boldsymbol w},{ \boldsymbol w},{\boldsymbol\varphi}\right)-\left(\pi_2,\nabla  \cdot{\boldsymbol\varphi}\right)
	=\left(\tfrac{G}{B},\varphi_3\right)+\delta\left(\sigma(w_3),\varphi_2\right)
	$$
for all ${\boldsymbol\varphi}\in \boldsymbol 
	W^{1,p}_0(\Sigma)$. Therefore
\begin{align}\label{difference1}\left(\boldsymbol\tau(D{ \boldsymbol u})-{\boldsymbol\tau}\left(D  { \boldsymbol w}\right),D \boldsymbol\varphi\right)&=\left(\boldsymbol\tau(D{ \boldsymbol u})-{\boldsymbol\tau}\left(D^\star{ \boldsymbol u}\right),D \boldsymbol\varphi\right)+
\left(\nabla^\star\cdot \boldsymbol\tau(D^\star{ \boldsymbol u})-\nabla\cdot \boldsymbol\tau(D^\star{ \boldsymbol u}),\boldsymbol\varphi\right)\nonumber\\
	&-{\cal R}e\left(a_\star\left(\boldsymbol u, 
	  \boldsymbol u,\boldsymbol\varphi\right)-a\left({ \boldsymbol w},{ \boldsymbol w},{\boldsymbol\varphi}\right)\right)+\left(\pi_1-\pi_2,\nabla\cdot \boldsymbol\varphi\right)-
	\delta\left(\sigma(w_3),\varphi_2\right)\nonumber\\
	&=I_1+I_2+I_3+I_4+I_5.\end{align}
$\circ$ Let us estimate the first term. By taking into account (\ref{estim_S3_2}), we have
$$\begin{array}{ll}\left|I_1\right|&=\left|\left({\boldsymbol\tau}\left(D^\star{ \boldsymbol u}\right)-{\boldsymbol\tau}\left(D  { \boldsymbol u}\right),D\boldsymbol\varphi\right)\right|\leq \left\|{\boldsymbol\tau}\left(D^\star{ \boldsymbol u}\right)-{\boldsymbol\tau}\left(D  { \boldsymbol u}\right)\right\|_{p'}
\left\|D\boldsymbol\varphi\right\|_{p}\vspace{1mm}\\	
&\leq \displaystyle (p-1) {\cal F}_1\left(\|D^\star{ \boldsymbol u}\|_{p}+\|D{ \boldsymbol u}\|_{p}\right)
\left\|D^\star{ \boldsymbol u}-D  { \boldsymbol u}\right\|_{p}
	\left\|D\boldsymbol\varphi\right\|_{p}.
	\end{array}$$
Since
	\begin{equation}\label{w_wstar}\left\|D^\star{ \boldsymbol u}-D  { \boldsymbol u}\right\|_{p}\leq \delta m \left(\|u_2\|_p+\|u_3\|_p\right),\end{equation}
 we deduce that
	\begin{equation}\label{difference_est3}\left|I_1\right|\leq F_ 1 \,\delta \left\|D\boldsymbol\varphi\right\|_{p},\end{equation}
where $F_1=m(p-1) {\cal F}_1\big(\|D^\star{ \boldsymbol u}\|_{p}+\|D{ \boldsymbol u}\|_{p}\big)\left(\|u_2\|_p+\|u_3\|_p\right)$. \vspace{2mm}\\
$\circ$ Estimate (\ref{diff_div}) together with the  Sobolev inequality (\ref{sobolev0}) and the Korn inequality (\ref{korn-l2}) yield 
	\begin{align}\label{difference_est4}\left|I_2\right|&\leq 
	\left\|\nabla^\star\cdot \boldsymbol\tau(D^\star{ \boldsymbol u})-\nabla\cdot \boldsymbol\tau(D^\star{ \boldsymbol u})\right\|_{p'}
	\left\|\boldsymbol\varphi\right\|_p\nonumber\\
	&\leq 4S_{2,p}m{\cal F}_1\big(\|D^\star{ \boldsymbol u}\|_p\big)
	\left(\left\|Du\right\|_{p}
	+\left\|D_{33}^\star{ \boldsymbol u}\right\|_{p}+\left\|D_{23}^\star{ \boldsymbol u}\right\|_{p}\right)\, \delta 
	\left\|\nabla \boldsymbol\varphi\right\|_{2}\leq  F_2\, \delta
	\left\|D\boldsymbol\varphi\right\|_{p},\end{align}
where $F_2=4\sqrt{6}m|\Sigma|^{\frac{1}{2}-\frac{1}{p}}\, S_{2,p}\, {\cal F}_1\big(\|D^\star{ \boldsymbol u}\|_p\big)\|D^\star{ \boldsymbol u}\|_p$. \vspace{2mm}\\
$\circ$ Similarly,
	\begin{align}\label{difference_est61}\tfrac{1}{{\cal R}e}\left|I_3\right|
&=\left|a_\star({ \boldsymbol u},{ \boldsymbol u},\boldsymbol \varphi)-a({ \boldsymbol w},{ \boldsymbol w},\boldsymbol \varphi)\right|
	\nonumber\\
&=\left|a_\star\left({ \boldsymbol u}-{ \boldsymbol w},{ \boldsymbol u},\boldsymbol \varphi\right)+a\left(\boldsymbol w,{ \boldsymbol u}-{ \boldsymbol w},\boldsymbol \varphi\right)+
	\left({ \boldsymbol w}\cdot \nabla^\star{ \boldsymbol u}-{ \boldsymbol w}\cdot \nabla  { \boldsymbol u},\boldsymbol \varphi\right)\right|\nonumber\\
&\leq\left|a_\star\left({ \boldsymbol u}-{ \boldsymbol w},{ \boldsymbol u},\boldsymbol \varphi\right)\right|
+\left|a\left(\boldsymbol w,{ \boldsymbol u}-{ \boldsymbol w},\boldsymbol \varphi\right)\right|+\left|
	\left({ \boldsymbol w}\cdot \nabla^\star{ \boldsymbol u}-{ \boldsymbol w}\cdot \nabla  { \boldsymbol u},\boldsymbol \varphi\right)
\right|\nonumber\\
&=\left|a_\star\left({ \boldsymbol u}-{ \boldsymbol w},{ \boldsymbol u},\boldsymbol \varphi\right)\right|+\left|a\left(\boldsymbol w,{ \boldsymbol u}-{ \boldsymbol w},\boldsymbol \varphi\right)\right|+
\left|\left(\tfrac{\delta}{B}\,w_2u_2,\varphi_3\right)-\left(\tfrac{\delta}{B}\,w_3u_3,\varphi_2\right)\right|\nonumber\\
&\leq\left|a_\star\left({ \boldsymbol u}-{ \boldsymbol w},{ \boldsymbol u},\boldsymbol \varphi\right)\right|+\left|a\left(\boldsymbol w,{ \boldsymbol u}-{ \boldsymbol w},\boldsymbol \varphi\right)\right|+
\displaystyle\left\|\tfrac{\delta}{B}|\boldsymbol w||\boldsymbol u||\boldsymbol \varphi|\right\|_1\nonumber\\
&\leq 	\left\|{ \boldsymbol u}-{ \boldsymbol w}\right\|_4
	\left\|\nabla^\star{ \boldsymbol u}\right\|_2\left\|\boldsymbol \varphi\right\|_4 +\left\|{ \boldsymbol w}\right\|_4
	\left\|\nabla\left({ \boldsymbol u}-{ \boldsymbol w}\right)\right\|_2\left\|\boldsymbol \varphi\right\|_4 + \delta m \left\|\boldsymbol w\right\|_4 \left\|\boldsymbol u\right\|_4
	\left\|\boldsymbol \varphi\right\|_2\nonumber\\
&\leq  \left(1+S_{2,2}\right)\left( S_{2,4}\right)^2
\left(\left\|\nabla\left({ \boldsymbol u}-{ \boldsymbol w}\right)\right\|_2\left(\left\|\nabla^\star{ \boldsymbol u}\right\|_2+\left\|\nabla{ \boldsymbol w}\right\|_2\right)+\delta m\left\|\nabla{ \boldsymbol w}\right\|_2\left\|\nabla{ \boldsymbol u}\right\|_2\right)
	\left\|\nabla \boldsymbol \varphi\right\|_2\nonumber\\
&\leq
F_3\left(\left\|D\left({ \boldsymbol u}-{ \boldsymbol w}\right)\right\|_{2}\left\|D\boldsymbol\varphi\right\|_2+\delta
	\left\|D\boldsymbol\varphi\right\|_{p}\right)\end{align}
with $F_3=\left(1+m\right)\sqrt{8} \left(1+S_{2,2}\right)\left(S_{2,4}\right)^2\left(\left\|D^\star{ \boldsymbol u}\right\|_{2}+\left\|D{ \boldsymbol w}\right\|_{2}+\left\|D\boldsymbol u\right\|_{2}\left\|D{ \boldsymbol w}\right\|_{2}\right)$.\vspace{2mm}\\
$\circ$ Let us now consider the term involving the pressure
	$$\left|I_4\right|=\left|\left(\pi_1-\pi_2,\nabla  \cdot
\boldsymbol\varphi\right)\right|\leq \left\|\pi_1-\pi_2\right\|_{p'}\left\|\nabla\cdot \boldsymbol\varphi\right\|_p.$$
Arguing as in the first part of the proof of Corollary \ref{pressure_1}, we can see that
	 $$\begin{array}{ll}\left\|\pi_1-\pi_2\right\|_{p'}
	&\leq   \tilde \kappa \left(\left\|\tau(D^\star{ \boldsymbol u})-\tau(D{ \boldsymbol w})\right\|_{p'}+\left\|\nabla^\star\cdot \tau(D^\star{ \boldsymbol u})-\nabla\cdot \tau(D^\star{ \boldsymbol u})\right\|_{p'}\right)\vspace{2mm}\\
	 &+\tilde \kappa {\cal R}e
	 \left(\left\|u\otimes u-w\otimes w\right\|_{p'}+
	\delta   \left\|u_3^2\right\|_{p'}\right)+\tilde \kappa \delta 
	\left\|\sigma(w_3)\right\|_{p'},\end{array}$$
where $\tilde \kappa$ depends only on $\Sigma$, $p$, $m$ and $n$.
Taking into account  (\ref{estim_S3_2}) and (\ref{w_wstar}), we have
	$$\begin{array}{ll}\left\|\tau(D^\star{ \boldsymbol u})-\tau(D{ \boldsymbol w})
	\right\|_{p'}&\leq 
	\left\|\boldsymbol\tau(D^\star{ \boldsymbol u})-\boldsymbol\tau(D{ \boldsymbol w})\right\|_{p'}\leq F_{4,1}\left\|D^\star{ \boldsymbol u}-D{ \boldsymbol w}\right\|_{p}\vspace{2mm}\\
	&\leq F_{4,1}\left(\left\|D^\star{ \boldsymbol u}-D{ \boldsymbol u}\right\|_{p}
	+\left\|D{ \boldsymbol u}-D{ \boldsymbol w}\right\|_{p}\right)\vspace{2mm}\\
	&\leq F_{4,2}\left(\delta+\left\|D{ \boldsymbol u}-D{ \boldsymbol w}\right\|_{p}
	\right),\end{array}$$
where $F_{4,1}=(p-1){\cal F}_1\big(\|D^\star{ \boldsymbol u}\|_{p}
+\|D  { \boldsymbol w}\|_{p}\big)$ and $F_{4,2}=F_{4,1}\left(1+ m\left(\|u_2\|_p+\|u_3\|_p\right)\right)$.
Moreover, by using (\ref{diff_div}), we obtain
$$\left\|\nabla^\star\cdot \tau(D^\star{ \boldsymbol u})-\nabla\cdot \tau(D^\star{ \boldsymbol u})
	\right\|_{p'}	\leq F_{4,3} \,\delta,$$
where $F_{4,3}= 4\sqrt{3}m {\cal F}_1\big(\|D^\star{ \boldsymbol u}\|_p\big)\|D^\star{ \boldsymbol u}\|_p $.
Similarly, by using the Sobolev inequality (\ref{sobolev0}) and the Korn inequality  (\ref{korn-l2})
$$\begin{array}{ll}\left\| u\otimes u-w\otimes w\right\|_{2}+\delta \left\|u_3^2\right\|_{2}
	&\leq	\left\|\left(u+w\right)
	\otimes \left( u-w\right)\right\|_{2}+\delta \left\|u_3\right\|^2_{4}\vspace{2mm}\\
	&\leq	\left\|u+w\right\|_{4}
	\left\|u-w\right\|_{4}+\delta \left\|u_3\right\|^2_{4}\vspace{2mm}\\
	&\leq S_{2,4} 
	\left\|\nabla\left(u-w\right)
	\right\|_2\left\|u+w\right\|_{4}+\delta\, \left\|u_3\right\|^2_{4}\vspace{2mm}\\
&\leq S_{2,4} \,
	\sqrt{2}
	\left\|D\left(\boldsymbol u-\boldsymbol w\right)
	\right\|_{2}\left\|u+w\right\|_{4}+\delta\, \left\|u_3\right\|^2_{4}\end{array}$$
and thus 
$$\begin{array}{ll}\left\| u\otimes u-w\otimes w\right\|_{p'}+\delta \left\|u_3^2\right\|_{p'}
	&\leq |\Sigma|^{\frac{1}{p'}-\frac{1}{2}} 
	\left(\left\| u\otimes u-w\otimes w\right\|_{2}+\delta \left\|u_3^2\right\|_{2}\right)\vspace{1mm}\\
	&\leq F_{4,4} \left(
	\left\|D\left(\boldsymbol u-\boldsymbol w\right)\right\|_p+\delta\right),\end{array}$$
where $F_{4,4}=|\Sigma|^{\frac{1}{p'}-\frac{1}{2}}\left(\sqrt{2}\, S_{2,4}\,
	|\Sigma|^{\frac{1}{2}-\frac{1}{p}}\left\|u+w\right\|_{4}+\left\|u_3\right\|^2_{4}\right)$. On the other hand, due to (\ref{est_psi_2}) we have
	$$\begin{array}{ll}
	\left\|\sigma(w_3)\right\|_{p'}\leq E_{\alpha,p}
	 \left\|\nabla w_3\right\|_{2}^\alpha=F_{4,5}.\end{array}$$
Combining these estimates, we deduce that
	\begin{align}\label{difference_est8}\left|I_4\right|&\leq \tilde\kappa\left(\left(F_{4,2}+{\cal R}e\, F_{4,4}\right)
\left\|D{ \boldsymbol u}-D{ \boldsymbol w}\right\|_{p}+ \left(F_{4,2}+F_{4,3}+{\cal R}e\,F_{4,4}+F_{4,5}\right)\delta\right) \left\|\nabla\cdot \boldsymbol\varphi\right\|_p\nonumber\\
	&\leq F_ 4\left(\left\|D\left({ \boldsymbol u}-{ \boldsymbol w}\right)\right\|_{p}+ 
	\delta\right)\left\|\nabla\cdot \boldsymbol\varphi\right\|_p.\end{align}
$\circ$ Finally, taking into account (\ref{est_psi_1}) we obtain
	\begin{align}\label{difference_est_psi}\left|I_5\right|&=\delta \left|\left(\sigma(w_3),
	\varphi_2\right)\right|\leq D_{2,\alpha} \,\delta
	\left\|\nabla w_3\right\|^\alpha_2
	\left\|\nabla\varphi_2\right\|_{2}\nonumber\\
	&\leq D_{2,\alpha}\,\delta
	\left\|\nabla w_3\right\|^\alpha_2
	\left\|\nabla\boldsymbol\varphi
	\right\|_{2}\leq \sqrt{2} D_{2,\alpha}\,\delta
	\left\|\nabla w_3\right\|^\alpha_2
	\left\|D\boldsymbol\varphi
	\right\|_{2}\nonumber\\
	&\leq F_5\,\delta
	\left\|D\boldsymbol\varphi
	\right\|_{p},\end{align}
where $F_5=\sqrt{2}\,D_{2,\alpha} |\Sigma|^{\frac{1}{2}-\frac{1}{p}}
	\left\|\nabla w_3\right\|^\alpha_2$.\vspace{2mm}\\
$\circ$ Combining (\ref{difference1}), (\ref{difference_est3})-(\ref{difference_est_psi}) yields
	$$\begin{array}{ll}\left|\left(\boldsymbol\tau(D{ \boldsymbol u})-{\boldsymbol\tau}\left(D  { \boldsymbol w}\right),D \boldsymbol\varphi\right)\right|&\leq	\left(F_1+F_2+{\cal R}e\, F_3+F_5\right)\, \delta\left\|D\boldsymbol\varphi\right\|_{p}+
	{\cal R}e\,F_3 \left\|D \left(\boldsymbol u-\boldsymbol w\right)\right\|_{2}\left\|D\boldsymbol\varphi\right\|_2\vspace{2mm}\\
	&+F_ 4\left(\left\|D\left({ \boldsymbol u}-{ \boldsymbol w}\right)\right\|_{p}+ 
	\delta\right)\left\|\nabla\cdot \boldsymbol\varphi\right\|_p.\end{array}$$
Setting $\boldsymbol\varphi=\boldsymbol u-\boldsymbol w$, taking into account $(\ref{estim_S3})$ and the estimates associated to 
$\boldsymbol u$ and $\boldsymbol w$, and using the Young inequality,  it follows that for every $\varepsilon>0$,
	$$\begin{array}{ll}\tfrac{1}{2^{p}(p-1)}
	\left\|D\left({ \boldsymbol u}-{ \boldsymbol w}\right)
	\right\|_{p}^p+
	\tfrac{1}{2}\left\|D\left({ \boldsymbol u}-{ \boldsymbol w}
	\right)\right\|_{2}^2
	&\leq \left({\boldsymbol\tau}\left(D{ \boldsymbol u}\right)
	-{\boldsymbol\tau}\left(D{ \boldsymbol w}\right),
	D\left({ \boldsymbol u}-{ \boldsymbol w}\right)\right)
	\vspace{2mm}\\
	&\leq	\left(F_1+F_2+{\cal R}e\, F_3+F_5\right)\, \delta\left\|D\left(\boldsymbol u-\boldsymbol w\right)\right\|_{p}\vspace{2mm}\\
	& \ +{\cal R}e\,F_3 \left\|D \left(\boldsymbol u-\boldsymbol w\right)\right\|_{2}^2\vspace{2mm}\\
	& \ +F_ 4 \left(\left\|D\left({ \boldsymbol u}-{ \boldsymbol w}\right)\right\|_{p}+ 
	\delta\right)\left\|\tfrac{\delta}{B}u_2\right\|_p\vspace{2mm}\\
	&\leq	C_1\left(\delta\left\|D\left({ \boldsymbol u}-{ \boldsymbol w}\right)\right\|_{p}+
	{\cal R}e \left\|D \left(\boldsymbol u-\boldsymbol w\right)\right\|_{2}^2+\delta^2\right)\vspace{2mm}\\
	& \leq 
	\varepsilon \left\|D\left({ \boldsymbol u}-{ \boldsymbol w}\right)\right\|_{p}^p+C_2(\varepsilon)\, \delta^{p'}\vspace{2mm}\\
	& \ +C_1\left({\cal R}e
\left\|D\left(\boldsymbol u-\boldsymbol w\right)\right\|_{2}^2
	+\delta^2\right),\end{array}$$
where $C_1$ is a positive constant only depending on $\Sigma$, $p$ and $m$. Observing that $\delta^{2}<\delta^{p'}$, choosing $\varepsilon= \tfrac{1}{2^{p+1}(p-1)}$ and assuming that $C_1{\cal R}e<\tfrac{1}{4}$, we deduce that 
	$$ \tfrac{1}{2^{p+1}(p-1)}\left\|D \left(\boldsymbol u-\boldsymbol w\right)\right\|_{p}^p+
	\tfrac{1}{4}\left\|D \left(\boldsymbol u-\boldsymbol w\right)\right\|_{2}^2\leq 
	C_3\,\delta^{p'}$$
and the claimed result is proved. $\hfill\Box$
\section{Shear-thinning flows}
\setcounter{equation}{0}
Let us now consider the case of shear-thinning fluids (corresponding to $p<2$). As for the shear-thickening fluids, we derive a  Korn inequality, establish some estimates on the convective term and on the extra stress tensor and prove existence and uniqueness results. As previously observed in Section \ref{section_shear_thic_introd}, we will restrict the exponent $p$ in order to ensure the uniqueness of the solution and carry out the approximation analysis with respect to $\delta$.
\subsection{On the Korn inequality}
Let us notice that if the classical Korn inequality can be applied to
the tensor $D^\star\boldsymbol u$ with $\delta=0$ or to $Du=\left(D_{ij}^\star\boldsymbol u\right)_{i=1,2}$ with $\delta>0$, this is no more necessarily the case if we consider $D^\star\boldsymbol u$ with $\delta>0$. The difficulty, basically related with the term $D_{23}^\star \boldsymbol u=
\tfrac{\partial u_3}{\partial x_2}-\tfrac{\delta}{B} u_3$, is overcome in the case $p=2$ by using the Hilbert setting and the fact that ${u_3}_{\mid \partial \Sigma}=0$. Indeed, since 
	$$\left(\tfrac{\partial u_3}{\partial x_2}, u_3\right)=0,$$ we obtain 
	$$\left\|D_{23}^\star \boldsymbol u\right\|_{2,B}^2=
	\left\|\tfrac{\partial u_3}{\partial x_2}-\tfrac{\delta}{B} u_3\right\|_{2,B}^2=\left\|\tfrac{\partial u_3}{\partial x_2}\right\|_{2,B}^2+\left\|\tfrac{\delta}{B} u_3\right\|_{2,B}^2.$$
This argument is one of the key points in the proof of the corresponding Korn inequality (see Section \ref{korn_thick}) but 
does not apply in the $L^p$ setting. The issue is overcome by using the Poincar\'e inequality (\ref{poincare2}) that involves only the first component $\tfrac{\partial u_3}{\partial x_1}$ of the gradient $\nabla u_3$. 
\begin{lemma} Let $ { \boldsymbol u}=(u_1,u_2,u_3) \in  { \boldsymbol W}^{1,p}_0(\Sigma)$ with $1<p<\infty$. Then
	\begin{equation}\label{korn-inequality} C_K\left\|\nabla^\star { \boldsymbol u}\right\|_{p,B}\leq 
	\left\|D^\star{ \boldsymbol u}\right\|_{p,B}\end{equation}
with  $C_K=\tfrac{C_{K,1}\left(n  m\right)^{-\frac{1}{p}}}{2(1+\delta m)}$, where $C_{K,1}$ is the classical Korn constant 
in $\boldsymbol W^{1,p}_0(\Sigma)$.
\end{lemma}
{\bf Proof.} Let us first observe that due to (\ref{poincare2}), we have
	\begin{equation}\label{korn_in1}\left\|\tfrac{\delta}{B}u_3\right\|_{p,B}\leq 
	\left\|\tfrac{\delta}{B}
	\tfrac{\partial u_3}{\partial x_1}\right\|_{p,B}\leq 2\delta m
	\left\|D_{13}^\star \boldsymbol u\right\|_{p,B}.\end{equation}	
It follows that
	$$\begin{array}{ll}
	\left\|D_{23}\boldsymbol u\right\|_{p,B}^p&\leq 
	\left(\left\|D_{23}^\star\boldsymbol u\right\|_{p,B}
	+\tfrac{1}{2}\left\|\tfrac{\delta}{B}u_3\right\|_{p,B}\right)^p
	\vspace{1mm}\\
	&\leq 
	\left(\left\|D_{23}^\star\boldsymbol u\right\|_{p,B}
	+\delta m
	\left\|D_{13}^\star \boldsymbol u\right\|_{p,B}\right)^p
	\vspace{1mm}\\
	&\leq 2^{p-1}\left(\left\|D_{23}^\star\boldsymbol u
	\right\|_{p,B}^p
	+\left(\delta m\right)^p
	\left\|D_{13}^\star \boldsymbol u\right\|_{p,B}^p\right)
	\end{array}$$
and thus
	\begin{align}\label{korn_in2}\left\|D\boldsymbol u\right\|_{p,B}^p&=
	\left\|D u\right\|_{p,B}^p+
	2\left\|D_{13} \boldsymbol u\right\|_{p,B}^p+
	2\left\|D_{23} \boldsymbol u\right\|_{p,B}^p\nonumber\\
	&\leq \left\|D u\right\|_{p,B}^p+
	2\left\|D_{13}^\star \boldsymbol u\right\|_{p,B}^p
	+2^{p}\left\|D_{23}^\star\boldsymbol u
	\right\|_{p,B}^p
	+\left(2\delta m\right)^p
	\left\|D_{13}^\star \boldsymbol u\right\|_{p,B}^p.\end{align}
On the other hand, taking into account the definition of $\nabla^\star$ and using the classical Korn inequality,  we can easily see that
	\begin{align}\label{korn_in3}
	\left\|\nabla^\star{ \boldsymbol u}\right\|_{p,B}^p&=
	\left\|\nabla{ \boldsymbol u}\right\|_{p,B}^p+
	\left\|\tfrac{\delta}{B}u_2\right\|_{p,B}^p+
	\left\|\tfrac{\delta}{B}u_3\right\|_{p,B}^p\nonumber\\
	&\leq n\left\|\nabla{ \boldsymbol u}\right\|_{p}^p+
	\left\|\tfrac{\delta}{B}u_2\right\|_{p,B}^p+
	\left\|\tfrac{\delta}{B}u_3\right\|_{p,B}^p\nonumber\\
	&\leq \tfrac{n}{C_{K,1}^p}\left\|D{ \boldsymbol u}\right\|_{p}^p+
	\left\|\tfrac{\delta}{B}u_2\right\|_{p,B}^p+
	\left\|\tfrac{\delta}{B}u_3\right\|_{p,B}^p\nonumber\\
	&\leq \tfrac{mn}{C_{K,1}^p}
	\left\|D{ \boldsymbol u}\right\|_{p,B}^p+
	\left\|\tfrac{\delta}{B}u_2\right\|_{p,B}^p+
	\left\|\tfrac{\delta}{B}u_3\right\|_{p,B}^p.\end{align}
Combining (\ref{korn_in1})-(\ref{korn_in3}), we deduce that
	$$\begin{array}{ll}&
	\left\|\nabla^\star{ \boldsymbol u}\right\|_{p,B}^p\vspace{2mm}\\
	&\leq \tfrac{mn}{C_{K,1}^p}
	\left(\left\|D u\right\|_{p,B}^p+
	2\left\|D_{13}^\star \boldsymbol u\right\|_{p,B}^p
	+2^{p}\left\|D_{23}^\star\boldsymbol u
	\right\|_{p,B}^p
	+\left(2\delta m\right)^p
	\left\|D_{13}^\star \boldsymbol u\right\|_{p,B}^p\right)
	\vspace{2mm}\\
	&+\left(2\delta m\right)^p
	\left\|D_{13}^\star \boldsymbol u\right\|_{p,B}^p+
	\left\|D_{33}^\star \boldsymbol u\right\|_{p,B}^p\vspace{2mm}\\
	&\leq \tfrac{mn}{C_{K,1}^p}\left(
	\left\|Du\right\|_{p,B}^p+
	\left(2+\left(2\delta m\right)^p\left(1+\tfrac{C_{K,1}^p}{mn}\right)\right)\left\|D_{13}^\star \boldsymbol u\right\|_{p,B}^p+
	2^p\left\|D_{23}^\star \boldsymbol u\right\|_{p,B}^p+
	\left\|D_{33}^\star \boldsymbol u\right\|_{p,B}^p\right)\vspace{2mm}\\
	&\leq \tfrac{\left(2+2\delta m\right)^p mn}{C_{K,1}^p}\left(
	\left\|D u\right\|_{p,B}^p+
	2\left\|D_{13}^\star \boldsymbol u\right\|_{p,B}^p+
	2\left\|D_{23}^\star \boldsymbol u\right\|_{p,B}^p+
	\left\|D_{33}^\star \boldsymbol u\right\|_{p,B}^p\right)\vspace{2mm}\\
	&=\tfrac{1}{C_{K}^p} \left\|D^\star \boldsymbol u\right\|_{p,B}^p
	\end{array}$$
and the claimed result is proved. $\hfill\Box$
\subsection{Estimates on the convective term and extra stress tensor}
We begin by a continuity property in $\boldsymbol W^{1,p}_0(\Sigma)$ of the trilinear form $a_\star$.
\begin{lemma} \label{convective_2_thin} Let ${ \boldsymbol u}$, ${ \boldsymbol v}$ and ${ \boldsymbol w}$ be in ${ \boldsymbol W}^{1,p}_0(\Sigma)$ with $\tfrac{3}{2}\leq p<2$. Then the following estimate holds
	$$\left|a_\star({ \boldsymbol u},{ \boldsymbol v},B{ \boldsymbol w})\right|\leq 
	\kappa_6
	\left\|D^\star{ \boldsymbol u}\right\|_{p,B}\left\|D^\star{ \boldsymbol v}\right\|_{p,B}
	\left\|D^\star{ \boldsymbol w}\right\|_{p,B},$$
with $\kappa_6=\tfrac{ n m^{\frac{3}{p}}}{C_K^3}\left( S_{p,2p'}\right)^2$  and where $C_K$ is the Korn constant given in $(\ref{korn-inequality})$.
\end{lemma}
{\bf Proof.}  H\"older's inequality and Sobolev's inequality with $r=2p'$ and $q=p$ show that
	$$\begin{array}{ll}\label{a_alpha}\left|a_\star({ \boldsymbol u},{ \boldsymbol v},B{ \boldsymbol w})\right|
	&\leq n\left\|{ \boldsymbol u}\right\|_{2p'}
	\left\|\nabla^\star { \boldsymbol v}\right\|_{p}
	\left\|{ \boldsymbol w}\right\|_{2p'}\vspace{2mm}\\
	&\leq n\left( S_{p,2p'}\right)^2
	\left\|\nabla  \boldsymbol u\right\|_{p}
	\left\|\nabla^\star { \boldsymbol v}\right\|_{p}
	\left\|\nabla \boldsymbol w\right\|_{p}\vspace{2mm}\\
	&\leq n m^{\frac{3}{p}}\left( S_{p,2p'}\right)^2
	\left\|\nabla  \boldsymbol u\right\|_{p,B}
	\left\|\nabla^\star { \boldsymbol v}\right\|_{p,B}
	\left\|\nabla \boldsymbol w\right\|_{p,B}.\end{array}$$
The estimate follows by using the Korn inequality (\ref{korn-inequality}). $\hfill\Box$\vspace{2mm}\\
In the remaining part of this section, we study some properties of the extra stress tensor and derive some associated estimates. We begin by an auxiliary result that will be useful in the sequel.
\begin{lemma}\label{desigualdade3} Let $1<p< 2$ and let $H_1\in L^{\frac{p}{2-p}}(\Sigma)$, $H_2\in L^1(\Sigma)$ and $H_3\in L^p(\Sigma)$ be non negative functions satisfying
	$$H_3(x)^2\leq H_1(x)H_2(x) 
	\qquad \mbox{for a.e.} \ x\in \Sigma.$$ 
Then,
   	 $$ \left\|H_3\right\|_p^2\leq \left\|H_1\right\|_{\frac{p}{2-p}}\left\|H_2\right\|_1.$$
\end{lemma}
{\bf Proof.} Taking into account the condition satisfied by $H_1$, $H_2$, $H_3$,  integrating and using the H\"{o}lder inequality, we obtain
	$$\begin{array}{ll}
	\left\|H_3\right\|_p^p&=\displaystyle
	\int_\Sigma \left(H_3(x)^2\right)^\frac{p}{2}dx
	\leq \int_\Sigma H_1(x)^\frac{p}{2}H_2(x)^\frac{p}{2}dx
	\vspace{2mm}\\
	&\leq \left\|H_1^\frac{p}{2}\right\|_{\frac{2}{2-p}}
  	  \left\|H_2^\frac{p}{2}\right\|_{\frac{2}{p}}
	=\left\|H_1\right\|_{\frac{p}{2-p}}^\frac{p}{2}\left\|H_2\right\|_{1}^\frac{p}{2}\end{array}$$
and the proof is complete.$\hfill\Box$\vspace{2mm}\\
The next result deals with continuity, coercivity and monotonicity results for the extra stress tensor $\boldsymbol\tau $.
\begin{proposition}\label{ineq1} Assume that $1<p<2$ and let  $\boldsymbol f, \boldsymbol g\in { \boldsymbol L}^{p}(\Sigma,\mathbb{R}^{3\times 3})$. Then the following estimates hold
	\begin{description}
 \item {Continuity.}
\begin{equation}\label{estim_S5}\left\|\left(1+|\boldsymbol f|^2\right)^{\frac{p-2}{2}}\boldsymbol g\right\|_{p',B}
	\leq \left\|\boldsymbol g\right\|_{p,B}^{p-1} \qquad \mbox{if} \ |\boldsymbol g|\leq |\boldsymbol f|,\end{equation}
\begin{equation}\label{estim_S7_2}\left\|{\boldsymbol\tau }\left(\boldsymbol f\right)-{\boldsymbol\tau }\left(\boldsymbol g\right)\right\|_{p',B}\leq  2C_{p}\left\|\boldsymbol f-\boldsymbol g\right\|_{p,B}^{p-1} \qquad \mbox{with} \ C_{p}=1+2^{\frac{2-p}{2}},\end{equation}
 \item {Coercivity.}
	\begin{equation}\label{estim_S6}\left({\boldsymbol\tau }\left(\boldsymbol f\right),B \boldsymbol f\right)\geq 
	\tfrac{2\|\boldsymbol f\|^2_{p,B}}{\left(\|B\|_1
	+\|\boldsymbol f\|_{p,B}^p\right)^{\frac{2-p}{p}}}, \end{equation}
 \item {Monotonicity.}
\begin{equation}\label{estim_S7}\left({\boldsymbol\tau }\left(\boldsymbol f\right)-
	{\boldsymbol\tau }\left(\boldsymbol g\right),B (\boldsymbol f-\boldsymbol g)\right)\geq 
	\tfrac{2(p-1)\|\boldsymbol f-\boldsymbol g\|^2_{p,B}}{\left(\|B\|_1
	+\|\boldsymbol f\|_{p,B}^p+\|\boldsymbol g\|_{p,B}^p
	\right)^{\frac{2-p}{p}}}.\end{equation}
\end{description}
\end{proposition}
{\bf Proof.}  Standard calculation show that if $|\boldsymbol g|\leq |\boldsymbol f|$, then
	$$\begin{array}{ll}
	\left\|\left(1+|\boldsymbol f|^2\right)^{\frac{p-2}{2}}
	\boldsymbol g
	\right\|_{p',B}^{p'}&=\displaystyle\int_\Sigma B
	\left((1+ |\boldsymbol f|^2)^{\frac{p-2}{2}}|\boldsymbol g|
	\right)^{p'}\,dx
	\vspace{2mm}\\
	&\leq 
	\displaystyle\int_\Sigma B\left((1+ |\boldsymbol g|^2)^{\frac{p-2}{2}}|\boldsymbol g|\right)^{p'}\,dx\vspace{2mm}\\
	&\leq
	\displaystyle\int_\Sigma B\left(( 
	|\boldsymbol g|^2)^{\frac{p-2}{2}}
	|\boldsymbol g|\right)^{p'}\,dx
	=\|\boldsymbol g\|_{p,B}^{p}\end{array}$$
which gives (\ref{estim_S5}). On the other  hand, due to  the monotonicity property in Lemma \ref{tensor_proper1}, we have
	\begin{align}\label{monot_B}\tfrac{1}{2(p -1)}\left({\boldsymbol\tau }\left(f\right)-
	{\boldsymbol\tau }\left(\boldsymbol g\right)\right):B(\boldsymbol f-\boldsymbol g)&\geq B\left(1+|\boldsymbol f|^{2}+|\boldsymbol g|^2\right)^{\frac{2-p}{2}}|\boldsymbol f-\boldsymbol g|^2\nonumber\\
	&= \left(B^{\frac{2}{p}}+|B^{\frac{1}{p}}\boldsymbol f|^{2}+|B^{\frac{1}{p}}\boldsymbol g|^2\right)^{\frac{2-p}{2}}\big|B^{\frac{1}{p}}(\boldsymbol f-\boldsymbol g)\big|^2.
	\end{align}
 Set
	$$\begin{array}{ll}H_1=\left(B^{\frac{2}{p}}+|B^{\frac{1}{p}}\boldsymbol f|^{2}+|B^{\frac{1}{p}}\boldsymbol g|^2\right)^{\frac{2-p}{2}},\vspace{2mm}\\
	H_2=\tfrac{1}{2(p -1)}\left({\boldsymbol\tau }\left(\boldsymbol f\right)-{\boldsymbol\tau }\left(\boldsymbol g\right):B(\boldsymbol f-\boldsymbol g)\right),\qquad 
	H_3=\big|B^{\frac{1}{p}}(\boldsymbol f-\boldsymbol g)\big|.\end{array}$$
Since $\boldsymbol f$ and $\boldsymbol g$ belong to $\boldsymbol L^{p}(\Sigma,\mathbb{R}^{3\times 3})$, it is easy to see that  $H_1\in L^{\frac{p}{2-p}}(\Sigma)$, $H_2\in
 L^1(\Sigma)$ and $H_3\in  L^p(\Sigma)$ and that due to  (\ref{monot_B}), we have
	$$H_3(x)^2\leq H_1(x)H_2(x) \qquad \mbox{for a.e.} \ x\in \Sigma.$$ 
Due to Lemma \ref{desigualdade3}, we obtain
	$$\begin{array}{ll}2(p -1) \left\|\boldsymbol f-\boldsymbol g
	\right\|_{p,B}^2
	&\leq\left\|\left(B^{\frac{2}{p}}+|B^{\frac{1}{p}}\boldsymbol f
	|^{2}+|B^{\frac{1}{p}}\boldsymbol g|^2\right)^{\frac{2-p}{2}}
	\right\|_{\frac{p}{2-p}}
	\left\|\left({\boldsymbol\tau }\left(\boldsymbol f\right)-{\boldsymbol\tau }\left(\boldsymbol g\right)\right):B(\boldsymbol f-\boldsymbol g)\right\|_1\\
	&\leq\left(\|B\|_1
	+\|\boldsymbol f\|_{p,B}^p+\|\boldsymbol g\|_{p,B}^p
	\right)^{\frac{2-p}{p}}
	\left({\boldsymbol\tau }\left(\boldsymbol f\right)-{\boldsymbol\tau }\left(\boldsymbol g\right),B(\boldsymbol f-\boldsymbol g)\right)\end{array}$$
which gives (\ref{estim_S7}). Estimate (\ref{estim_S6}) can be obtained very similarly by using the coercivity condition in Lemma \ref{tensor_proper1}. Finally, by taking into account Lemma \ref{tensor_proper1}, we have
$$\left\|{\boldsymbol\tau }\left(\boldsymbol f\right)-{\boldsymbol\tau }\left(\boldsymbol g\right)\right\|_{p',B}
	\leq \displaystyle 2C_{p}\left\||\boldsymbol f-\boldsymbol g|^{p-1}\right\|_{p',B} =2C_{p}\left\|\boldsymbol f-\boldsymbol g\right\|^{p-1}_{p,B}$$
which gives estimate (\ref{estim_S7_2}). $\hfill\Box$
\begin{proposition} Assume that $\tfrac{3}{2}\leq p< 2$ and let  $\boldsymbol f\in { \boldsymbol L}^p(\Sigma,\mathbb{R}^{3\times 3})$. Then the following estimate holds
	\begin{equation}\label{diff_div_thin}\left\|\nabla^\star\cdot \boldsymbol\tau (\boldsymbol f)-\nabla\cdot \boldsymbol\tau (\boldsymbol f)\right\|_{p'}\leq 4\delta m
	\left(\left\|f\right\|_{p}^{p-1}
	+\left\|f_{33}\right\|_{p}^{p-1}+
	\left\|f_{23}\right\|_{p}^{p-1}\right).\end{equation}
where $f=(f_{ij})_{ij=1,2}$.
\end{proposition}
{\bf Proof.} Taking into account (\ref{estim_S5}), the estimate can be derived by following step by step the proof of Proposition \ref{prop_div_diff}.$\hfill \Box$
\subsection{Existence and uniqueness of shear-thinning flows}
\begin{proposition} \label{existence_state3} Assume that $\tfrac{3}{2}\leq p<2$ and let $\boldsymbol u$ be a weak solution of $(\ref{weak_formulation})$. Then, estimates
 $(\ref{main_estimates_1_thin})$-$(\ref{est_u_thin})$ hold.
\end{proposition}
{\bf Proof.} The proof is split into three steps.\vspace{1mm}\\
\underline{\it Step 1.} Let us set ${\boldsymbol\varphi}={ \boldsymbol u}$ in the weak formulation (\ref{weak_formulation}) and use Lemma \ref{prop_a1} and (\ref{estim_S6}) 
to obtain
	\begin{equation}\label{coerciv_B}\tfrac{2\left\|D^\star{ \boldsymbol u}\right\|_{p,B}^2}{\left(\|B\|_1
	+\left\|D^\star{ \boldsymbol u}\right\|_{p,B}^{p}\right)^{\frac{2-p}{p}}}
	\leq  \left({\boldsymbol\tau }\left(D^\star{ \boldsymbol u}\right),BD^\star{ \boldsymbol u}\right)=\left(G,u_3\right).\end{equation}
On the other hand, classical arguments together with Poincar\'e inequality (\ref{poincare}) give
	$$\left(G,u_3\right)\leq 2\kappa_1
	\left\|u_3\right\|_{p,B}\leq 2\kappa_1
	\left\|\tfrac{\partial u_3}{\partial x_1}\right\|_{p,B},$$
where $\kappa_1=\tfrac{m^{\frac{1}{p}}}{2}|G||\Sigma|^{\frac{1}{p'}}$. Observing that
	$$\left\|D^\star{ \boldsymbol u}\right\|^p_{p,B}=
	\int_\Sigma B\left(|D^\star{ \boldsymbol u}|^2\right)^{\frac{p}{2}}\,dx\geq 
	\left\|\tfrac{\partial u_3}{\partial x_1}\right\|^p_{p,B},$$
we deduce that
	\begin{equation}\label{sec_m_B}\left(G,u_3\right)\leq 2\kappa_1
	\left\|D^\star{ \boldsymbol u}\right\|_{p,B}.\end{equation}
Due to (\ref{coerciv_B}) and (\ref{sec_m_B}), we have
	$$\left\|D^\star{ \boldsymbol u}\right\|_{p,B}\leq 
	\kappa_1\left(\|B\|_1
	+\left\|D^\star{ \boldsymbol u}\right\|^{p}_{p,B}\right)^{\frac{2-p}{p}}$$
and thus
	\begin{equation}\label{est_s5gen}\left\|D^\star{ \boldsymbol u}\right\|_{p,B}^{\frac{p}{2-p}}\leq 
	\kappa_1^{\frac{p}{2-p}}\left(\|B\|_1
	+\left\|D^\star{ \boldsymbol u}\right\|^{p}_{p,B}\right).
	\end{equation}
The Young inequality yields
	\begin{equation}\label{est_s6gen}
	\kappa_1^{\frac{p}{2-p}}\left\|D^\star{ \boldsymbol u}\right\|_{p,B}^{p}\leq (2-p)
	\left\|D^\star{ \boldsymbol u}\right\|^{\frac{p}{2-p}}_{p,B}+
	(p-1)\kappa_1^{\frac{p}{(2-p)(p-1)}}\end{equation}
and by combining (\ref{est_s5gen}) and (\ref{est_s6gen}), we deduce that
	$$(p-1)\left\|D^\star{ \boldsymbol u}\right\|^{\frac{p}{2-p}}_{p,B}\leq 
	\kappa_1^{\frac{p}{2-p}}\|B\|_1+(p-1)\kappa_1^{\frac{p}{(2-p)(p-1)}}.$$
Consequently
	$$\left\|D^\star{ \boldsymbol u}\right\|_{p,B}\leq \kappa_1 \left(\tfrac{\|B\|_1}{p-1}+\kappa_1
	^{p'}\right)^{\frac{2-p}{p}} $$
and estimate  $(\ref{main_estimates_1_thin})$ is proved.\vspace{1mm}\\
\underline{\it Step 2.} Let us now prove (\ref{main_estimates_2_thin}). Similar arguments together with 
the coercivity property and (\ref{sec_m_B}) show that
	\begin{align}\label{estima_M1}\left\|D^\star{ \boldsymbol u}\right\|_{p,B}^p&=\displaystyle\int_ {\Sigma_{\boldsymbol u}}
	B|D^\star{ \boldsymbol u}|^p\,dx+
	\int_{\Sigma\setminus \Sigma_{\boldsymbol u}}B|D^\star{ \boldsymbol u}|^p\,dx\vspace{2mm}\nonumber\\
	&\leq \displaystyle \int_ {\Sigma_{\boldsymbol u}}
	\tfrac{|B^\frac{1}{p}D^\star{ \boldsymbol u}|^2}{|B^\frac{1}{p} D^\star{ \boldsymbol u}|^{2-p}}\,dx+
	\|B\|_1\nonumber\vspace{2mm}\\
	& \leq 2^{\frac{2-p}{2}} \displaystyle  \int_ {\Sigma_{\boldsymbol u}}\tfrac{|B^\frac{1}{p}D^\star{ \boldsymbol u}|^2}{\left(B^\frac{2}{p}+|B^\frac{1}{p}D^\star{ \boldsymbol u}|^2\right)^{\frac{2-p}{2}}}\,dx+\|B\|_1\nonumber\vspace{2mm}\\
	&\leq 2^{\frac{2-p}{2}}\int_ {\Sigma_{\boldsymbol u}}
{\boldsymbol\tau }(D^\star{ \boldsymbol u}):BD^\star{ \boldsymbol u}\,dx+ \|B\|_1\nonumber\vspace{2mm}\\
	&\leq 2^{\frac{2-p}{2}} \left({\boldsymbol\tau }(D^\star{ \boldsymbol u}),BD^\star{ \boldsymbol u}\right)+ \|B\|_1= 2^{\frac{2-p}{2}}\left(G,u_3\right)+
	 	\|B\|_1\nonumber\\
	& \leq 2^{\frac{2-p}{2}}\kappa_1\left\|D^\star{ \boldsymbol u}\right\|_{p,B}+ \|B\|_1,\end{align}
where  $\Sigma_{\boldsymbol u}=\left\{x\in \Sigma\mid B^\frac{1}{p}|D^\star{ \boldsymbol u}(x)|\geq 1\right\}$.
The Young  inequality yields
	\begin{equation}\label{estima_M4}2^{\frac{2-p}{2}}\kappa_1
	 \left\|D^\star{ \boldsymbol u}\right\|_{p,B}\leq 
	\tfrac{1}{p'}\left(2^{\frac{2-p}{2}}\kappa_1\right)^{p'}+\tfrac{1}{p}
	\left\|D^\star{ \boldsymbol u}\right\|_{p,B}^p\end{equation}
and the claimed result follows by combining (\ref{estima_M1}) and (\ref{estima_M4}).\vspace{2mm}\\
\underline{\it Step 3.} By taking into account (\ref{korn_in1}), we have
	$$\begin{array}{ll}\left\|\nabla u_3\right\|_{p,B}&\leq \left\|\tfrac{\partial u_3}{\partial x_1}\right\|_{p,B}
	+\left\|\tfrac{\partial u_3}{\partial x_2}-\tfrac{\delta}{B}u_3\right\|_{p,B}
	+\left\|\tfrac{\delta}{B}u_3\right\|_{p,B}\vspace{2mm}\\
	&\leq 2\left(1+\delta m\right)\|D_{13}^\star \boldsymbol u\|_{p,B}+2\|D_{23}^\star \boldsymbol u\|_{p,B}\vspace{2mm}\\
	&\leq 
	2^{2-p}\left(1+\delta m \right)\|D^\star \boldsymbol u\|_{p,B}\end{array}$$
and estimate (\ref{est_u3_thin}) is then a consequence of (\ref{main_estimates_1_thin}). Let us finally prove (\ref{est_u_thin}). Setting $\boldsymbol\varphi=u$ in (\ref{weak_formulation}) we obtain
	\begin{equation}\label{est_u_thin1}\int_\Sigma \left(B^{\frac{2}{p}}+ 
	|B^{\frac{1}{p}}D^\star \boldsymbol u|^2\right)^{\frac{p-2}{2}}|B^{\frac{1}{p}}D u|^2dx
	\leq \tfrac{1}{2}{\cal R}e\delta	\left(u_3^2,u_2\right).\end{equation}
Unlike the proof of estimates (\ref{main_estimates_1_thin}) and 
(\ref{main_estimates_2_thin}), the coercivity property is not immediat and cannot be used. Let us then set 
	$$\begin{array}{ll}H_1=\left(B^{\frac{2}{p}}+ |B^{\frac{1}{p}}D^\star\boldsymbol u|^2\right)^{\frac{2-p}{2}}, \vspace{2mm}\\
H_2=\left(B^{\frac{2}{p}}+ |B^{\frac{1}{p}}D^\star\boldsymbol u|^2\right)^{\frac{p-2}{2}}|B^{\frac{1}{p}} Du|^2, 
	\qquad H_3=|B^{\frac{1}{p}}D u|.\end{array}$$
It is easy to see that $H_1\in \boldsymbol L^{\frac{p}{2-p}}(\Sigma)$ and that $H_3\in \boldsymbol L^p(\Sigma)$. Moreover, since 
$|D^\star \boldsymbol u|\geq |Du|$, we have
	$$H_2\leq \left(B^{\frac{2}{p}}+ |B^{\frac{1}{p}}D u|^2\right)^{\frac{p-2}{2}}|B^{\frac{1}{p}} Du|^2\leq 
	|B^{\frac{1}{p}} Du|^p\in L^1(\Sigma).$$
Using  Lemma \ref{desigualdade3}, we deduce that 
	\begin{align}\label{est_u_thin2}\left\|Du\right\|_{p,B}^2&\leq \left\|\big(B^{\frac{2}{p}}+ |B^{\frac{1}{p}}D^\star\boldsymbol u|^2\big)^{\frac{2-p}{2}}\right\|_{\frac{p}{2-p}}\left\|\big(B^{\frac{2}{p}}+ |B^{\frac{1}{p}}D^\star\boldsymbol u|^2\big)^{\frac{p-2}{2}}|B^{\frac{1}{p}} Du|^2\right\|_1\nonumber\\
	&\leq \big(\|B\|_1+ 
\|D^\star\boldsymbol u\|_{p,B}^p\big)^{\frac{2-p}{p}}  \left\|\big(B^{\frac{2}{p}}+ |B^{\frac{1}{p}}D^\star\boldsymbol u|^2\big)^{\frac{p-2}{2}}|B^{\frac{1}{p}} Du|^2\right\|_1.\end{align}
On the other hand, taking into account (\ref{est_psi_1_thin})  with $\alpha=2$ and $q=p$ and the classical Korn inequality, we have
	\begin{align}\label{est_u_thin3}
	\left|\left(u_3^2,u_2\right)\right|&\leq  
	\left(S_{p,2p'}\right)^3
	\left\|\nabla u_3\right\|_{p}^2  \left\|\nabla u_2\right\|_{p}
	\leq \tfrac{\left(S_{p,2p'}\right)^3}{C_{K,1}}
	\left\|\nabla u_3\right\|_{p}^2  \left\|D u\right\|_{p}
	\nonumber\\
	&\leq \tilde\kappa	\left\|\nabla u_3\right\|_{p,B}^2  \left\|D u\right\|_{p,B},\end{align}
where $\tilde\kappa=\tfrac{1}{C_{K,1}}\left(m^{\frac{1}{p}}S_{p,2p'}\right)^3$ .
Combining (\ref{est_u_thin1})-(\ref{est_u_thin3}) and taking into account 
(\ref{main_estimates_2_thin}), we get
	$$\begin{array}{ll}\|D u\|_{p,B}
	&\leq \tfrac{\tilde\kappa}{2}
	 \left\|\nabla u_3\right\|_{p,B}^2 \left(\|B\|_1+ 
\|D^\star\boldsymbol u\|_{p,B}^p\right)^{\frac{2-p}{p}}\delta {\cal R}e
\vspace{2mm}\\
	&\leq\tfrac{\tilde\kappa}{2}
	 \left\|\nabla u_3\right\|_{p,B}^2 \left(\tfrac{2p-1}{p-1}\|B\|_1+ 
\left(2^{\frac{2-p}{2}}\kappa_1\right)^{p'}\right)^{\frac{2-p}{p}}\delta {\cal R}e\vspace{2mm}\\
	&\leq\tilde\kappa
	 \left\|\nabla u_3\right\|_{p,B}^2 \left(\tfrac{\|B\|_1}{p-1}+ 
\kappa_1^{p'}\right)^{\frac{2-p}{p}}\delta {\cal R}e\end{array}$$
and the conclusion follows from estimate (\ref{est_u3_thin}). $\hfill\Box$\vspace{2mm}\\
{\bf Proof of Theorem \ref{main2}.} Let $\boldsymbol u^k$ be a standard Galerkin approximation. Arguments similar to those used in the proof of Proposition \ref{existence_state3} show that
	$$\left\|D^\star \boldsymbol u^k\right\|_{p,B}^p\leq p' \|B\|_1+\left(2^{\frac{2-p}{2}}\kappa_1\right)^{p'},$$
and  the sequence $\left(D^\star \boldsymbol u^k\right)_k$ is then bounded in $\boldsymbol L^p(\Sigma)$.  Taking into account 
the Korn inequality (\ref{korn-inequality}), we deduce that $\left(\nabla^\star \boldsymbol u^k\right)_k$ is bounded in $\boldsymbol L^p(\Sigma)$ and thus  $\left(\nabla \boldsymbol u^k\right)_k$ is bounded in $\boldsymbol L^p(\Sigma)$. Moreover, 
 the continuity property $(\ref{estim_S5})$ implies that
	$$\left\|\boldsymbol\tau \left(D^\star\boldsymbol u^k\right)\right\|_{p',B}\leq
	 2\left\|D^\star\boldsymbol u^k\right\|_{p,B}^{p-1}$$
and the sequence $\left(B^{\frac{1}{p'}}\boldsymbol\tau \left(D^\star\boldsymbol u^k\right)\right)_k$ is bounded in $\boldsymbol L^{p'}(\Sigma)$. There then exist a subsequence, still indexed by $k$,  $\boldsymbol u\in \boldsymbol V_B^p$ and 
$\widetilde {\boldsymbol\tau }\in \boldsymbol  L^{p'}(\Sigma)$ such that $\left(\nabla \boldsymbol u^k\right)_k$ converges to 
$\boldsymbol u$ weakly in $\boldsymbol  L^{p}(\Sigma)$  and $\left(B^{\frac{1}{p'}}\boldsymbol\tau \left(D^\star\boldsymbol u^k\right)\right)_k$ converges to $\widetilde {\boldsymbol\tau }$ weakly in $ \boldsymbol  L^{p'}(\Sigma)$. Moreover, since $p >\tfrac{4}{3}$, by using compactness results on Sobolev spaces, we deduce that $\left(\boldsymbol u^k\right)_k$ strongly converges to $\boldsymbol u$ in $L^{p'}(\Sigma )$. Taking into account these convergence results, we deduce that for every $\boldsymbol\varphi \in 
{\cal V}_B$, we have
	 $$\begin{array}{ll}\left|a_\star\left(\boldsymbol u^k,\boldsymbol u^k,B\boldsymbol\varphi\right)
	-a_\star\left(\boldsymbol u,\boldsymbol u,B\boldsymbol\varphi\right)\right|&\leq\left|a\left(\boldsymbol u^k-\boldsymbol u,\boldsymbol u^k,B\boldsymbol\varphi\right)\right|
	+\left|a_\star\left(\boldsymbol u,\boldsymbol u^k-\boldsymbol u,B\boldsymbol\varphi\right)\right|\vspace{2mm}\\
	&=\left|a_\star\left(\boldsymbol u^k-\boldsymbol u,\boldsymbol u^k,B\boldsymbol\varphi\right)\right|
	+\left|a_\star\left(B\boldsymbol u,\boldsymbol\varphi,\boldsymbol u^k-\boldsymbol u\right)\right|\vspace{2mm}\\
 	&\leq \left(\left\|\nabla^\star \boldsymbol u^k\right\|_{p,B}\left\|\boldsymbol\varphi\right\|_{\infty}+
	\left\|\boldsymbol u\right\|_{p,B} \left\|\nabla \boldsymbol\varphi\right\|_{\infty}\right)
	 \left\|\boldsymbol u^k-\boldsymbol u\right\|_{p',B}\vspace{2mm}\\
	& \longrightarrow 0 \qquad \mbox{when} \ k\rightarrow +\infty.\end{array}$$
Moreover, by passing to the limit in 
	$$\left({\boldsymbol\tau }(D^\star{ \boldsymbol u}^k),B D^\star{\boldsymbol\varphi}\right)+{\cal R}e\,a_\star\left({ \boldsymbol u}^k, { \boldsymbol u}^k,
	B{\boldsymbol\varphi}\right)=\left(G,\varphi_3\right)\qquad  \mbox{for all} \ {\boldsymbol\varphi}\in {\cal V}_B,$$
we obtain
	$$\left(\widetilde {\boldsymbol\tau },B^{\frac{1}{p}} D^\star{\boldsymbol\varphi}\right)+{\cal R}e\, 
	a_\star\left({ \boldsymbol u},{ \boldsymbol u},
	B{\boldsymbol\varphi}\right)=\left(G,\varphi_3\right) \qquad \mbox{for all} \ {\boldsymbol\varphi}\in {\cal V}_B$$
and by using the fact that ${\cal V}_B$ is dense in $\boldsymbol V_B^p$, we deduce that 
	$$\left(\widetilde {\boldsymbol\tau },B^{\frac{1}{p}} D^\star{\boldsymbol\varphi}\right)+{\cal R}e\, 
	a_\star\left({ \boldsymbol u},{ \boldsymbol u},
	B{\boldsymbol\varphi}\right)=\left(G,\varphi_3\right) \qquad \mbox{for all} \ {\boldsymbol\varphi}\in \boldsymbol V_B^p.$$
The rest of the proof for the existence of a weak solution is very similar to the first step in the proof of Theorem \ref{main1} and is omitted. To prove the uniqueness result, let us assume that ${ \boldsymbol u}$ and ${ \boldsymbol v}$ are two weak solutions of $(\ref{equation})$. Setting ${\boldsymbol\varphi}={ \boldsymbol u}-{ \boldsymbol v}$ in the corresponding weak formulation and taking into account Lemma \ref{prop_a1} 
and (\ref{estim_S7}), we obtain 
\begin{align}
	\label{uniqueness1_1}\tfrac{\left\|D^\star({ \boldsymbol u}-{ \boldsymbol v})\right\|_{p,B}^2}{
	\left(\|B\|_1+\left\|D^\star{ \boldsymbol u}\right\|_{p,B}^{p}
	+\left\|D^\star{ \boldsymbol v}\right\|_{p,B}^{p}\right)^{\frac{2-p}{p}}}
	&\leq\left({\boldsymbol\tau }(D^\star{ \boldsymbol u})-
	{\boldsymbol\tau }(D^\star{ \boldsymbol v}),BD^\star\left({ \boldsymbol u}-{ \boldsymbol v}\right)\right)\nonumber\\
	&= {\cal R}e\, a_\star\left({ \boldsymbol v},{ \boldsymbol v},B\left({ \boldsymbol u}-{ \boldsymbol v}\right)\right)-{\cal R}e \,a_\star\left({ \boldsymbol u},{ \boldsymbol u},B\left({ \boldsymbol u}-{ \boldsymbol v}\right)\right)\nonumber\\
	&=-{\cal R}e \,a\left({ \boldsymbol u}-{ \boldsymbol v},{ \boldsymbol v},B\left({ \boldsymbol u}-{ \boldsymbol v}\right)\right).
\end{align}
Lemma \ref{convective_2_thin} and estimate $(\ref{main_estimates_1_thin})$ then yield
	\begin{align} \label{visco_convective3_1}
	\left|a_\star\left({ \boldsymbol u}-{ \boldsymbol v},{ \boldsymbol v},{ \boldsymbol u}-{ \boldsymbol v}\right)\right|&\leq \kappa_6
	\left\|D^\star\left({ \boldsymbol u}-{ \boldsymbol v}\right)\right\|_{p,B}^2
	\left\|D^\star { \boldsymbol v}\right\|_{p,B}\nonumber\vspace{0mm}\\	
	&\leq
	 \kappa_1\kappa_6\left(\tfrac{\|B\|_1}{p-1}+\kappa_1
	^{p'}\right)^{\frac{2-p}{p}}
	\left\|D^\star\left({ \boldsymbol u}-{ \boldsymbol v}\right)\right\|_{p,B}^2. \end{align}
On the other hand, by taking into account estimate (\ref{main_estimates_2_thin}), we have
	\begin{align}\label{tensor_suplem_1}
	\left(\|B\|_1+\left\|D^\star{ \boldsymbol u}\right\|_{p,B}^{p}
	+\left\|D^\star{ \boldsymbol v}\right\|_{p,B}^{p}
	\right)^{\frac{2-p}{p}}
&\leq  \left(\tfrac{3p-1}{p-1}\|B\|_1+\left(2^{\frac{2-p}{2}}
\kappa_1\right)^{p'}\right)^{\frac{2-p}{p}}\nonumber\\
&\leq  2\left(\tfrac{\|B\|_1}{p-1}+\kappa_1^{p'}\right)^{\frac{2-p}{p}}.\end{align}
 By combining (\ref{uniqueness1_1}), (\ref{visco_convective3_1}) and (\ref{tensor_suplem_1}), we deduce that
	$$\left(\tfrac{1}{2\left(\frac{\|B\|_1}{p-1}+\kappa_1^{p'}\right)^{\frac{2-p}{p}}}-{\cal R}e\,\kappa_1\kappa_6
	\left(\tfrac{\|B\|_1}{p-1}+\kappa_1^{p'}\right)^{\frac{2-p}{p}}\right)
	\left\|D^\star\left({ \boldsymbol u}-{ \boldsymbol v}\right)\right\|_{p,B}^2\leq 0$$
and thus ${ \boldsymbol u}= { \boldsymbol v}$ if condition (\ref{restriction_control}) is satisfied.$\hfill\Box$\vspace{2mm}\\
{\bf Proof of Corollary \ref{pressure_1_thin}.} Arguing as in the shear-thickening case (cf. the proof of Corollary \ref{pressure_1}), we can prove that 
	$$\left\|\pi\right\|_{p'}\leq C\left(\left\|D u\right\|_{p,B}^2+
	\delta \left\|\nabla u_3\right\|_{p,B}^2+
	\left\|\tau \left(D^\star \boldsymbol u\right)\right\|_{p'}+
	\left\|\nabla^\star\cdot\tau \left(D^\star \boldsymbol u\right)-
	\nabla\cdot\tau \left(D^\star \boldsymbol u\right)\right\|_{p'}\right),$$
where $\tau=(\tau_{ij})_{i,j=1,2}$.
On the other hand, by taking into account (\ref{estim_S5}) and (\ref{diff_div_thin}), we obtain
	$$\left\|\tau \left(D^\star \boldsymbol u\right)\right\|_{p'}= 2\left\|\left(1+|D^\star \boldsymbol u|^2\right)^{\frac{p-2}{2}}Du\right\|_{p'}\leq 2 \left\|Du\right\|_p^{p-1},$$
and 
	\begin{equation}\label{div_tau_pres}\left\|\nabla^\star\cdot\tau \left(D^\star \boldsymbol u\right)-
	\nabla\cdot\tau 
	\left(D^\star \boldsymbol u\right)\right\|_{p'}\leq 4\delta m 
	\left(2\left\|Du\right\|_p^{p-1}+
	\left\|\tfrac{\delta}{B}u_2\right\|_p^{p-1}\right).\end{equation}
The conclusion follows by combining the three inequalities.\hfill\Box
\subsection{$\delta$-approximation}
{\bf Proof of Proposition \ref{w_exist_thin}.} Even though the idea is similar to the one used in the shear-thickening case, the lack of coercivity of the stress tensor when splitting the system and considering the equations for $(0,0,w_3)$ and $(w_1,w_2,0)$ generates an additional difficulty. Setting $\phi=(0,0,w_3)$ in the corresponding weak formulation, we obtain
	$$\int_\Sigma\left(1+|D\boldsymbol w|^2\right)^{\frac{p-2}{2}}
	|\nabla w_3|^2\,dx=\left(\tfrac{G}{B},w_3\right).$$
Let
	$$H_1=\left(1+|D\boldsymbol w|^2\right)^{\frac{2-p}{2}},\qquad 
	H_2=\left(1+|D\boldsymbol w|^2\right)^{\frac{p-2}{2}}
	|\nabla w_3|^2,
	\qquad H_3=|\nabla w_3|.$$
Due to Lemma \ref{desigualdade3} , we have
	$$\tfrac{\left\|\nabla w_3\right\|^2_p}{\left(|\Sigma|+
	\left\|D\boldsymbol w\right\|_p^p\right)^{\frac{2-p}{p}}}\leq 
	\left(\tfrac{G}{B},w_3\right).$$
On the other hand, standard arguments together with the Poincar\'e inequality (\ref{poincare}) and the Korn inequality give
	$$\left|\left(\tfrac{G}{B},w_3\right)\right|\leq c_1
	\left\|\nabla w_3\right\|_p,$$
where $c_1=m|G||\Sigma|^{\frac{1}{p'}}$.
Combining these inequalities yields
	\begin{equation}\label{dw3_thin}
	\tfrac{\left\|\nabla w_3\right\|_p}{\left(|\Sigma|+
	\left\|D\boldsymbol w\right\|_p^p\right)^{\frac{2-p}{p}}}
	\leq c_1.\end{equation} 
Similarly, by setting $\boldsymbol\varphi=(w,0)=(w_1,w_2,0)$ in the corresponding weak formulation, we obtain
	$$2\int_\Sigma\left(1+|D\boldsymbol w|^2\right)^{\frac{p-2}{2}}
	|D w|^2\,dx=\delta \left(\sigma(w_3),w_2\right)$$
and thus
	$$\tfrac{\|Dw\|_{p}^2}{\left(|\Sigma|+
	\left\|D\boldsymbol w\right\|_p^p\right)^{\frac{2-p}{p}}}\leq \tfrac{\delta}{2}\left(\sigma(w_3),w_2\right).$$
Estimate (\ref{est_psi_1_thin}) together with the Korn inequality give
	$$\begin{array}{ll}
	\tfrac{\|Dw\|_{p}^2}{\left(|\Sigma|+
	\left\|D\boldsymbol w\right\|_p^p\right)^{\frac{2-p}{p}}}&\leq  \tfrac{D_{p,\alpha}}{2}\, c_0\delta
	\left\|\nabla w_3\right\|_p^\alpha
	\left\|\nabla w_2\right\|_p\vspace{2mm}\\
	&\leq\tfrac{D_{p,\alpha}}{2}\, c_0\delta
	\left\|\nabla w_3\right\|_p^\alpha
	\left\|\nabla w\right\|_p\leq 
	c_2 c_0\delta
	\left\|\nabla w_3\right\|_p^\alpha
	\left\|Dw\right\|_p\end{array}$$
where $c_2=\tfrac{D_{p,\alpha}}{2C_{K,1}}$. Consequently, we have
	\begin{equation}\label{dw_thin}\tfrac{\|Dw\|_{p}}{\left(|\Sigma|+
	\left\|D\boldsymbol w\right\|_p^p\right)^{\frac{2-p}{p}}}\leq 
	c_0c_2\, \delta
	\left\|\nabla w_3\right\|_p^\alpha.
	\end{equation}
Combining (\ref{dw3_thin}) and (\ref{dw_thin}), it follows that
	$$\begin{array}{ll}
	\tfrac{\left\|D\boldsymbol w\right\|_p^p}{\left(|\Sigma|+
	\left\|D\boldsymbol w\right\|_p^p\right)^{2-p}}
	&\leq c_1^p+\left(c_0c_2\right)^p
	\left\|\nabla w_3\right\|_p^{\alpha p}\vspace{2mm}\\
	&\leq 
	c_1^p+\left(c_0c_1^{\alpha}c_2\right)^p
	\left(|\Sigma|+
	\left\|D\boldsymbol w\right\|_p^p\right)^{(2-p)\alpha}
	\end{array}$$
yielding to
	$$\begin{array}{ll}\left\|D\boldsymbol w\right\|_p^p&\leq 
	c_1^p \left(|\Sigma|+
	\left\|D\boldsymbol w\right\|_p^p\right)^{2-p}+
	\left(c_0c_1^{\alpha}c_2\right)^p
	\left(|\Sigma|+
	\left\|D\boldsymbol w\right\|_p^p\right)^{(2-p)(\alpha+1)}
	\vspace{2mm}\\
	&\leq c_3^p	\left(1+
	\left\|D\boldsymbol w\right\|_p^p\right)^{(2-p)(\alpha+1)}
	\end{array}$$
and 
	\begin{equation}\label{Dw_global_1}
	\left\|D\boldsymbol w\right\|_p^{\frac{p}{(2-p)(\alpha+1)}}
	\leq c_3^{\frac{p}{(2-p)(\alpha+1)}}
	\left(1+\left\|D\boldsymbol w\right\|_p^{p}\right),
	\end{equation}
where $c_3=\left(c_1+
	c_0c_1^{\alpha}c_2 \right)
	\left(1+|\Sigma|\right)^{\frac{(2-p)(\alpha+1)}{p}}$. By using the Young inequality, we deduce that for $\alpha<\tfrac{p-1}{2-p}$ we have
	\begin{align}\label{Dw_global_2}
	c_3^{\frac{p}{(2-p)(\alpha+1)}}\left\|D\boldsymbol w\right\|_p^{p} &\leq (2-p)(\alpha+1)
	 \left\|D\boldsymbol u\right\|_p^{\frac{p}{(2-p)(\alpha+1)}}
	\nonumber\\
	&+\left(1-(2-p)(\alpha+1)\right)
	c_3^{\frac{p}{\left(1-(2-p)(\alpha+1)\right)(2-p)(\alpha+1)}}. \end{align}
Combining (\ref{Dw_global_1}) and (\ref{Dw_global_2}), we obtain
	\begin{equation}\label{dw3w_thin}
	\left\|D\boldsymbol w\right\|_p\leq
	 c_3\left(\tfrac{1}{1-(2-p)(\alpha+1)}+
	c_3^{\frac{1}{1-(2-p)(\alpha+1)}}
	\right)^{\frac{(2-p)(\alpha+1)}{p}}.\end{equation}
The conclusion follows from (\ref{dw3_thin}), (\ref{dw_thin}) and 
(\ref{dw3w_thin}). $\hfill\Box$\vspace{2mm}\\
{\bf Proof of Proposition \ref{delta_vs_zero_thin}.} The ideas of the proof are similar to the ones used in the shear-thickening case. Indeed, by arguing as in the proof of Proposition \ref{delta_vs_zero}, we obtain
\begin{align}\label{difference1_thin}\left(\boldsymbol\tau(D{ \boldsymbol u})-{\boldsymbol\tau}\left(D  { \boldsymbol w}\right),D \boldsymbol\varphi\right)&=\left(\boldsymbol\tau(D{ \boldsymbol u})-{\boldsymbol\tau}\left(D^\star{ \boldsymbol u}\right),D \boldsymbol\varphi\right)+
\left(\nabla^\star\cdot \boldsymbol\tau(D^\star{ \boldsymbol u})-\nabla\cdot \boldsymbol\tau(D^\star{ \boldsymbol u}),\boldsymbol\varphi\right)\nonumber\\
	&-{\cal R}e\left(a_\star\left(\boldsymbol u, 
	  \boldsymbol u,\boldsymbol\varphi\right)-a\left({ \boldsymbol w},{ \boldsymbol w},{\boldsymbol\varphi}\right)\right)+\left(\pi_1-\pi_2,\nabla\cdot \boldsymbol\varphi\right)-
	\delta\left(\sigma(w_3),\varphi_2\right)\nonumber\\
	&=I_1+I_2+I_3+I_4+I_5.\end{align}
$\circ$ By taking into account (\ref{estim_S7_2}), we have
\begin{align}\label{difference_est3_thin}\left|I_1\right|&\leq \left\|{\boldsymbol\tau }\left(D^\star{ \boldsymbol u}\right)-{\boldsymbol\tau }\left(D  { \boldsymbol u}\right)\right\|_{p'}
\left\|D\boldsymbol\varphi\right\|_{p}\nonumber\\
	&\leq	2C_{p}\left\|D^\star{ \boldsymbol u}-D  { \boldsymbol u}\right\|_{p}^{p-1}
	\left\|D\boldsymbol\varphi\right\|_{p}\leq F_ 1 \,\delta^{p-1} \left\|D\boldsymbol\varphi\right\|_{p},
	\end{align}
where $F_1=2C_{p}m ^{p-1}\left(\|u_2\|_p+\|u_3\|_p\right)^{p-1}$. \vspace{2mm}\\
$\circ$ Estimate (\ref{diff_div_thin})  together with the Sobolev inequality (\ref{sobolev0}) and the Korn inequality (\ref{korn-inequality}) yield 
	\begin{align}\label{difference_est4_thin}\left|I_2\right|&\leq 
	\left\|\nabla^\star\cdot \boldsymbol\tau (D^\star{ \boldsymbol u})-\nabla\cdot \boldsymbol\tau (D^\star{ \boldsymbol u})\right\|_{p'}
	\left\|\boldsymbol\varphi\right\|_p\nonumber\\
	&\leq 4mS_{p,p}
	\left(\left\|Du\right\|_{p}^{p-1}
	+\left\|D_{33}^\star{ \boldsymbol u}\right\|_{p}^{p-1}+\left\|D_{23}^\star{ \boldsymbol u}\right\|_{p}^{p-1}\right)\, \delta 
	\left\|\nabla \boldsymbol\varphi\right\|_{p}\nonumber\\
	 &\leq 12mS_{p,p}
	\left\|D^\star{ \boldsymbol u}\right\|_{p}^{p-1}\, \delta 
	\left\|\nabla \boldsymbol\varphi\right\|_{p}\leq  F_2\, \delta
	\left\|D\boldsymbol\varphi\right\|_{p},\end{align}
where $F_2=\tfrac{12m}{C_{K,1}} S_{p,p}\, \left\|D^\star{ \boldsymbol u}\right\|_{p}^{p-1}$. \vspace{2mm}\\
$\circ$ The convective term is estimated similarly
\begin{align}\label{difference_est61_thin}\tfrac{1}{{\cal R}e}\left|I_3\right|
&=\left|a_\star({ \boldsymbol u},{ \boldsymbol u},\boldsymbol \varphi)-a({ \boldsymbol w},{ \boldsymbol w},\boldsymbol \varphi)\right|\nonumber\\
&\leq\left|a_\star\left({ \boldsymbol u}-{ \boldsymbol w},{ \boldsymbol u},\boldsymbol \varphi\right)\right|+\left|a\left(\boldsymbol w,{ \boldsymbol u}-{ \boldsymbol w},\boldsymbol \varphi\right)\right|+
\displaystyle\left|\left(\boldsymbol w\cdot \nabla^\star \boldsymbol u-\boldsymbol w\cdot \nabla \boldsymbol u,\boldsymbol \varphi\right)\right|\nonumber\\
&\leq\left|a_\star\left({ \boldsymbol u}-{ \boldsymbol w},{ \boldsymbol u},\boldsymbol \varphi\right)\right|+\left|a\left(\boldsymbol w,{ \boldsymbol u}-{ \boldsymbol w},\boldsymbol \varphi\right)\right|+
\displaystyle\left\|\tfrac{\delta}{B}|\boldsymbol w||\boldsymbol u||\boldsymbol \varphi|\right\|_1\nonumber\\
&\leq 	\left\|{ \boldsymbol u}-{ \boldsymbol w}\right\|_{2p'}
	\left\|\nabla^\star{ \boldsymbol u}\right\|_p\left\|\boldsymbol \varphi\right\|_{2p'} +\left\|{ \boldsymbol w}\right\|_{2p'}
	\left\|\nabla\left({ \boldsymbol u}-{ \boldsymbol w}\right)\right\|_p\left\|\boldsymbol \varphi\right\|_{2p'} + \delta m \left\|\boldsymbol w\right\|_{2p'} \left\|\boldsymbol u\right\|_{2p'}
	\left\|\boldsymbol \varphi\right\|_p\nonumber\\
&\leq  \left(1+S_{p,p}\right)\left( S_{p,2p'}\right)^2
\left(\|\nabla\left({ \boldsymbol u}-{ \boldsymbol w}\right)\|_p\left(\|\nabla^\star{ \boldsymbol u}\|_p+\|\nabla{ \boldsymbol w}\|_p\right)+\delta m\|\nabla{ \boldsymbol w}\|_p\|\nabla{ \boldsymbol u}\|_p\right)
	\left\|\nabla \boldsymbol \varphi\right\|_p\nonumber\\
	&\leq F_3\left(\delta
	+\left\|D\left({ \boldsymbol u}-{ \boldsymbol w}\right)\right\|_{p}\right)\left\|D \boldsymbol\varphi\right\|_p\end{align}
with $F_3=\tfrac{1+m}{C_K^3}\left(1+S_{p,p}\right)\left( S_{p,2p'}\right)^2\left(\left\|D^\star \boldsymbol u\right\|_p+\left\|D\boldsymbol w\right\|_p+\left\|D\boldsymbol u\right\|_p\left\|D\boldsymbol w\right\|_p\right)$, and where $C_K$ is the Korn constant given in $(\ref{korn-inequality})$.\vspace{2mm}\\
$\circ$ The estimate associated to $\left|I_4\right|$ may be obtained with slight modifications in the proof given in Corollary \ref{pressure_1} by observing that 
$$\begin{array}{ll}\left\|\pi_1-\pi_2\right\|_{p'}
	&\leq   \left\|\tau (D^\star{ \boldsymbol u})-\tau (D{ \boldsymbol w})\right\|_{p'}+
	\left\|\nabla^\star\cdot\tau \left(D^\star \boldsymbol u\right)-
	\nabla^\star\cdot\tau
	 \left(D^\star \boldsymbol u\right)\right\|_{p'}
	\vspace{2mm}\\
	 &+{\cal R}e
	 \left(\left\|u\otimes u-w\otimes w\right\|_{p'}+
	\delta   \left\|u_3^2\right\|_{p'}\right)+\delta \left\|\sigma(w_3)\right\|_{p'}.\end{array}$$
Taking into account (\ref{div_tau_pres}) and (\ref{estim_S7_2}), we have
	$$\left\|\nabla^\star\cdot\tau \left(D^\star \boldsymbol u\right)-
	\nabla^\star\cdot\tau \left(D^\star \boldsymbol u\right)\right\|_{p'}\leq F_{4,1}\,\delta$$
and
$$\begin{array}{ll}\left\|\tau (D^\star{ \boldsymbol u})-\tau (D{ \boldsymbol w})\right\|_{p',B}&\leq 
	\left\|\tau (D^\star{ \boldsymbol u})-\tau (D{ \boldsymbol w})\right\|_{p'}
	\leq 2C_{p}\left\|D^\star{ \boldsymbol u}-D{ \boldsymbol w}\right\|_{p}^{p-1}\vspace{2mm}\\
	&\leq 2C_{p}\left(\left\|D^\star{ \boldsymbol u}-D{ \boldsymbol u}\right\|_{p}^{p-1}
	+\left\|D{ \boldsymbol u}-D{ \boldsymbol w}\right\|_{p}^{p-1}\right)\vspace{2mm}\\
	&\leq F_{4,2}\left(\left\|D{ \boldsymbol u}-D{ \boldsymbol w}\right\|_{p}^{p-1}
	+\delta^{p-1}\right),\end{array}$$
where $F_{4,2}= 2C_{p}+F_1$. Moreover,  by using the Sobolev inequality (\ref{sobolev0}) with $(q,r)=(2,4)$ and the classical Korn inequality , we obtain
$$\begin{array}{ll}\left\| u\otimes u-w\otimes w\right\|_{p'}+\delta \left\|u_3^2\right\|_{p'}

	&\leq \left\|\left(u+w\right)

	\otimes \left( u-w\right)\right\|_{p'}+\delta \left\|u_3\right\|^2_{2p'}\vspace{2mm}\\

	&\leq 	\left\|u+w\right\|_{2p'}

	\left\|u-w\right\|_{2p'}+\delta \left\|u_3\right\|^2_{2p'}\vspace{2mm}\\

	&\leq  S_{p,2p'} 

	\left\|\nabla\left(u-w\right)

	\right\|_p\left\|u+w\right\|_{2p'}+\delta\, \left\|u_3\right\|^2_{2p'}\vspace{2mm}\\

&\leq   S_{p,2p'} 

	\left\|\nabla\left(\boldsymbol u-\boldsymbol w\right)
	\right\|_p\left\|u+w\right\|_{2p'}+\delta\, \left\|u_3\right\|^2_{2p'}\vspace{2mm}\\
&\leq \tfrac{S_{p,2p'}}{C_K} 
	\left\|D\left(\boldsymbol u-\boldsymbol w\right)
	\right\|_p\left\|u+w\right\|_{2p'}+\delta\, \left\|u_3\right\|^2_{2p'}\vspace{2mm}\\
	&\leq F_{4,3} \left(
	\left\|D\left(\boldsymbol u-\boldsymbol w\right)\right\|_p+\delta\right),\end{array}$$
where $F_{4,3}=\tfrac{S_{p,2p'}}{C_{K,1}}  \left\|u+w\right\|_{2p'}+ \left\|u_3\right\|^2_{2p'}$. Due to (\ref{est_psi_2_thin}), we have
	$$\left\|\sigma(w_3)\right\|_{p'}\leq c_0 E_{\alpha,p}
	 \left\|\nabla w_3\right\|_{p}^\alpha=F_{4,4}.$$
Combining these estimates, we deduce that
	\begin{align}\label{difference_est8_thin}\left|I_4\right|&\leq \left(F_{4,2}\left\|D{ \boldsymbol u}-D{ \boldsymbol w}\right\|_{p}^{p-1}+{\cal R}e\, F_{4,3}
\left\|D{ \boldsymbol u}-D{ \boldsymbol w}\right\|_{p}\right) \left\|\nabla \cdot \boldsymbol\varphi\right\|_p\nonumber\\
&+\left( \left(F_{4,1}+{\cal R}e\,F_{4,3}\right)\delta+F_{5,2}\delta^{p-1}\right) \left\|\nabla \cdot \boldsymbol\varphi\right\|_p\nonumber\\
	&\leq F_ 4\left(\left\|D\left({ \boldsymbol u}-{ \boldsymbol w}\right)\right\|_{p}
	+\left\|D\left({ \boldsymbol u}-{ \boldsymbol w}\right)\right\|_{p}^{p-1}+ 
	\delta^{p-1}\right)\left\|\nabla \cdot \boldsymbol\varphi\right\|_p.\end{align}
$\circ$ Finally, taking into account (\ref{est_psi_1_thin}) we obtain
	\begin{align}\label{difference_est_psi_thin}\left|I_5\right|&=\delta \left|\left(\sigma(w_3),
	\varphi_2\right)\right|\leq D_{p,\alpha} \,\delta
	\left\|\nabla w_3\right\|^\alpha_p
	\left\|\nabla\varphi_2\right\|_{p}\nonumber\\
	&\leq D_{p,\alpha}\,\delta
	\left\|\nabla w_3\right\|^\alpha_2
	\left\|\nabla\boldsymbol\varphi
	\right\|_{p}\leq F_5\,\delta
	\left\|D\boldsymbol\varphi
	\right\|_{p},\end{align}
where $F_5=\tfrac{D_{p,\alpha}}{C_{K,1}} 
	\left\|\nabla w_3\right\|^\alpha_p$. \vspace{2mm}\\
$\circ$ Combining (\ref{difference_est3_thin})-(\ref{difference_est_psi_thin}), and taking into the estimates associated to 
$\boldsymbol u$ and $\boldsymbol w$, we deduce that
	$$\begin{array}{ll}&\left|\left(\boldsymbol\tau (D{ \boldsymbol u})-{\boldsymbol\tau }\left(D  { \boldsymbol w}\right),D \boldsymbol\varphi\right)\right|\vspace{2mm}\\
	&\leq	\left(F_1\, \delta^{p-1}+F_2\, \delta+{\cal R}e\,F_3\, \delta+F_5\, \delta\right)\left\|D\left({ \boldsymbol u}-{ \boldsymbol w}\right)\right\|_{p}+{\cal R}e\,F_3 \left\|D \left(\boldsymbol u-\boldsymbol w\right)\right\|_{p}\left\|D\boldsymbol\varphi\right\|_p\vspace{2mm}\\
	& +F_ 4\left(\left\|D\left({ \boldsymbol u}-{ \boldsymbol w}\right)\right\|_{p}
	+\left\|D\left({ \boldsymbol u}-{ \boldsymbol w}\right)\right\|_{p}^{p-1}+ 
	\delta^{p-1}\right)\left\|\nabla \cdot \boldsymbol\varphi\right\|_p \vspace{2mm}\\
	&\leq	C_1\left(\delta^{p-1}\left\|D^\star\left({ \boldsymbol u}-{ \boldsymbol w}\right)\right\|_{p,B}+
	\delta\left\|D^\star\left({ \boldsymbol u}-{ \boldsymbol w}\right)\right\|_{p,B}^{p-1}+
	{\cal R}e\left\|D^\star \left(\boldsymbol u-\boldsymbol w\right)\right\|_{p,B}^2+\delta^p\right),\end{array}$$
where $C_1$ depends only on $p$, $\Sigma$, $m$, $n$, $\alpha$ and $c_0$. Setting $\boldsymbol\varphi=\boldsymbol u-\boldsymbol w$ and taking into account the estimates associated to $\boldsymbol u$ and $\boldsymbol w$, we deduce that
	$$\begin{array}{ll}&
	\left({\boldsymbol\tau }\left(D{ \boldsymbol u}\right)
	-{\boldsymbol\tau }\left(D{ \boldsymbol w}\right),
	D\left({ \boldsymbol u}-{ \boldsymbol w}\right)\right)
	\vspace{2mm}\\
	&\leq \left(F_1\, \delta^{p-1}+F_2\, \delta+{\cal R}e\,F_3\, \delta+F_5\, \delta\right)\left\|D\left({ \boldsymbol u}-{ \boldsymbol w}\right)\right\|_{p}+{\cal R}e\,F_3 \left\|D \left(\boldsymbol u-\boldsymbol w\right)\right\|_{p}^2\vspace{2mm}\\
	& +F_4\left(\left\|D\left({ \boldsymbol u}-{ \boldsymbol w}\right)\right\|_{p}
	+\left\|D\left({ \boldsymbol u}-{ \boldsymbol w}\right)\right\|_{p}^{p-1}+ 
	\delta^{p-1}\right)\left\|\tfrac{\delta}{B}u_2\right\|_p \vspace{2mm}\\
	&\leq	C_2\left(\delta^{p-1}\left\|D\left({ \boldsymbol u}-{ \boldsymbol w}\right)\right\|_{p}+
	\delta\left\|D\left({ \boldsymbol u}-{ \boldsymbol w}\right)\right\|_{p}^{p-1}+
	{\cal R}e\left\|D \left(\boldsymbol u-\boldsymbol w\right)\right\|_{p}^2+\delta^p\right),\end{array}$$
where $C_2$ depends only on $p$, $\Sigma$, $m$, $n$, $\alpha$ and $c_0$. Using the Young inequality,  it follows that for every $\varepsilon>0$, we have
	\begin{align}\label{est_fin}&\left({\boldsymbol\tau }\left(D{ \boldsymbol u}\right)
	-{\boldsymbol\tau }\left(D{ \boldsymbol w}\right),
	D\left({ \boldsymbol u}-{ \boldsymbol w}\right)\right)\nonumber\\
	& \leq 
	\left(\varepsilon+C_2\,{\cal R}e\right) \left\|D^\star{ \boldsymbol u}-D^\star{ \boldsymbol w}\right\|_{p,B}^2+C_3(\varepsilon)\left(\delta^{2(p-1)}+\delta^{\frac{2}{3-p}}\right)+C_1\,\delta^p.\end{align}
On the other hand, by taking into account $(\ref{estim_S7})$ and the estimates associated to $\boldsymbol u$ and $\boldsymbol w$, we deduce that there exists a constant $C_4$ depending on $p$, $\Sigma$, $G$ and $m$, $\alpha$ and $c_0$, but independent of $\delta$, such that
	\begin{equation}\label{difference7_thin}
	\left({\boldsymbol\tau }\left(D{ \boldsymbol u}\right)-{\boldsymbol\tau }\left(D{ \boldsymbol w}\right),
	D\left({ \boldsymbol u}-{ \boldsymbol w}\right)\right)\geq \tfrac{2\left\|D\left({ \boldsymbol u}-{ \boldsymbol w}\right)\right\|_{p}^2}{
	\left(|\Sigma|+\|D{ \boldsymbol u}\|_{p}^p+\|D{ \boldsymbol w}\|_{p}\right)^{\frac{2-p}{p}}}
	\geq C_4 \left\|D\left({ \boldsymbol u}-{ \boldsymbol w}\right)\right\|_{p}^2.\end{equation}
Combining (\ref{est_fin}) and (\ref{difference7_thin}),  observing that $\delta^{p}<\delta^{2(p-1)}$ and $\delta^{\frac{2}{3-p}}<\delta^{2(p-1)}$, choosing $\varepsilon= \tfrac{C_4}{2}$ and assuming that ${\cal R}e<\tfrac{C_4}{2C_2}$, we deduce that 
	$$\left\|D\left(\boldsymbol u-\boldsymbol w\right)\right\|_{p,B}\leq 
	C_5\,\delta^{p-1}$$
and the claimed result is proved. $\hfill\Box$

\appendix
\section{Toroidal coordinate system} \label{appendix_a}
\setcounter{equation}{0}
Let us consider the new coordinate system, in the variables $
(\widetilde x_1,\widetilde x_2 ,\widetilde x_3)$, given by the 
transformations $(\widetilde x_1,\widetilde x_2,\widetilde x_3)\mapsto (\widetilde y_1,\widetilde y_2,\widetilde y_3)$ satisfying
(\ref{transf_1}) and (\ref{transf_2}).
Let $M$ be a generic point such that 
	$$M=\widetilde y_1(\widetilde x_1,\widetilde x_2 ,\widetilde x_3)\,\mathbf{e}_{1}
	+\widetilde y_2(\widetilde x_1,\widetilde x_2 ,\widetilde x_3)	\,\mathbf{e}_{2}+
	\widetilde y_3(\widetilde x_1,\widetilde x_2 ,\widetilde x_3)
	\,\mathbf{e}_{3}.$$
 Defining the scale vectors 
                $$h_{1}=\left| \tfrac{\partial M}{\partial \widetilde x_1}
		\right|, \qquad\quad h_{2}=\left|
		  \tfrac{\partial M}{\partial \widetilde x_2 }\right|,\quad\qquad  
		  h_{3}=\left| \tfrac{\partial M}{\partial \widetilde  x_3}
		  \right|$$
and the local basis 
	$$\mathbf{a}_{1}=\tfrac{1}{h_{1}}
	   \tfrac{\partial M}{\partial \widetilde x_1}, \quad\qquad 
	   \mathbf{a}_{2}=\tfrac{1}{h_{2}}
	   \tfrac{\partial M}{\partial \widetilde x_2},\quad
	    \qquad \mathbf{a}_{3}=\tfrac{1}{h_{3}}
	    \tfrac{\partial M}{\partial \widetilde  x_3},$$
we obtain, 
	$$h_{1}=1, \quad\qquad h_{2}=1, \quad\qquad 
	     h_{3}=1+\tfrac{1}{R }\widetilde x_2$$
and
	$$\mathbf{a}_{1}=\mathbf{e}_{3},\qquad\quad
	     \mathbf{a}_{2}=\cos \left(\tfrac{\widetilde x_3}{R }\right)
	     \mathbf{e}_{1}+\sin \left(\tfrac{\widetilde x_3}{R }\right)
	\mathbf{e}_{2},\quad\qquad
	     \mathbf{a}_{3}=-\sin \left(\tfrac{\widetilde x_3}{R }\right)\mathbf{e}_{1}
	     +\cos \left(\tfrac{\widetilde x_3}{R }\right)\mathbf{e}_{2}.$$
The corresponding derivatives of the vector basis are given by 
         $$\begin{array}{lll}
	 \frac{\partial \mathbf{a}_{1}}{\partial \widetilde x_1}=0, \quad  &
	  \quad \frac{\partial \mathbf{a}_{1}}{\partial \widetilde x_2}=0, \quad  &  \quad 
	  \frac{\partial \mathbf{a}_{1}}{\partial \widetilde x_3}=0
	\medskip , \\ 
	  \frac{\partial \mathbf{a}_2}{\partial \widetilde x_1}=0,
 \quad &
\quad	   \frac{\partial \mathbf{a}_2}{\partial \widetilde x_2 }=0, \quad &
	  \quad \frac{\partial \mathbf{a}_2}{\partial \widetilde x_3}
	  =\frac{\mathbf{a}_{3}}{R}\medskip , \\ 
	  \frac{\partial \mathbf{a}_{3}}{\partial \widetilde x_1}=0, \quad & \quad 
	  \frac{\partial \mathbf{a}_3}{\partial \widetilde x_2}=0,
	   \quad & \quad 
	\frac{\partial \mathbf{a}_3}{\partial \widetilde x_3}=-
	   \frac{\mathbf a_2}{R }.
	   \end{array}
	   $$
$\bullet$  The gradient of a scalar function $\widetilde \psi$ in the rectangular toroidal coordinates is then given by 
         \begin{equation}\label{grad_scalar}
	\widetilde\nabla \widetilde \psi =
	 \tfrac{1}{h_1}
	 \tfrac{\partial \widetilde \psi}{\partial \widetilde x_1}\, \mathbf a_1+
	 \tfrac{1}{h_2}\tfrac{\partial \widetilde \psi}
	 {\partial \widetilde x_2}\,\mathbf a_2+\tfrac{1}{h_3}
	 \tfrac{\partial \widetilde \psi}{\partial \widetilde x_3}\mathbf a_3
	=\tfrac{\partial \widetilde \psi}{\partial \widetilde x_1}\,\mathbf a_1+
	 \tfrac{\partial \widetilde \psi}{\partial \widetilde x_2}\,\mathbf a_2
	+ \tfrac{1}{B}
	 \tfrac{\partial \widetilde \psi}{\partial \widetilde x_3}\,\mathbf a_3,
	\end{equation}
where $B=1+\tfrac{1}{R}\, \widetilde x_2$. \vspace{2mm}\\
$\bullet$ The gradient of a vector  $\widetilde{\boldsymbol v}\equiv (\widetilde v_{1},\widetilde v_{2},\widetilde v_{3})$ is defined by
	$$\widetilde \nabla \widetilde{\boldsymbol v}=\left(\tfrac{1}{h_1}\mathbf a_1\otimes \tfrac{\partial }{\partial \widetilde x_1}+\tfrac{1}{h_2}
	   \mathbf a_2\otimes \tfrac{\partial }{\partial \widetilde x_2}+
	   \tfrac{1}{h_3}\mathbf a_3\otimes 
	   \tfrac{\partial }{\partial \widetilde x_3}\right)
	   \left(\mathbf a_1\widetilde v_1+\mathbf a_2 \widetilde v_2
	   +\mathbf a_3\widetilde v_3\right)$$
that is
	\begin{align}\label{gradient}\widetilde \nabla \widetilde{\boldsymbol v}&=\displaystyle\sum_{j=1}^3\sum_{i=1}^2 \tfrac{\partial \widetilde v_j}{
		\partial \widetilde x_i}\, \mathbf a_i\otimes \mathbf a_j
		+\tfrac{1}{RB}
	       \left(\widetilde v_2\,\mathbf a_3\otimes\mathbf a_3-
		\widetilde v_3\,
	       \mathbf a_3\otimes\mathbf a_2\right)+\tfrac{1}{B}\sum_{j=1}^3
	       \tfrac{\partial \widetilde v_j}{\partial \widetilde x_3} \,
	       \mathbf a_3\otimes\mathbf a_j\nonumber\\
		&=\left( \begin{array}{lll}
	    \frac{\partial \widetilde v_{1}}{\partial \widetilde x_1} \ & \
	     \frac{\partial \widetilde v_2}{\partial \widetilde x_1} \ & \
	      \frac{\partial \widetilde v_3}{\partial \widetilde x_1}\vspace{3mm}\\ 
	      \frac{\partial \widetilde v_1}{\partial \widetilde x_2} \ & \
	     \frac{\partial \widetilde v_2}{\partial \widetilde x_2} \ & \
	        \frac{\partial \widetilde v_3}{\partial \widetilde x_2 }\vspace{3mm} \\ 
		\frac{1}{B}\frac{\partial \widetilde v_{1}}{\partial \widetilde x_3} \ & \ 
		\frac{1}{B}\left(\frac{\partial \widetilde v_{2}}{\partial \widetilde x_3}
		-\frac{\widetilde v_3}{R}\right)\ & \ \frac{1}{B}\left(\frac{\partial \widetilde v_{3}}{\partial \widetilde x_3}
		+\frac{\widetilde v_2}{R}\right)
		\end{array}\right).\end{align}
$\bullet$  Similarly, if $\widetilde{\boldsymbol v}\equiv (\widetilde v_{1},\widetilde v_{2},\widetilde v_{3})$ and $\widetilde{\boldsymbol w}\equiv (
\widetilde w_{1},\widetilde w_{2},\widetilde w_{3})$ are two vectors, then the convective term is defined by\vspace{3mm}\\
          \lefteqn{\hspace{1cm} \widetilde {\boldsymbol v}\cdot \widetilde \nabla \widetilde{\boldsymbol w}= \left(\mathbf a_1 \widetilde v_1+
	   \mathbf a_2 \widetilde v_2+\mathbf a_3 \widetilde v_3\right)\cdot 
	   \Big(\sum_{i,j}\mathbf a_i\otimes \mathbf a_j
	   \big(\widetilde \nabla\widetilde {\boldsymbol w}\big)_{ij}\Big)
	   =\sum_{i}\mathbf a_i\Big(\sum_j \widetilde v_j
	    \big(\widetilde \nabla\widetilde {\boldsymbol w}\big)_{ji}\Big),}\vspace{3mm}\\
that is
	\begin{align}\label{convective}\widetilde {\boldsymbol v}\cdot \widetilde \nabla \widetilde{\boldsymbol w}&=\displaystyle \sum_{j=1}^3\sum_{i=1}^2 \widetilde v_i \tfrac{\partial \widetilde w_j}{\partial \widetilde x_i} \mathbf a_j+\tfrac{\widetilde v_3}{RB}\left(\widetilde w_2\mathbf a_3-
	   \widetilde w_3\mathbf a_2\right)+\tfrac{\widetilde v_3}{B}
	   \tfrac{\partial \widetilde {\boldsymbol w}}{\partial \widetilde x_3}\nonumber\\
	    &=\left( \begin{array}{l}
	   \widetilde v_1\tfrac{\partial \widetilde w_1}{\partial \widetilde x_1}+
	    \widetilde v_2\tfrac{\partial \widetilde w_1}{\partial \widetilde x_2} +\tfrac{\widetilde v_3}{B}
	    \tfrac{\partial \widetilde w_1}{\partial \widetilde x_3}
	   \vspace{3mm} \\ 
	    \widetilde v_1\tfrac{\partial \widetilde w_2}{\partial \widetilde x_1}+
	    \widetilde v_2\tfrac{\partial \widetilde w_2}{\partial \widetilde x_2} +\tfrac{\widetilde v_3}{B}
	    \tfrac{\partial \widetilde w_2}{\partial \widetilde x_3}
	   -\tfrac{1}{RB}\widetilde v_3\widetilde w_3\vspace{3mm} \\ 
	   \widetilde v_1\tfrac{\partial \widetilde w_3}{\partial \widetilde x_1}+
	    \widetilde v_2\tfrac{\partial \widetilde w_3}{\partial \widetilde x_2} +\tfrac{\widetilde v_3}{B}
	    \tfrac{\partial \widetilde w_3}{\partial \widetilde x_3}+
	  \tfrac{1}{RB}\widetilde v_3\widetilde w_2
	   \end{array}\right).\end{align}
$\bullet$ Finally, the divergence of a vector $\widetilde{\boldsymbol v}\equiv (\widetilde v_{1},\widetilde v_{2},\widetilde v_{3})$ is defined by
\begin{align}\label{div_vector}\widetilde \nabla \cdot \widetilde{\boldsymbol v}&=\tfrac{1}{h_1h_2h_3}\left( \tfrac{
	   \partial }{\partial \widetilde x_1}\left( h_2h_3\,\widetilde v_1\right) 
	   +\tfrac{\partial }{\partial \widetilde x_2}
	   \left( h_3h_1\,\widetilde v_2\right) 
	   +\tfrac{\partial }{\partial \widetilde x_3}
	\left( h_1h_2\,\widetilde v_3 \right)
	    \right)\nonumber\\
	&=\tfrac{\partial \widetilde v_1}{\partial \widetilde x_1}
	   +\tfrac{\partial \widetilde v_2}{\partial \widetilde x_2}+
	   \tfrac{\widetilde v_2}{RB}+
	   \tfrac{1}{B}\tfrac{\partial \widetilde v_3}{\partial \widetilde x_3}.\end{align}
 The divergence of a tensor $\widetilde{\mathbf S}$ is defined by \vspace{3mm}\\
         $$\widetilde\nabla\cdot\widetilde{\mathbf S}=
	    \left(\mathbf a_1\tfrac{\partial }{\partial \widetilde x_1}+
	   \mathbf a_{2}
	   \tfrac{\partial }{\partial \widetilde x_2}
		+\mathbf a_3\tfrac{1}{B}
	   \tfrac{\partial }{\partial \widetilde x_3}\right)\cdot 
	\Big(\sum_{i,j}
	   \mathbf a_i\otimes \mathbf a_j \widetilde S_{ij}\Big)$$
which gives
          \begin{align}\label{div_tensor}
	  \widetilde\nabla\cdot\widetilde{\mathbf S}
	&= \displaystyle \sum_{j=1}^3\sum_{i=1}^2 \tfrac{\partial \widetilde S_{ij}}{\partial \widetilde x_i}\, \mathbf a_j+	
	\tfrac{1}{B}\sum_{j=1}^3 \tfrac{\partial \widetilde S_{3j}}{\partial \widetilde x_3}\, \mathbf a_j+\tfrac{1}{RB}
	\left(\widetilde S_{21}\mathbf a_1+
	\left(\widetilde S_{22}-\widetilde S_{33}\right)
	\mathbf a_2+\left(\widetilde S_{32}+\widetilde S_{23}
	\right)\mathbf a_3\right)\nonumber\\
	&=\left(\begin{array}{ll}
	    \tfrac{\partial \widetilde S_{11}}{\partial \widetilde x_1}+
	    \tfrac{\partial \widetilde S_{21}}{\partial \widetilde x_2}+
	    \tfrac{1}{B}\tfrac{\partial
	    \widetilde S_{31}}{\partial \widetilde x_3}
	+\tfrac{\widetilde S_{21}}{RB}\vspace{3mm}\\
             \tfrac{\partial \widetilde S_{12}}{\partial 
	\widetilde x_1}+
	    \tfrac{\partial \widetilde S_{22}}{\partial \widetilde x_2}+
	    \tfrac{1}{B}\tfrac{\partial \widetilde S_{32}}{\partial \widetilde x_3}+
	    \tfrac{\widetilde S_{22}-\widetilde S_{33}}{RB}\vspace{3mm}\\
	    \tfrac{\partial \widetilde S_{13}}{\partial \widetilde x_1}+
	    \tfrac{\partial \widetilde S_{23}}{\partial \widetilde x_2}+
	    \tfrac{1}{B}\tfrac{\partial \widetilde S_{33}}{\partial \widetilde x_3}+
	    \tfrac{\widetilde S_{32}+\widetilde S_{23}}{RB}\end{array}\right).
	\end{align}

\section{Dimensionless system}\label{appendix_b}
\setcounter{equation}{0}
For the confort of the reader, the dimensionless equation (\ref{equation}) is derived hereafter. To simplify the notation, we consider the fully developed case and will assume that (\ref{fdev_velocity}) and (\ref{fdev_pressure}) are fulfilled. Taking into account (\ref{convective}), (\ref{grad_scalar}) and (\ref{adim_var}), we can easily see that
	$$\begin{array}{lll}
	 \widetilde {\boldsymbol u}\cdot \widetilde \nabla \widetilde{\boldsymbol u}&=
	   &\left( \begin{array}{l}
	   \widetilde u_1\tfrac{\partial \widetilde u_1}{\partial \widetilde x_1}+
	    \widetilde u_2\tfrac{\partial \widetilde u_1}{\partial \widetilde x_2} 
	   \vspace{3mm} \\ 
	    \widetilde u_1\tfrac{\partial \widetilde u_2}{\partial \widetilde x_1}+
	    \widetilde u_2\tfrac{\partial \widetilde u_2}{\partial \widetilde x_2}
	   -\tfrac{1}{RB}\widetilde u_3^2
	\vspace{3mm} \\ 
	   \widetilde u_1\tfrac{\partial \widetilde u_3}{\partial \widetilde x_1}+
	    \widetilde u_2\tfrac{\partial \widetilde u_3}{\partial \widetilde x_2}+
	  \tfrac{1}{RB}\widetilde u_3\widetilde u_2
	   \end{array}\right)\vspace{4mm}\\
	&=\tfrac{U_0^2}{r_0}&
	   \underbrace{\left( \begin{array}{l}
	   u_1\tfrac{\partial u_1}{\partial x_1}+
	    u_2\tfrac{\partial u_1}{\partial x_2} 
	   \vspace{3mm} \\ 
	    u_1\tfrac{\partial u_2}{\partial x_1}+
	    u_2\tfrac{\partial u_2}{\partial x_2}
	   -\tfrac{\delta}{B}u_3^2
	\vspace{3mm} \\ 
	   u_1\tfrac{\partial u_3}{\partial x_1}+
	   u_2\tfrac{\partial u_3}{\partial x_2}+
	 \tfrac{\delta}{B}u_3u_2
	   \end{array}\right)}_{{\boldsymbol u}\cdot \nabla^\star {\boldsymbol u}}
	\end{array}$$
and
	$$\widetilde\nabla \widetilde \pi =
	\left(\begin{array}{ll}
	\tfrac{1}{r_0}\tfrac{\partial \widetilde \pi }{\partial x_1}\vspace{3mm}\\
	 \tfrac{1}{r_0}\tfrac{\partial \widetilde \pi }{\partial x_2}\vspace{3mm}\\
	 -\tfrac{\widetilde G}{B}\end{array}\right)=
	\tfrac{\mu U_0}{r_0^2}\underbrace{\left(\begin{array}{ll}
	\tfrac{\partial \pi }{\partial x_1}\vspace{3mm}\\
	 \tfrac{\partial \pi }{\partial x_2}\vspace{3mm}\\
	 -\tfrac{G}{B}\end{array}\right)}_{\nabla^\star \pi}$$
with $B=1+\delta x_2$. 
Similarly, by taking into account (\ref{gradient}), we obtain
	$$\begin{array}{lll}\widetilde D \widetilde{\boldsymbol u}&=
	&\left( \begin{array}{lll}
	    \frac{\partial \widetilde u_{1}}{\partial \widetilde x_1} \ & \
	     \tfrac{1}{2}\left(\frac{\partial \widetilde u_2}{\partial \widetilde x_1}+
	\frac{\partial \widetilde u_1}{\partial \widetilde x_2}\right) \ & \
	       \tfrac{1}{2}\,\frac{\partial \widetilde u_3}{\partial \widetilde x_1}\vspace{3mm}\\ 
	         \tfrac{1}{2}\left(\frac{\partial \widetilde u_2}{\partial \widetilde x_1}+
	\frac{\partial \widetilde u_1}{\partial \widetilde x_2}\right)  \ & \
	     \frac{\partial \widetilde u_2}{\partial \widetilde x_2} \ & \
	        \tfrac{1}{2}\left(\frac{\partial \widetilde u_3}{\partial \widetilde x_2 }-\frac{\widetilde u_3}{RB}\right)
	\vspace{3mm} \\ 
		 \tfrac{1}{2}\,\frac{\partial \widetilde u_3}{\partial \widetilde x_1} \ & \ 
		\tfrac{1}{2}\left(\frac{\partial \widetilde u_3}{\partial \widetilde x_2 }-\frac{\widetilde u_3}{RB}\right)
		 \ & \ \frac{\widetilde u_2}{RB}
		\end{array}\right)\vspace{3mm}\\
	&=\tfrac{U_0}{r_0} &
	\underbrace{\left( \begin{array}{lll}    \frac{\partial u_{1}}{\partial x_1} \ & \
	     \tfrac{1}{2}\left(\frac{\partial u_2}{\partial x_1}+
	\frac{\partial u_1}{\partial x_2}\right) \ & \
	       \tfrac{1}{2}\,\frac{\partial u_3}{\partial x_1}\vspace{3mm}\\ 
	         \tfrac{1}{2}\left(\frac{\partial u_2}{\partial x_1}+
	\frac{\partial u_1}{\partial x_2}\right)  \ & \
	     \frac{\partial u_2}{\partial x_2} \ & \
	        \tfrac{1}{2}\left(\frac{\partial u_3}{\partial x_2 }-\frac{\delta}{B} u_3\right)
	\vspace{3mm} \\ 
		 \tfrac{1}{2}\,\frac{\partial u_3}{\partial x_1} \ & \ 
		\tfrac{1}{2}\left(\frac{\partial u_3}{\partial x_2 }-\frac{\delta}{B}u_3\right)
		 \ & \ \frac{\delta}{B}u_2
		\end{array}\right)}_{D^\star \boldsymbol u}.
	\end{array}$$
It follows that $\big|\widetilde D \widetilde{\boldsymbol u}\big|^2=\left(\tfrac{U_0}{r_0}\right)^2\left|D^\star \boldsymbol u\right|^2$ and
	$$\widetilde{\mathbf S}\big(\widetilde D \widetilde{\boldsymbol u}\big)\equiv 2\mu\left(1+\big|\widetilde D \widetilde{\boldsymbol u}\big|^2\right)^{\frac{p-2}{2}} \widetilde D \widetilde{\boldsymbol u}=
	\tfrac{\mu U_0}{r_0}\underbrace{2\left(1+\left(\tfrac{U_0}{r_0}\right)^2\left|D^\star \boldsymbol u\right|^2\right)^{\frac{p-2}{2}} D^\star \boldsymbol u}_{\boldsymbol\tau\left(D^\star \boldsymbol u\right)}.$$
Due to (\ref{div_tensor}), to the fact that $\widetilde{\mathbf S}$ is symmetric and that 
	$$\tfrac{\partial \widetilde S_{3i}}{\partial \widetilde x_3}=0 \qquad i=1,2,3$$
we deduce that
	$$\begin{array}{lll}
	   \widetilde\nabla \cdot  \big(\widetilde{\mathbf S}\big(\widetilde D \widetilde {\boldsymbol u}\big)\big) &= &\left(\begin{array}{ll}
	    \tfrac{\partial  \widetilde S_{11}}{\partial  \widetilde x_1}+
	    \tfrac{\partial  \widetilde S_{12}}{\partial  \widetilde x_2}+\tfrac{1}{RB}\widetilde S_{12}\vspace{3mm}\\
             \tfrac{\partial  \widetilde S_{12}}{\partial  \widetilde x_1}+
	    \tfrac{\partial  \widetilde S_{22}}{\partial  \widetilde x_2}+
	    \tfrac{1}{RB}\left(\widetilde S_{22}- \widetilde\tau_{33}\right)\vspace{3mm}\\
	    \tfrac{\partial  \widetilde S_{13}}{\partial  \widetilde x_1}+
	    \tfrac{\partial  \widetilde S_{23}}{\partial  \widetilde x_2}+
	    \tfrac{2}{RB}\widetilde S_{23}\end{array}\right)\vspace{3mm}\\ 
	&=\tfrac{1}{r_0}&\left(\begin{array}{ll}
	    \tfrac{\partial  \widetilde S_{11}}{\partial  x_1}+
	    \tfrac{\partial  \widetilde S_{12}}{\partial  x_2}+\tfrac{\delta}{B}\widetilde S_{12}\vspace{3mm}\\
             \tfrac{\partial  \widetilde S_{12}}{\partial  x_1}+
	    \tfrac{\partial  \widetilde S_{22}}{\partial  x_2}+
	    \tfrac{\delta}{B}\left(\widetilde S_{22}- \widetilde S_{33}\right)\vspace{4mm}\\
	    \tfrac{\partial  \widetilde S_{13}}{\partial  x_1}+
	    \tfrac{\partial  \widetilde S_{23}}{\partial x_2}+
	    \tfrac{2\delta}{B}	
	\widetilde S_{23}\end{array}\right)\vspace{3mm}\\ 
	&=\tfrac{\mu U_0}{r_0^2}&\underbrace{\left(\begin{array}{ll}
	    \tfrac{\partial  \tau_{11}}{\partial  x_1}+
	    \tfrac{\partial  \tau_{12}}{\partial  x_2}+\tfrac{\delta}{B}\tau_{12}\vspace{4mm}\\
             \tfrac{\partial  \tau_{12}}{\partial  x_1}+
	    \tfrac{\partial  \tau_{22}}{\partial  x_2}+
	    \tfrac{\delta}{B}\left(\tau_{22}- \tau_{33}\right)\vspace{3mm}\\
	    \tfrac{\partial  \tau_{13}}{\partial  x_1}+
	    \tfrac{\partial  \tau_{23}}{\partial x_2}+
	    \tfrac{2\delta}{B}	\tau_{23}\end{array}\right)}_{
	\nabla^\star\cdot\left(\boldsymbol\tau\left(D^\star \boldsymbol u\right)\right)},\end{array}$$
where we dropped the dependance on $\widetilde D \widetilde{\boldsymbol u}$ and $D^\star \boldsymbol u$. Taking into account these identities and substituting in equation $(\ref{equation_dim})_1$, we obtain
	$$\tfrac{\rho U_0^2}{r_0}\, \boldsymbol u\cdot \nabla^\star \boldsymbol u +\tfrac{\mu U_0}{r_0^2} \, \nabla^\star \pi=
	\tfrac{\mu U_0}{r_0^2} \, \nabla^\star\cdot\left(\boldsymbol\tau\left(D^\star \boldsymbol u\right)\right)
	 $$
which, by multiplying by $\tfrac{r_0^2}{\mu U_0}$, gives equation $(\ref{equation})_1$. Finally, by taking into account (\ref{div_vector}), we obtain 
	$$\widetilde \nabla \cdot \widetilde{\boldsymbol u}=\tfrac{\partial \widetilde u_1}{\partial \widetilde x_1}
	   +\tfrac{\partial \widetilde u_2}{\partial \widetilde x_2}+
	   \tfrac{\widetilde u_2}{RB}=\tfrac{U_0}{r_0}
	\underbrace{\left(\tfrac{\partial u_1}{\partial x_1}
	   +\tfrac{\partial u_2}{\partial x_2}+
   \tfrac{\delta}{B}\, u_2\right)}_{\nabla^\star \cdot \boldsymbol u},$$
showing that $(\ref{equation_dim})_2$ implies $(\ref{equation})_2$.


\begin{thebibliography}{10} 

\bibitem{amrouche} C. Amrouche, V. Girault, Decomposition of vector spaces and application to the Stokes problem in arbitrary dimension.
 Czechoslovak Mathematical Journal, 44 (1994) 109-140.

\bibitem{berger} A. A. Berger, L. Talbot, L.-S. Yao,  Flow in curved pipes, Ann. Rev. Fluid Mech., 15 (1983) 461-512.

\bibitem{Rob-coscia} V. Coscia, A. M. Robertson, Existence and uniqueness of steady, fully developed flows of second order fluids in curved
pipes, Math. Models Methods Appl. Sci. 6 (1) (2001) 1055-1071.

\bibitem{daska} P. Daskopoulos, A. M. Lenhoff, Flow in curved ducts:
bifurcation structure for stationary ducts, J. Fluid Mech., 203 (1989) 125-148.

\bibitem{dean1} W. R. Dean, Note on the motion of fluid in curved pipe,
Philos. Mag., 20 (1927) 208-223.

\bibitem{dean2} W. R. Dean, The streamline motion of fluid in curved pipe,
Philos. Mag., 30 (1928) 673-695.

\bibitem{eustice1} J. Eustice, Flow of water in curved pipes, Proc. R. Soc. Lond. A 84 (1910) 107-118.

\bibitem{eustice2} J. Eustice, Experiments of streamline motion in curved pipes, Proc. R. Soc. Lond. A 85 (1911) 119-131.

\bibitem{fan} Y. Fan, R. I. Tanner, N. Phan-Thien, Fully developed viscous and viscoelastic flows in curved pipes, 
 J. Fluid Mech., 440 (2001) 327-357.

\bibitem{galdirobertson} G. P. Galdi, A. M. Robertson, On flow of a Navier-Stokes fluid in curved pipes. Part I: Steady flow, Applied Math. Letters, 18 (2005) 1116-1124.

\bibitem{galdi} G. P. Galdi, An introduction to the mathematical theory of the Navier-Stokes equations, Vol. I and II, Springer Tracts in Natural Philosophy 38, 39, 2nd edition, Springer-Verlag, New York, 1998.

\bibitem{ito} H. Ito, Flow in curved pipes, JSME Int. J., 30 (1987) 543-552.

\bibitem{Rob-jichote} W. Jitchote, A. M. Robertson, Flow of second order fluids in curved pipes, J. Non-Newtonian Fluid Mech. 90, 2000,  91-116.

\bibitem{lady1} { O. A. Ladyzhenskaya}, New equations for the description of motion of viscous incompressible fluids and solvability in the large of boundary value problems for them, Proc. Stek. Inst. Math. 102 (1967) 95-118.

\bibitem{lady2} { O. A. Ladyzhenskaya}, On some modifications of the Navier-Stokes equations for large gradients of Velocity, Zap. Nau\v{c}n. Sem. Leningrad. Otdel. Mat. Inst. Steklov (LOMI) 7 (1968) 126-154.

\bibitem{lady} { O. A. Ladyzhenskaya}, The mathematical theory of viscous incompressible flow, Gordon and Beach, New York, 1969.

\bibitem{lions} { J.-L. Lions}, Quelques m\'ethodes de r\'esolution des probl\`emes aux limites non lin\'eaires, Dunod, Gauthier-Villars, Paris, 1969.

\bibitem{necas} {\sc J. Ne\v{c}as, J. M\'alek, J. Rokyta, M. Ru\v{z}i\v{c}ka},   Weak and measure-valued solutions to evolutionary partial differential equations, Applied Mathematics and Mathematical Computation, Vol. 13, Chapmann and Hall, London, 1996.

\bibitem{Rob} A. M. Robertson, On viscous flow in curved pipes of non-uniform cross section, Inter. J. Numer. Meth. fluid., 22 (1996) 771-798.

\bibitem{Rob-mul} A. M. Robertson, S. J. Muller, Flow of Oldroyd-B fluids in curved pipes of circular and annular cross-section, Int. J. Non-Linear  Mechanics, 31 (1996) 1-20.

\bibitem{soh} W. Y. Soh, S. A. Berger, Fully developed flow in a curved pipe of arbitrary curvature ratio, Int. J. Numer. Meth. Fluids, 7 (1987) 733-755.

\bibitem{topa} H. C. Topakoglu, M. A. Ebadian, On the steady laminar flow of an incompressible viscous fluid in a curved pipe of elliptical cross-section, J. Fluid Mech., 158 (1985) 329-340.

\bibitem{yang} Z. H. Yang, H. B. Keller, Multiple laminar flows through curved pipes, Appl. Numer. Math., 2, (1986) 257-271.

\end{thebibliography}
\end{document}